%% file: fatou.tex
\documentstyle[fullpage,12pt]{article}
\input psfig
\begin{document}

\newtheorem{lem}{Lemma}[section]
\newtheorem{prop}{Proposition}
\newtheorem{con}{Construction}[section]
\newtheorem{defi}{Definition}[section]
\newcommand{\hf}{\hat{f}}
\newtheorem{fact}{Fact}[section]
\newenvironment{proof}{{\bf Proof:}\newline}{\begin{flushright}{\em
Q.E.D.}\end{flushright}}
\newcommand{\ex}{{\cal EX}}
\newcommand{\Cr}{{\bf Cr}}
\newcommand{\dist}{\mbox{dist}}
\newcommand{\mod}{\mbox{mod }}
\newenvironment{double}{\renewcommand{\baselinestretch}{2}\protect\Large
\protect\normalsize}{}

\title{Hyperbolicity is Dense in the Real Quadratic Family}
\author{Grzegorz \'{S}wi\c{a}tek\\ Mathematics  Department\\Princeton
University\\ Princeton, NJ 08544, USA
\thanks{The author gratefully acknowledges partial support from  
NSF Grant \#DMS-9206793 and the Sloan Foundation.}}
\date{March 5, 1994}
\maketitle

\thispagestyle{empty}
\input imsmark.tex
\SBIMSMark{1992/10}{July 1992}{revised March 1994}

\begin{abstract}
It is shown that for non-hyperbolic real quadratic polynomials  
topological and quasisymmetric conjugacy classes are the same.

By quasiconformal
rigidity, each class has only one representative in the quadratic
family, which proves that hyperbolic maps are dense. 
\end{abstract}                                   

\section{Fundamental concepts}
\subsection{Introduction}
\paragraph{Statement of the results.}
{\bf Dense Hyperbolicity Theorem}\begin{em}
In the real quadratic family 
\[ f_{a}(x) = ax(1-x)\; ,\, 0<a\leq 4\]
the mapping $f_{a}$ has an attracting cycle, and thus is hyperbolic on
its Julia set, for an open and dense set of parameters $a$.
\end{em}  

\medskip
What we  actually prove is this:\newline
{\bf Main Theorem}
\begin{em}
Let $f$ and $\hat{f}$ be two real quadratic polynomials with a
bounded forward critical orbit and no attracting or indifferent cycles.  
Then, if they are topologically conjugate, the conjugacy extends to a 
quasiconformal conjugacy between
their analytic continuations to the complex plane.    
\end{em} 
\bigskip

\subparagraph{Derivation of the Dense Hyperbolicity Theorem.}
We show that the Main Theorem implies the Dense Hyperbolicity Theorem.
Quasiconformal conjugacy classes of normalized complex quadratic
polynomials are
known to be either points or open (see~\cite{mss}.) We remind the reader,
see~\cite{mithu}, that the kneading sequence is aperiodic for a real
quadratic polynomial precisely when this polynomial has no attracting
or indifferent periodic orbits. Therefore, by the Main Theorem,
topological conjugacy classes of real quadratic polynomials with
aperiodic kneading sequences are either points or open in the space of
real parameters $a$. On the other hand, it is an elementary
observation that the set of
polynomials with the same aperiodic kneading sequence in the real 
quadratic family is also closed. So, for every aperiodic kneading
sequence there is at most one polynomial in the real quadratic family
with this kneading sequence. 

Next, between two parameter values $a_{1}$ and $a_{2}$ for which
different kneading sequences occur, there is a parameter $a$ so that
$f_{a}$ has a periodic kneading sequence. So, the only way the Dense
Hyperbolicity Theorem could fail is if there were an interval filled
with polynomials without attracting periodic orbits and yet with
periodic kneading sequences. Such polynomials would all have to be
parabolic (have indifferent periodic orbits). It well-known, however,
by the work of~\cite{dohu}, that there are only countably many such
polynomials. The Dense Hyperbolicity Theorem follows.

\subparagraph{Consequences of the theorems.}
The Dense Hyperbolicity Conjecture had a long history. In a paper from
1920, see~\cite{fat}, Fatou expressed the
belief that ``general'' (generic in today's language?) rational
maps are expanding on the Julia set. Our result may be regarded as
progress in the verification of his conjecture. More recently, the
fundamental work of Milnor and Thurston, see~\cite{mithu}, showed the 
monotonicity of the kneading invariant in the
quadratic family. They also conjectured that the
set of parameter values for which attractive periodic orbits
exist is dense, which means that the kneading sequence is strictly
increasing unless it is periodic. The Dense Hyperbolicity Theorem
implies  Milnor and Thurston's conjecture. Otherwise, we would have an
interval in the parameter space filled with polynomials with an
aperiodic kneading sequence, in a clear violation of the Dense
Hyperbolicity Theorem.

\subparagraph{The Main Theorem and other results.}
Yoccoz, ~\cite{yoc}, proved that a non-hyperbolic quadratic polynomial
with a fixed non-periodic kneading sequence is unique
up to an affine conjugacy unless it is infinitely renormalizable. Thus,
we only need to prove our Main Theorem if the maps are infinitely
renormalizable. However, our approach automatically gives a proof for
all non-hyperbolic polynomials, so we provide an independent argument. 
The work of~\cite{miszczu} proved the Main Theorem for infinitely
renormalizable polynomials of bounded combinatorial type. 
The paper~\cite{kus} proved the Main Theorem for some infinitely
renormalizable quadratic polynomials not covered by~\cite{miszczu}. 
A different approach to Fatou's conjecture was taken in a recent
paper~\cite{macu}. That work proves that there is no invariant line
field on the Julia set of an infinitely renormalizable real
polynomial. This result implies that there is no non-hyperbolic
component of the Mandelbrot set containing this real polynomial in its
interior, however it is not known if it can also imply the Dense
Hyperbolicity Theorem.

\paragraph{Beginning of the proof.}
Our method is based on the direct
construction of a quasiconformal conjugacy and relies on 
techniques developed in~\cite{yours} and~\cite{kus}.

The pull-back construction shown in~\cite{miszczu} allows one to pass
from quasisymmetric conjugacy classes on the real line to
quasiconformal conjugacies in the complex plane. Thus, the Main
Theorem is reduced to conjugacies on the real line.

We can reduce the proof of the Main
Theorem, to the
following Reduced Theorem:\newline
{\bf Reduced Theorem}
\begin{em}
Let $f$ and $\hat{f}$ be two real quadratic 
polynomials with the same aperiodic kneading sequence and bounded
forward critical orbits. Normalize them to the
form $x\rightarrow ax(1-x)$. Then the conjugacy between  
$f$ and $\hat{f}$ on the interval $[0,1]$ is quasisymmetric.    
\end{em}

This paper relies heavily on~\cite{yours}, which describes the
inducing process on the real line, and~\cite{kus} which worked out
many ideas and estimates that we use. 

\paragraph{Questions remaining.} 
Does the Main Theorem remain true for quadratic S-unimodal maps? I
believe that the
difficulties here are only technical in nature and the answer should
be affirmative. However, in that case the Dense Hyperbolicity Theorem
will not follow from the Main Theorem.

Is the Main Theorem true for unimodal polynomials of higher even degree?
The present proof uses degree $2$ in the proof of Theorem 4. 
Theorem 4 seems an irreplaceable element of the proof. 

Our result implies that the polynomials for which a homoclinic
tangency occurs (the critical orbit meets a repelling periodic orbit)
are dense in the set of non-hyperbolic polynomials. On the other hand,    
in view of~\cite{invmes}, those correspond to density points in the
parameter space of the set of polynomials with an absolutely continuous
invariant measure. It seems reasonable to conjecture that the set of
non-hyperbolic maps without an absolutely continuous invariant measure
has $0$ Lebesgue measure. This problem was posed by J. Palis
during the workshop in Trieste in June 1992. We do not know of any
recent progress in solving this problem. 
 
\paragraph{Acknowledgements.}
Many ideas of the proof come from~\cite{yours} and~\cite{kus} which
were done jointly with Michael Jakobson. His
support in the preparation of the present paper was also
crucial. Another important source of my ideas were Dennis
Sullivan's lectures which I heard in New York in 1988.  I am also
grateful to Jean-Christophe Yoccoz
for pointing out to certain deficiencies of the original draft.
Jacek Graczyk helped me with many discussions as well as by giving the
idea of the proof of Proposition~\ref{prop:16gp,1}. 

\subsection{Outline of the paper.}
In order to prove the Reduced Theorem we apply the inducing
construction, essentially similar to the one used in~\cite{kus}, to
$f$ and $\hf$. We also develop the technique for constructing
quasiconformal ``branchwise equivalences'' in a parallel pull-back
construction. The infinitely renormalizable case is treated by
constructing a ``saturated map'' on each stage of renormalization,
together with a uniformly quasisymmetric branchwise equivalence,
and sewing them to get the quasisymmetric conjugacy. These are the
same ideas as used in~\cite{kus}. 
   
The rest of section 1 is devoted to defining and introducing main
concepts of the proof. We also reduce the Reduced Theorem to an even
simpler Theorem 1. Theorem 1 allows one to eliminate renormalization
from the picture and proceed in almost the same way in renormalizable
and non-renormalizable cases. 

The results of section 2 are summarized in Theorem 2. 
Theorem 2 represents the beginning stage in the construction of
induced mappings and branchwise equivalences. Our main technique of
``complex pull-back'', introduced later in section 3, may not
immediately apply to high renormalizations of a polynomial, since
those are not known to be complex polynomial-like in the sense
of~\cite{dohu}. For this reason we are forced to proceed mostly by
real methods introduced in~\cite{kus}. This section also contains an
important new lemma about nearly parabolic S-unimodal mappings, i.e.
Proposition~\ref{prop:16gp,1}.     

In section 3 we introduce our powerful tool for constructing
quasiconformal branchwise equivalences. This combines certain ideas
of~cite{kus} (internal marking) with complex pull-back similar to a
construction used in~\cite{brahu}. The main features of the
construction are  described by Theorem 3. We then proceed to prove
Theorem 4. Theorem 4 describes the conformal geometry of our so-called
``box case'', which is somewhat similar to the persistently
recurrent case studied by~\cite{yoc}. Another proof of Theorem 4 can
be found in~\cite{indue}. However, the proof we give is simpler once
we can apply our technique of complex pull-back of branchwise
equivalences.

In section 4 we apply the complex pull-back construction to the
induced objects obtained by Theorem 2. Estimates are based on Theorems
3 and 4. The results of this section are given by Theorem 5. 

Section 5 concludes the proof of Theorem 1 from Theorem 5. The
construction of saturated mappings follows the work of~\cite{kus}
quite closely. 

The Appendix contains a result related to Theorem 4 and illustrates
the technique of separating annuli on which the work of~\cite{indue},
referenced from this paper, is based. The result of the Appendix is
not, however, an integral part of the proof of our main theorems. 

To help the reader (and the author as well) to cope with the size of
the paper, we tried to make all sections, with the exception of
section 1, as independent as possible. Cross-section references are
mostly limited to the Theorems so that, hopefully, each section can be
studied independently.  

\subsection{Induced mappings}
We define a class of unimodal mappings. 

\begin{defi}\label{defi:9kp,1}
For $\eta>0$, we define the class $\cal F_{\eta}$ to comprise all
unimodal mappings of the
interval $[0,1]$ into itself normalized so that $0$ is a fixed point
which satisfy these conditions:
\begin{itemize}
\item
Any $f\in {\cal F}$ can be written as $h(x^{2})$ where
$h$ is a polynomial defined on a set containing $[0,1]$ with range 
$(-1-\eta, 1+\eta)$. 
\item
The map $h$ has no critical values except on the real line.
\item
The Schwarzian derivative of $h$ is non-positive.
\item
The mapping $f$ has no attracting or indifferent periodic cycles.
\item
The critical orbit is recurrent.
\end{itemize}
We also define 
\[ {\cal F} := \bigcup_{\eta>0} {\cal F}_{\eta}\; .\]
\end{defi}

We observe that class $\cal F$ contains all infinitely renormalizable
quadratic polynomials and their renormalizations, up to an affine
change of coordinates. 

\paragraph{Induced maps}
The method of inducing was applied to the study of unimodal maps first
in~\cite{invmes}, then in~\cite{gujo}.  
In~\cite{kus} and~\cite{yours} an elaborate approach was developed to
study induced maps, that is, transformations defined to be iterations
of the original unimodal map restricted to pieces of the domain. We
define a more general and abstract notion in this work, namely:

\begin{defi}\label{defi:704a,1}
A {\em generalized induced map $\phi$ on an interval $J$} is assumed
to satisfy the following conditions:
\begin{itemize}
\item
the domain of $\phi$, called $U$, is an open and dense subset of $J$,
\item
$\phi$ maps into $J$, 
\item
restricted to each connected component $\phi$ is a polynomial with all
critical values on the real line and with negative Schwarzian
derivative,
\item
all critical points of $\phi$ are of order $2$ and each connected
component of $U$ contains at most one critical point of $\phi$.  
\end{itemize}
\end{defi}

A restriction of a generalized induced map to a connected component of
its domain will be called a {\em branch} of $\phi$. Depending on
whether the domain of this branch contains the critical point or not,
the branch will be called {\em folding} or {\em monotone}. Domains of
branches  of $\phi$ will also be referred to as {\em domains of
$\phi$}, not to be confused with the domain of $\phi$ which is $U$.    
In most cases  generalized induced maps should be thought of as
piecewise iterations of a mapping from $\cal F$. If they do not arise
in this way, we will describe them as {\em artificial maps}.

\subparagraph{The fundamental inducing domain.}
By the assumption that all periodic orbits are repelling, every $f\in
{\cal F}$ has a fixed point $q>0$.  

\begin{defi}\label{defi:30xa,1}
If $f\in{\cal F}$, we define the {\em fundamental inducing domain} of
$f$. Consider the first return time of the critical point to the
interval $(-q,q)$. If it is not equal to $3$, or it is equal to $3$
and there is a periodic point of period $3$, then the fundamental
inducing domain is $(-q,q)$. Otherwise, there is a periodic point $q'<0$ of
period $2$ inside $(-q,q)$. Then, the fundamental inducing domain is
$(q',-q')$.  
\end{defi}

\paragraph{Branchwise equivalences.}
\begin{defi}\label{defi:25na,1}
Given two generalized induced mappings on $J$ and $\hat{J}$
respectively,  a {\em branchwise
equivalence} between them is an
orientation preserving homeomorphism of $J$ onto $\hat{J}$  which
maps the domain $U$ of the first map onto the domain $\hat{U}$ of the
second map.
\end{defi}

So the notion of a branchwise equivalence is independent of the
dynamics, only of domains of the generalized induced mappings.  

\subsection{Conjugacy between renormalizable maps}
\paragraph{The Real K\"{o}be Lemma.}
Consider a diffeomorphism $h$ onto its image $(b,c)$. Suppose that 
its has an extension $\tilde{h}$ onto a larger image $(a,d)$ which is
still a diffeomorphism. Provided that  $\tilde{h}$ has  negative Schwarzian
derivative, and $\frac{|a-b|\cdot|c-d|}{|c-a|\cdot|d-b|}\geq
\epsilon$, we will say that $h$ is 
{\em $\epsilon$-extendible.} The following holds for
$\epsilon$-extendible maps:
\begin{fact}\label{fa:9ga,1}
There is a function $C$ of $\epsilon$ only so that
$C(\epsilon)\rightarrow 0$ as $\epsilon\rightarrow 1$ and 
for every $h$ defined on an
interval $I$ and $\epsilon$-extendible,
\[ |{\cal N}h|\cdot |I| \leq C(\epsilon) \; .\]
\end{fact}
\begin{proof}
Apart from the limit behavior as $\epsilon$ goes to $1$, this fact
is proved in~\cite{mest}, Theorem IV.1.2 . The asymptotic behavior can
be obtained from Lemma 1 of~\cite{guke} which says that if $\tilde{h}$
maps the unit interval into itself, then 
\begin{equation}\label{equ:3gp,1}
{\cal N}h(x) \leq \frac{2h'(x)}{\mbox{dist}\:(\{0,1\}, h(x))} \; .
\end{equation}
The normalization condition can be satisfied by pre- and
post-composing $\tilde{h}$ with affine maps. This will not change
${\cal N}\tilde{h}\cdot |I|$, so we just assume that $\tilde{h}$ is
normalized.  
Since we are interested in $\epsilon$ close to $1$, the  denominator 
of~(\ref{equ:3gp,1}) is large and $h'(x)$ is no more than 
\[ \exp C(\frac{1}{2})\frac{|h(I)|}{|I|}\; .\]
 As $|h(I)|$ goes to $0$
with $\epsilon$ growing to $1$, we are done. 
\end{proof}

\paragraph{Properties of renormalization.}
A mapping  $f\in{\cal F}$ will be called {\em renormalizable} provided
that a {\em restrictive interval} exists for $f$. An open interval $J$
symmetric with respect to $0$   
will be called {\em restrictive} if for some $n>1$ intervals 
$J$, $f(J)$, $f^{n}(J)$ are disjoint, whereas $f^{n}(J)\subset J$. 
These definitions are broadly used in literature and can be traced
back at least to~\cite{guke1}.
Given an $f$, the notions of a locally maximal and maximal restrictive
interval will be used which are self-explanatory. Observe that if $J$
is locally maximal, then $f^{\partial J}\subset\partial J$. 

If $f$ is renormalizable, $J$ is its maximal restrictive interval and 
$n$ is the first return time form $J$ into itself, we can consider 
$f^{n}$ restricted to $J$ which will be called the {\em first
renormalization} of $f$. If $f$ is in $\cal F$, 
we define its {\em renormalization 
sequence} $f_{0}, f_{1}, \cdots, f_{\omega}$. Here $\omega$ can be
finite or infinity meaning that the sequence is infinite. The
definition is inductive.  
$f_{0}$ is $f$. If $f_{i}$ is renormalizable, then $f_{i+1}$
is the first renormalization of $f_{i}$. If $f_{i}$ is
non-renormalizable, the sequence ends. The original mapping $f$ is
called infinitely renormalizable if $\omega=\infty$, finitely
renormalizable if $0<\omega<\infty$ and non-renormalizable if
$\omega=0$.  
  
\subparagraph{Distortion in renormalization sequences.}
\begin{fact}\ref{fa:9ga,1}
Let $f\in {\cal F}_{\eta}$ and $f_{i}$ be the renormalization
sequence. For every $\eta>0$ there is a $\tilde{\eta}>0$ so that for
every $i$ $f_{i}$ belongs to ${\cal F}_{\tilde{\eta}}$ after an
affine change of coordinates.
\end{fact}
\begin{proof}
A similar estimate appeared in~\cite{miszczu}. Our version appears
as a step in the proof of Lemma VI.2.1 of~\cite{mest}. 
\end{proof} 

\paragraph{Saturated mappings.}
Let us assume that we have a topologically conjugate pair $f$ and
$\hf$, $f,\hf\in {\cal F}$. Let $f_{i}$ and $\hf_{i}$ be the
corresponding renormalization sequences. As a consequence of $f$ and
$\hf$ being conjugate, $f_{i}$ and $\hf_{i}$ are conjugate for each
$i$. Also, both renormalization sequences are of the same length. 

\begin{defi}\label{defi:29xa,2}
Let a renormalization sequence $f_{i}$ be given, and let $i<\omega$. Then, we
define the {\em saturated map} $\phi_{i}$ of $f_{i}$ as a generalized
induced  map
(Definition~\ref{defi:704a,1}) on the
fundamental inducing domain of $f_{i}$. The domain of $\phi_{i}$ is the
the backward orbit of the fundamental inducing domain $J$ of $f_{i+1}$ under
$f_{i}$. Restricted to a connected set of points whose first entry
time into $J$ is $j$, the mapping is $f^{j}$.
\end{defi}

\begin{defi}\label{29xa,3}
If $i<\omega$, a {\em saturated branchwise equivalence}
$\upsilon_{i}$ is any  branchwise equivalence between the
saturated maps. If $i=\omega$, the saturated branchwise is also
defined 
equivalence and is just the topological conjugacy on the fundamental
inducing domain of $f_{i}$. 
\end{defi}  \noindent

In this situation, we have a following fact:
\begin{fact}\label{fa:9ga,2}
Let $f$ and $\hf$ be a topologically conjugate pair of renormalizable
mappings with their renormalization sequences. Assume the existence of
a $K>0$  and a sequence of saturated branchwise equivalences $\upsilon_{i}$, 
$i\leq \omega$ which satisfy these estimates.
\begin{itemize}
\item
Every $\upsilon_{i}$ is $K$-quasisymmetric.
\item
For $i<\omega$ every domain of $\upsilon_{i}$ is adjacent to two other
domains, and for any pair of adjacent domains the ratio of their
lengths is bounded by $K$. 
\item
For every $i<\omega$ all branches of the corresponding saturated
maps are at least $1/K$-extendible.   
\end{itemize}

Then, the topological conjugacy between $f$ and $\hf$ is
quasisymmetric with a norm
bounded as a function of $K$ only.    
\end{fact}
\begin{proof}
This is a direct consequence of Theorem 2 of~\cite{kus}.
\end{proof}
\paragraph{Further reduction of the problem.}
We will prove the following theorem:\newline
{\bf Theorem 1}\newline\begin{em}
Suppose that $f$ and $\hf$ are both in ${\cal F}_{\eta}$ for some
$\eta>0$ and are topologically conjugate. Then, there is a bound $K$
depending only on $\eta$ so that there is a $K$-quasisymmetric 
saturated branchwise equivalence $\upsilon_{0}$. In addition, if
$\omega>0$, the all branches of the saturated maps $\phi_{0}$ and
$\hat{\phi}_{0}$ are $1/K$-extendible. 
\end{em}  

\subparagraph{Theorem 1 implies the Reduced Theorem.}
We check the hypotheses of Fact~\ref{fa:9ga,2}. To check that $\upsilon_{i}$ is
uniformly quasiconformal, we apply Theorem 1 to $f_{i}$ (after an
affine change of coordinates the resulting map is in $\cal F$). By
Fact~\ref{fa:9ga,1}, all these mappings belong to some ${\cal
F}_{\tilde{\eta}}$ where $\tilde{\eta}$ only depends on $\eta$.   
So, Theorem 1 implies that all saturated mappings are uniformly
quasisymmetric. In the same way we derive the uniform extendibility of
saturated maps $\phi_{i}$ and $\hat{\phi}_{i}$. It is a well-known
fact (see~\cite{guke1}) that preimages of the fundamental inducing domain of
$f_{i+1}$, or of any neighborhood of $0$, by $f_{i}$ are dense. 

We still need to check the condition regarding adjacent domains of
$\phi_{i}$ and $\hat{\phi}_{i}$. We only do the check for $\phi_{i}$,
since it is the same in the phase space of the other mapping. 
Inside the domain of $f_{i+1}$ every domain of $\phi_{i}$ is adjacent
to two others with
comparable lengths. For this, see the proof of Proposition 1
in~\cite{kus} where the preimages of the
fundamental inducing domain are explicitly constructed. The
computations done there are also applicable in our case since
the ``distortion norm'' used there is bounded in terms of
$\tilde{\eta}$ by the Real K\"{o}be Lemma. 

Denote by $W$  the maximal restrictive interval of $f_{i}$ (so $W$ is
the domain of $f_{i+1}$). Outside of $W$, consider a connected
component of point with the same first
entry time $j_{0}$ into the enlargement of $W$ 
with scale $1+\tilde{\eta}$. By this definition and the Real K\"{o}be
Lemma, the distortion of $f^{j}_{0}$ on $f^{-j_{0}}_{i}(W)$ is
bounded in terms of $\tilde{\eta}$. Now every pair of
adjacent domains of $\phi_{i}$ can be obtained as the image under some
$f^{-j_{0}}$ in this form of a pair of adjacent domains from the
interior of the restrictive interval of $f_{i}$. It follows that the
condition of Fact~\ref{fa:9ga,2} regarding the ratio of lengths of
adjacent domains is satisfied for any pair. Now, the Reduced Theorem
follows from Theorem 1.      

\subparagraph{Theorem 1 for non-renormalizable mappings.}
Theorem 1 also has content for non-renormalizable mappings. It follows
from~\cite{yoc} if $f$ and $\hf$ are in the quadratic family, or
renormalizations of quadratic polynomials.  Our Theorem 1 is
marginally more general.  

\vskip 1in
\goodbreak
\vskip -1in

\subsection{Box mappings} 
\paragraph{Real box mappings.}
\begin{defi}\label{defi:1,1}
A generalized induced map $\phi$ on an interval $J$ symmetric with
respect to $0$ is called a {\em box mapping},
or a {\em real box mapping}, if 
its branches satisfy these conditions. For every branch consider the
smallest interval $W$ symmetric with respect to $0$ which contains the
range of this branch. We will say that $W$ is a {\em box} of $\phi$,
and we will also say that the branch {\em ranges through $W$}. 
For every branch of $\phi$, it
is assumed that it ranges through an interval which is contained in
$J$, and that the  endpoints of this interval 
are {\em not} in the domain $U$ of $\phi$. In addition, it is assumed
that all branches are monotone except maybe the one whose domain
contains $0$ and that there is a monotone branch mapping onto $J$.  
Also, the central domain of $\phi$ is always considered a box of
$\phi$.   
\end{defi}

Every box mapping has a {\em box structure} which is simply the
collection of all its boxes ordered by inclusion. 

\paragraph{Complex box mappings.}
\begin{defi}\label{defi:10ga,1}
Let $\phi$ be a real box mapping. We will say that $\Phi$ is a {\em
complex box mapping} and a {\em complex extension} of $\phi$ provided
that the following holds:
\begin{em}
\begin{itemize}
\item
$\Phi$ is defined on an open set $V$ symmetric with respect to the real
axis. We assume that  connected components of $V$ are topological
disks, call them {\em domains} of $\Phi$, and refer the restriction of
$\Phi$ to any of its domains as a {\em branch} of $\Phi$, or simply a
{\em complex branch}.   
\item
If a domain of $\Phi$ has a non-empty intersection $W$ with the real line, then
$W$ must be a domain of $\phi$ in the sense of Definition~\ref{defi:1,1}. 
Moreover, the complex branch defined on this domain is the analytic
continuation of the branch of $\phi$ defined on $W$, and this analytic
continuation has no critical points in the closure of the domain 
except on the real line.
\item
Every domain of $\phi$ which belongs to a branch that {\em does not}
map onto the entire $J$ is the intersection of the real
line with a domain of $\Phi$. If this is true for those domains
belonging to long monotone branches as well, we talk of a {\em complex
box mapping with diamonds.}
\item
If two domains of $\Phi$ are analytic continuations of branches of $\phi$
that range through the same box of $\phi$, then the corresponding two
branches of $\Phi$ have the same image.
\item
If a domain of $\Phi$ is disjoint with the real line, then the branch
defined there is univalent and shares its range with some branch of
$\Phi$ whose domain intersects the real line.
\item
The boundary of the range of any branch of $\Phi$ is disjoint with
$V$.
\end{itemize}
\end{em}
\end{defi}

We see that a complex box mapping also has a box structure defined as
the set of the ranges of all its branches. These {\em complex boxes}
are in a one-to-one correspondence with the boxes of the real mapping
$\phi$. Also, we will freely talk of monotone or folding branches for
complex mappings, rather than univalent or degree 2.

\paragraph{Special types of box mappings.}
The definitions below have the same wording for real and complex box
mappings. So we just talk of a ``box mapping'' with specifying whether
it is real or complex. 

\begin{defi}\label{defi:9mp,1}
A box mapping $\phi$ is said to be of type I if it satisfies these
conditions. 
\begin{em}
\begin{itemize}
\item
The box structure
contains three boxes, which are denoted 
$B_{0}\supset B'\supset B$.
\item
$B$ equals the domain of the central branch (i.e. the only domain of
$\phi$ which contains $0$.)
\item
Every monotone branch maps onto $B$ or $B_{0}$. Depending on which
possibility occurs, we talk of {\em long} and {\em short} monotone
branches. 
\item
Every short monotone branch can be continued analytically to a
diffeomorphism onto $B'$ and the domain of such a continuation is
either compactly contained in $B'$, or compactly contained in the
complement of $B'$. Restricted to the real line, this continuation has
negative Schwarzian derivative.  
\item
The central branch is folding and ranges through $B'$. It also has an
analytic continuation of degree $2$ which maps onto a larger set
$B''\subset B_{0}$
which compactly contains $B'$. The domain of this continuation is
compactly contained in $B'$. The Schwarzian derivative of this
continuation restricted to the real line is negative.  
\item
The closure of the union of domains of all short monotone branches is
disjoint with the boundary of $B'$.
\item
If the closures of two domains domains intersect, at least one of these domains
is long monotone.
\end{itemize}
\end{em}
\end{defi}

A real box mapping of type I is shown on Figure~\ref{fig:4xp,1}. 
A complex box mapping of type I is shown on Figure~\ref{fig:5xa,1}.

\begin{defi}\label{defi:9mp,2}
A box mapping $\phi$ is said to be of type  II if the following is
satisfied.
\begin{em}
\begin{itemize}
\item
The box structure contains three boxes, $B_{0}\subset B'\supset B$. 
$B$ denotes the domain of the central branch. 
\item
All monotone branches of $\phi$ map either onto $B_{0}$ or $B'$, and
are accordingly classified as long or short. 
\item
The central branch is folding, and it has an analytic continuation 
whose domain is compactly contained in $B'$, the range is some
$B''\subset B_{0}$ which compactly contains $B'$, the degree of this
continuation is $2$, and the Schwarzian derivative of its restriction
to the real line is negative. 
\item
The closure of the union of all domains of short monotone branches is
disjoint with the boundary of $B'$.
\item
If the closures of two domains domains intersect, at least one of these domains
is long monotone.
\end{itemize}
\end{em}
\end{defi}

\begin{defi}\label{defi:11ga,1}
A box mapping will be called {\em full} if it has only two boxes,
$B_{0}\supset B$. Also, the central branch must be folding and range
through $B_{0}$. 
\end{defi}

\paragraph{Hole structures.}
Given a complex box mapping $\Phi$, the domains which intersect the real line
and such that the range of $\Phi$ restricted to any of them is less
than the largest box will be called {\em holes}. The central domain is
also considered a hole, and described as the central hole. Other
domains of $\Phi$ which intersects the real line will be called {\em
diamonds.}  Boxes, holes and diamonds can then be studied purely
geometrically, without any reference to the dynamics. So we consider
{\em hole structures} which are collections of boxes and holes, and
{\em hole structures with diamonds} which also include diamonds. 
Given a hole structure, the number of its boxes minus one will be
called the {\em rank} of the structure.  
 
\subparagraph{Some geometry of hole structures.}
Although the facts which we will state now have an easy generalization
for hole structures of any rank, we formally restrict ourselves to
structures of rank $0$.
We want to specify what bounded geometry is for a hole structure. A
positive number $K$ is said to provide the bound for a hole structure
if the series of estimates listed below are satisfied. Here we adopt
the notation $I$ to mean the box of the structure. 
\begin{enumerate}
\item
All holes are strictly contained in $I$, moreover, between any hole
and $I$ there is an annulus of modulus at least $K^{-1}$.
\item
All holes and $I$ are at least $K$-quasidisks. The same estimate holds        
for half-holes and half-$I$, that is, regions in which a half of the
quasidisks was cut off along the real line.
\item
Consider a hole and its intersection $(a,b)$ with the real line. Then,
there are points $a' \leq a$ and $b' \geq b$ with the property that
the hole is enclosed in the diamond of angle $\pi/2-K^{-1}$ based on
$a',b'$. Furthermore, for different holes the corresponding intervals
$(a',b')$ are disjoint.
\end{enumerate}

\subparagraph{The ``mouth lemma''.}
\begin{defi}\label{defi:15fa,2}
A {\em diamond neighborhood} $D(\theta)$ of an interval $J$ given by an angle
$0<\theta<\pi$ is the union of two regions symmetric with respect to
the real line and defined as follows. Each region is bounded by a
circular arc which intersects the real line in the endpoints of $J$.
We adopt the convention that $\theta$ close to $0$ given a very thin
diamond, while $\theta$ close to $\pi$ gives something called a
``butterfly'' in~\cite{miszczu}.    
\end{defi}

We know consider an enriched hole structure of rank $0$. 
\begin{defi}\label{defi:15fp,1}
Given a rank $0$ hole structure, 
in the part of $J$ which is not covered by the holes we arbitrarily
choose a set of disjoint open intervals. For each interval, we define
a diamond neighborhood. Each diamond has a neighborhood of 
modulus $K_{h}^{-1}$ which is contained in the box.  Next, each
diamond neighborhood
contains two symmetrical arcs. It is assumed that each of them joins
the endpoints in the interval and except for them is disjoint with the
line, and that each arc is $K_{h}$-quasiarc.
This object will be referred to
as a {\em hole structure with diamonds}.
\end{defi}
A bound for a hole structure with diamonds is the greater of the norm
for the hole structure without diamonds and $K_{h}$. 

For a hole structure with diamonds, consider the ``teeth''. i.e. the
curve which consists of
the upper halves of the boundaries of the holes and the diamonds, and
of points of the real line which are outside of all holes and
diamonds. 
Close this curve with the ``lip'' which means the upper half
of the boundary of $I$. The whole curve can justly be called
the ``mouth''.
\begin{lem}\label{lem:11na,1}
For a bounded hole structure with diamonds, the mouth is a quasidisk
with the norm uniform with
respect to the bound of the map.
\end{lem}
\begin{proof}
The proof is based of the three point property of Ahlfors (see~\cite{ahl}.)  
A Jordan curve is said to have this property if for any pair of its
points their Euclidean distance is comparable with the
smaller of the diameters of the arcs joining them. Moreover, a uniform
bound in the three point property implies an estimate on the
distortion of the quasicircle. Our conditions of boundedness were set 
precisely to make
it easy to verify the three point property. One has to consider
various choices of the two points. If they are both on the same tooth,
or both on the lip, the property follows immediately from the fact
that teeth are uniform quasicircles. An interesting situation is when
one point is on one tooth and the other one on another tooth. We
notice that it is enough to check the diameter of a simpler arc which
goes along either tooth to the real line, than takes a shortcut along the
line, and climbs the other tooth (left to the reader.)
Suppose
first that both teeth are boundaries of holes. Consider the one which
is on the left. We first consider the arc which goes to the $b'$. By   
conditions 2 and 3 of the boundedness of hole structures, the
diameter of this arc is not only comparable to the distance from the
point to $b'$, but even to the distance from the projection of the
point to the line. The estimate follows. The case when one of the
teeth is a diamond is left to the reader. Another interesting case is
when one point is on the lip and another one on a tooth. If the tooth
is a diamond, the estimate follows right away from the choice of
angles. If the tooth is a hole, we have to use condition 1. We use
the bounded modulus to construct a uniformly bounded quasiconformal map
straightens the lip and the tooth to round circles without changing
the modulus too much. The estimate also follows.
\end{proof}

\subparagraph{Extension lemma.}
We are ready to state our main result which in the future will be
instrumental in extending real branchwise equivalences to the complex domain.
Given two hole structures with diamonds, a homeomorphism $h$ is said
to establish the equivalence between them if it is order preserving,
and for any hole, box, or diamond of one structure, it transforms its
intersection with the real axis onto the intersection between a hole,
box or diamond respectively of the other structure with the real axis.
\begin{lem}\label{lem:11na,2}
For a uniform constant $K$, 
let two $K$-bounded hole structures with diamonds be with an
equivalence establishing $K$-quasisymmetric homeomorphism $h$. We want
to extend $h$ to
the complex plane
with prescribed behavior on the boundary of each tooth and on the lip.
This prescribed behavior is restricted by the following condition. For
each tooth or the lip consider the closed curve whose base is the
interval on the real line, and the rest is the boundary of the tooth
or the lip, respectively. The map is already predefined on this curve.
We demand that it maps onto the corresponding object of the other hole
structure, and that it extends to a $K$-quasiconformal
homeomorphism. Then, $h$ has a global
quasiconformal extension with the prescribed behavior and whose
quasiconformal norm is uniformly bounded in terms of $K$.
\end{lem}
\begin{proof}
We first fill the teeth with a uniformly quasiconformal map which
is something we assumed was possible. Next, we extend to the lower
half-plane by~\cite{abu}. Then, we want an extension to the region
above the lip and above the real line. This is certainly possible,
since the boundary of this region is a uniform quasicircle (because
the lip was such, use the three point property to check the rest).
But the map can be defined below this curve by filling the lip and
using~\cite{abu} in the lower half-plane. Then one uses reflection
(see~\cite{ahl}.) Now, the map is already defined outside of the
mouth, one uses Lemma~\ref{lem:11na,1} and quasiconformal reflection
again to finish the proof.    
\end{proof}
\subsection{Introduction to inducing}\label{sec:13gp,1}
\paragraph{Standard extendibility of real box mappings.}
For monotone branches we have the notion of extendibility which is
compatible with our statement of the Real K\"{o}be Lemma
(Fact~\ref{fa:9ga,1}.) That is, a monotone branch with the range equal
to $(b,c)$ will be deemed {\em $\epsilon$-extendible} if it has an
extension as a diffeomorphism with negative Schwarzian derivative onto
a larger interval $(a,d)$ and so that 
\[ \frac{|a-b|\cdot |d-c|}{|c-a|\cdot |d-b|} \geq \epsilon\; .\]
This amounts to saying that the length of $(b,c)$ in the
Poincar\'{e} metric of $(a,d)$ is no more than $-\log\epsilon$. 
\footnote{See~\cite{mest} for a discussion of
the Poincar\'{e} metric on the interval.}  The interval $(a,d)$ will
be called the {\em margin of extendibility} while the domain of the
extension will be referred to as the {\em collar of extendibility}. 

We can also talk of extendibility for the central folding branch. To
this end, we represent the folding branch as $h(x^{2})$. Suppose that
the central branch ranges through a box $(b,c)$. Then the folding
branch is considered $\epsilon$-extendible provided that $h$ has an
extension as a diffeomorphism with negative Schwarzian derivative onto
a larger interval $(a,d)$ and the Poincar\'{e} length of $(b,c)$
inside $(a,d)$ is at most $-\log\epsilon$. Again, $(a,d)$ is described
as the margin of extendibility and its preimage by the extension of
the folding branch is the collar of extendibility.  

\begin{defi}\label{defi:10gp,1}
Consider a box mapping $\phi$ on $J$. {\em Standard
$\epsilon$-extendibility} for $\phi$ means that the following
properties hold.   
\begin{em}
\begin{itemize}
\item
There is an interval $I\supset J$ so that the Poincar\'{e} length of
$J$ in $I$ is no more than $-\log\epsilon$. Every branch of $\phi$
which ranges through $J$ is extendible with the margin of
extendibility equal to $I$. Furthermore, for each such branch the
collar of extendibility also has Poincar\'{e} length not exceeding
$\epsilon$ in $I$; moreover if the domain of the branch is compactly
contained in $J$, so is the the collar of extendibility.
\item
For every other box $B$ in the box structure of $\phi$, there is an
interval $I_{B}\supset \overline{B}$ which serves as the margin of
extendibility for all branches that range through $B$. $I_{B}$ must be
contained in any box larger than $B$. 
If the domain
of a branch is contained in $B$, the collar of extendibility of this
branch must be contained in $I_{B}$, moreover, if this domain is
compactly contained in $B$, the collar is contained in $B$.  
\end{itemize}
\end{em}
\end{defi}
We will call a choice of $I$ and $I_{B}$ the {\em extendibility
structure} of $\phi$.    
\paragraph{Refinement at the boundary.}
Consider a generalized induced map $\phi$ on $J$. If we look at an
endpoint of $J$, clearly one of two possibilities
occurs. Either the domains of $\phi$ accumulate at the endpoint, or
there is one branch whose domain has the endpoint on its boundary. In the
first case, we will say that $\phi$  is {\em infinitely refined} at its
(left, right) boundary. In the second case the branch adjacent to the
the endpoint will be called the (left, right) {\em external branch}. 

Assume now that an external branch ranges through $J$, and that its
common endpoint with $J$ is repelling fixed point of the branch.  
Call the external branch $\zeta$. We will describe the construction of
{\em refinement at a boundary}.
Namely, we can define $\phi^{1}$ as $\phi$ outside of the
domain of $\zeta$, and as $\phi\circ\zeta$ on the domain of $\zeta$. 
Inductively, $\phi^{i+1}$ can be defined as $\phi$ outside of the
domain of $\zeta$, and $\phi^{i}\circ\zeta$ on the domain of $\zeta$. 
The endpoint is repelling. Thus, the external branches of $\phi^{i}$
adjacent to the
boundary  point will shrink at an exponential rate. We can either stop
at some moment and get a version of $\phi$ {\em finitely refined at
the boundary} or proceed to an 
$L^{\infty}$ limit of this construction to obtain a mapping infinitely
refined at this boundary. There is an analogous process for the other
external branch, $\zeta'$ (it is exists). Namely, define $(\phi^{i})'$
as $\phi$ outside the domain of $\zeta'$ and $\phi^{i-1}\circ\zeta'$
on the domain of $\zeta'$.  

We call the mapping obtained in this
process a {\em version of $\phi$ refined at the boundary} (finitely or
infinitely.) 
Observe that the refinement at the boundary does not change either the
box structure of the map or its extendibility structure.  

\paragraph{Operations of inducing.}
Loosely speaking, inducing means that some branches of  a box mapping are 
being replaced
 by compositions with other branches. Below we
give precise definitions of five procedures of inducing. 

\subparagraph{Simultaneous monotone pull-back.}
Suppose that a box mapping $\phi$ on $J$ is given with a set of
branches of $\phi$ which range through $J$. We denote
$\phi'=\phi_{b}$ where $\phi_{b}$ is a version of $\phi$ refined at
the boundary. Then on each branch
$\zeta$ from this set replace $\phi$ with $\phi^{t}_{b}\zeta$ where
$\phi^{t}_{b}$ is $\phi_{b}$ with the central branch branch replaced by the
identity.   Observe that this operation does not change the box or 
extendibility structures. This is a tautology except for the one
which is just the restriction of
$\zeta$ to the preimage of the central domain. This is extendible with
margin $I$, so the standard extendibility can be satisfied with the
same $\epsilon$.    

\subparagraph{Filling-in.}
To perform this operation, we need a box mapping $\phi$ and $\phi'$
with the same box and extendibility structures as $\phi$ and with a choice
of boxes $B_{1}$ smaller than $B_{0}:=J$ but bigger than the central
domain $B$.  
We denote $\phi_{0}:=\phi'$ and proceed inductively by considering all
monotone branches of $\phi$ which map onto $B_{1}$. On the domain
of each such branch $\zeta$, we replace $\zeta$ with $\phi^{t}_{i}
\zeta$. Here, $\phi^{t}_{i}$ denotes $\phi_{i}$ whose central branch was
replaced by the identity. 
The resulting map is $\phi_{i+1}$. In this way we proceed to the
$L^{\infty}$ limit, and this limit $\phi_{\infty}$ is the outcome of the
filling-in. Note that
$B_{1}$ drops from the box structure of $\phi_{\infty}$, but otherwise the
new map has the the same boxes with the same margins of extendibility
as $\phi$. Thus, if $\phi$ satisfied standard
$\epsilon$-extendibility, so does $\phi_{\infty}$. 

\subparagraph{Simple critical pull-back.}
For this, we need a box mapping $\phi$ whose central branch is folding,
and another box mapping $\phi'$ about which we assume that it has the same
box structure as $\phi$ and is not refined at either boundary.  
Let $\phi^{t}$ denote $\phi'$ whose central
branch was replaced by the identity, and $\psi$ be the central branch
of $\phi$. Then the outcome is equal to $\phi$ modified on the central
domain by substituting $\psi$ with $\phi^{t}\circ\psi$. 

\subparagraph{Almost parabolic critical pull-back.}
The arguments of this operation are the same as for the simple
critical pull-back. In addition, we assume that the critical value of
$\psi$ is in the central domain of $\psi$ (and $\psi'$), but $0$ is
not in the range of $\psi$ from the real line. Also, it is assumed
that $\psi^{l}(0)\notin B$ so some $l>0$, and let $l$ denote the
smallest integer with this property. Then, for each point $x\in B$ we
define the {\em exit time} $e(x)$ as the smallest $j$ so that
$\psi^{j}(x)\notin B$. Clearly, $e(x)\leq l$ for every $x\in B$. Then
the outcome of this procedure is $\psi$ modified on the subset of
points with exit times less than $l$ by replacing
$\psi$ with 
\[ x\rightarrow \psi'\circ\psi^{e(x)}(x) \; .\]
On the set of point with exit time equal to $l$, which form a
neighborhood of $0$, the map is unchanged. 

\subparagraph{Critical pull-back with filling-in.}
This operation takes a box mapping $\phi$ which must be of type I or
full, and another mapping $\phi'$ with the same box structure and not
refined at either boundary. 
Let $\psi$ denote the central branch of $\phi$, and let $\chi$ be the
generic notation for a short monotone branch of $\phi$. We define
inductively two sequences, $\phi_{i}$ and $\phi^{r}_{i}$. Let
$\phi_{0}:=\phi$ and $\phi_{0}^{r}$ be $\phi'$ with the central branch
replaced by the identity. Then $\phi_{i+1}$ is obtained from
$\phi$ by replacing the central branch $\psi$ with $\phi_{i}^{r}\psi$
and every short monotone branch $\chi$ of $\phi$ with
$\phi^{r}_{i}\circ\psi\circ\chi$. Then $\phi_{i+1}^{r}$ is obtained in
a similar way, but only those short monotone and folding branches of
$\phi$ are replaced whose domains are contained in the range of
$\psi$. Others are left unchanged. At the end, the central branch of
the map obtained in this way is replaced by the identity. The outcome
is the $L^{\infty}$ limit of the sequence $\phi_{i}$. This is the most
complicated procedure and is described in~\cite{kus}, and called
filling-in, in the case when $\phi$ is full. We note here that
$\phi_{\infty}$ is either of type I, full or Markov depending on the
position of the critical value of $\psi$. If the critical value was in
a long monotone domain of $\phi'$, then the outcome is a full mapping.
If it was in a short monotone domain of $\phi'$ or in the central
domain, the result is a type I map. If the critical value was not in the
domain of $\phi'$, the result is a Markov map, that is a mapping whose
all branches are monotone and onto $J$. 

\paragraph{Boundary refinement.} 
The five basic operations described above obviously fall in two
classes. The first two involve composing with monotone branches and do
not affect standard extendibility. The last three involve composing
with folding branches and can affect standard extendibility. Each of
these last three operations can be preceded by a process called {\em
boundary refinement.} To do the boundary refinement, assume that a box
mapping $\phi$ is given whose central branch $\psi$ is folding. Then look at
the set of ``bad'' branches of $\phi$. For the steps of simple
critical pull-back and critical pull-back with filling-in, a branch of
$\phi$ is
considered ``bad'' if it ranges through $J$, its domain is in the
range of $\psi$ and its collar of standard extendibility contains the
critical value. For almost parabolic critical pull-back the last
condition is weakened so that the branches which contain
$\psi^{i}(0)$, $0<i\leq l$ are also considered bad. Then boundary
refinement is defined as the simultaneous monotone pull-back on the
set of bad branches  with $\phi'$ that may vary with the branch. 
There are two possibilities here. If we do {\em infinite} boundary
refinement, as $\phi'$ we use the version of $\phi$ infinitely refined
at both endpoints. If we do {\em minimal} boundary refinement, we use
the version of $\phi$ refined at the endpoint of the side of the
critical value only, and refined to the minimal depth which will
enforce standard extendibility.

This gives the mapping $\phi'$ which then enters
the corresponding procedure of critical pull-back. 
\begin{lem}\label{lem:12gp,1}
Let $\phi$ be a type I or full mapping with standard
$\epsilon$-extendibility. Apply critical pull-back with filling-in
preceded by boundary refinement. Then, the resulting map
$\phi_{\infty}$ also has standard $\epsilon$-extendibility.
\end{lem}
\begin{proof}
Boundary refinement assures us that for every $i$ the 
branches of $\phi_{i}$ or $\phi^{r}_{i}$ which range through $J$ are 
extendible with the same
margin as in the original $\phi$. This ends the proof if
$\phi_{\infty}$ if full. Otherwise, let $B$ denote the
central domain of
$\phi_{\infty}$ (the same as the central domain of $\phi_{i}$, any
$i>0$), and $B'$ the central domain of $\phi$. Observe that all short
branches of $\phi_{\infty}$ are extendible with margin $B'$. Thus, one
can set $I_{B'}$ inherited from the extendibility structure of $\phi$,
and $I_{B}:=B'$. One easily checks that the conditions of standard   
extendibility hold. 
\end{proof}

\vskip 1in
\goodbreak
\vskip -1in

\section{Initial inducing}
\subsection{Real branchwise equivalences}
\paragraph{Statement of the result.}
Our definitions follow~\cite{kus}. 
\begin{defi}\label{defi:7ya,1}
Let $\upsilon_{1}$ be a quasisymmetric and order-preserving homeomorphism
of the line onto itself. Let $\upsilon_{2}$ be a quasisymmetric
order-preserving homeomorphism from an interval $J$ onto another
interval, say $J'$.
We will say that $\upsilon_{2}$ {\em replaces} $\upsilon_{1}$ on an interval
$(a,b)$ with distortion $L$ if the following mapping $\upsilon$ is an
$L$-quasisymmetric order preserving homeomorphism of the line onto
itself:
\begin{itemize}
\item
outside of $(a,b)$, the mapping $\upsilon$ is the same as
$\upsilon_{1}$,
\item
inside $(a,b)$, $\upsilon$ has the form
\[ \upsilon= A'\circ \upsilon_{2} \circ A^{-1}\]
where $A$ and $A'$ are affine and map $J$  onto $[a,b]$ and $J'$ onto
$\upsilon_{1}([a,b])$ respectively.
\end{itemize}
\end{defi}

\begin{defi}\label{defi:13gp,1}
We will say that a branchwise equivalence $\upsilon$ between 
box mappings on $J$ satisfies the {standard
replacement condition} with distortion $K$ provided that
\begin{itemize}
\item
$\upsilon$ restricted to any domain of the box mapping
replaces $\upsilon$ on $J$,
\item
$\upsilon$ restricted to any box of this mapping replaces $\upsilon$
on $J$.
\end{itemize}
\end{defi}

\begin{defi}\label{defi:13yp,1}
A box map $\phi$ is said to be 
{\em $\alpha$-fine}
if for every domain $D$ and box $B$ so that $D\subset B$ but
$\overline{D}\not\subset B$ 
\[ \frac{|D|}{\dist(D,\partial B)} \leq \alpha \; .\]
\end{defi}

\begin{defi}\label{defi:14gp,1}
A branch is called {\em external in the box $B$} provided that the
domain $D$ of the branch is contained in $B$, but the closure of $D$
is not contained in $D$. We will say that the map {\em is not} refined 
at the boundary of $B$ provided that an external branch exists at each
endpoint of $B$.
\end{defi}
 
\begin{prop}\label{prop:13gp,1}
Suppose that box mappings $\phi$ and $\phi'$
are given which will enter one of the five inducing operations
defined in section~\ref{sec:13gp,1}. 
Assume further that $\hat{\phi}$ is topologically conjugate to $\phi$, 
and $\hat{\phi'}$ is topologically conjugate to $\phi'$, with the
same topological conjugacy. Let $\upsilon_{0}$ and $\upsilon'_{0}$ be the
branchwise equivalences between $\phi$ and $\hat{\phi}$ as well as
$\phi'$ and $\hat{\phi}'$ respectively. Suppose that all these branchwise
equivalences coincide outside of $J$.   
Assume in addition that 
\begin{em}
\begin{itemize}
\item
 all branches of $\phi$, $\phi'$, $\hat{\phi}$, $\hat{\phi}'$ are
$\epsilon$-extendible, 
\item
all four box mappings are $\alpha$-fine, 
\item
$\upsilon_{0}$ and $\upsilon'_{0}$ satisfy
the standard replacement condition with distortion $K$ and are
$Q$-quasisymmetric,
\item
$\phi'$ and $\hat{\phi'}$ are not refined at the boundary of any box
smaller than $B_{0}$, 
\item
if $D$ is the domain of an external branch of in a box $B$ (for
$\phi'$ or $\hat{\phi'}$), then  $|D|/|B|\geq \epsilon_{1}$.   
\end{itemize}
\end{em}
Then for all $\epsilon$, $K$, $Q$, $\alpha$ and $\epsilon_{1}$ there are
bounds $L_{1}$ and $L_{2}$ so that the following holds. Perform one of
the five inducing operations on $\phi$ and $\phi'$ as well as
$\hat{\phi}$ and $\hat{\phi'}$. Call the outcomes $\phi_{\infty}$ and
$\hat{\phi}_{\infty}$, or $\phi_{\infty,b}$ and
$\hat{\phi}_{\infty,b}$ for their versions refined at the boundary.
Then there are branchwise equivalences $\upsilon_{\infty}$ and
$\upsilon_{\infty,b}$ between $\phi_{\infty}$ and
$\hat{\phi}_{\infty}$, or $\phi_{\infty,b}$ and
$\hat{\phi}_{\infty,b}$ respectively which satisfy
\begin{itemize}
\item
$\upsilon_{\infty}$ and $\upsilon_{\infty,b}$ are the same as
$\upsilon$ outside of $J$,
\item
$\upsilon_{\infty}$ and $\upsilon_{\infty,b}$ satisfy the standard
replacement condition with distortion $L_{1}$,
\item
$\upsilon_{\infty}$ and $\upsilon_{\infty,b}$ are both
$L_{2}$-quasisymmetric.
\end{itemize}
\end{prop}

\paragraph{Some technical material.}
The proof of Proposition~\ref{prop:13gp,1} splits naturally in five
cases depending on the kind of inducing operation involved. Among
these the almost parabolic critical pull-back stands out as the
hardest. Proposition~\ref{prop:13gp,1} in all remaining cases follows rather
straightforwardly from the work of~\cite{kus}. We begin by recalling
the technical tools of~\cite{kus}. 

\begin{defi}\label{defi:13gp,2}
Given an interval $I$ on the real line, its {\em diamond neighborhood}
is a set symmetrical with respect to the real axis constructed as
follows. In the upper half-plane, the diamond neighborhood is the
intersection of a round disk with the upper half-plane, while on the
intersection of the same disk with the real line is $I$. For a diamond
neighborhood, its {\em height} refers to twice the Hausdorff distance
between the neighborhood and $I$ divided by the length of $I$.
\end{defi}

This is slightly different from the definition used in~\cite{kus}, but
this is not an essential difference in any argument. The nice property
of diamond neighborhoods is the way they are pull-back by polynomial
diffeomorphisms. 

\begin{fact}\label{fa:13gp,1}
Let $h$ be a polynomial which is a diffeomorphism from an interval $I$
onto $J$. 
Assume that $h$ preserves the real line and that all its critical
values are on the real line. Suppose also that $h$ is still a
diffeomorphism from a larger interval $I'\supset I$ onto a larger
interval $J'\supset J$ so that the Poincar\'{e} length of $J$ inside
$J'$ does not exceed $-\log\epsilon$. 
Then, there are constants $K_{1}>0$ and $K_{2}$ depending on
$\epsilon$ only so that the distortion of $h$ (measured as
$|h''/(h')^{2}|$) on the diamond neighborhood of $I$ with height $K_{1}$ is
bounded by a $K_{2}$. 
\end{fact}
\begin{proof}
Observe that the inverse branch of $h$ which maps $J'$ onto $I'$ is
well defined in the entire slit plane ${\bf C}\setminus J'$. Since $h$
is a local diffeomorphism on $I$, the inverse branch can be defined on
a diamond neighborhood of $J'$ with sufficiently small height. As we
gradually increase the height of this neighborhood we see that the
inverse branch can be defined until the boundary of the neighborhood
hits a critical value of $h$. That will never occur, though, since all
critical values are on the real line by assumption. Now, the diamond
neighborhood of $J$ with height $1$ is contained in ${\bf C}\setminus
J'$ with modulus bounded away from $0$ in terms of $\epsilon$. So, by
K\"{o}be's distortion lemma, the preimage of this diamond neighborhood
by the inverse branch of $h$ contains a diamond neighborhood of  $I$
of definite height, and the distortion of $h$ there is bounded.   
\end{proof}
\paragraph{Pull-back of branchwise equivalences.}
\begin{defi}\label{defi:8fa,1}
A {\em pull-back ensemble} is the conglomerate of the following objects:
\begin{em}
\begin{itemize}
\item
two equivalent induced mappings, one in the phase space of $f$, the
other in the phase space of $\hat{f}$, with a branchwise equivalence
$\cal Y$ and a pair of distinguished branches $\Delta$ and
$\hat{\Delta}$ whose domains, $D$ and $\hat{D}$ respectively,
correspond by $\cal Y$,
\item
a branchwise equivalence $\Upsilon$ which must map the critical value
of $\Delta$ to the critical value of $\hat{\Delta}$ in case if
$\Delta$ is folding.
\end{itemize}
\end{em}

The following objects and notations will always be associated with a
pull-back ensemble.      
If $\Delta$ is monotone, let $(-Q,Q')$, $\frac{6}{5}Q>Q'\geq Q$, be
its range. If $\Delta$ is folding, make $Q'=Q$ and define $(-Q,Q')$ to
be the smallest interval symmetrical with respect to $0$
which contains the range of $\Delta$. Analogously, we define
$(\hat{Q},\hat{Q}')$ to be the range of $\hat{\Delta}$ or the maximal
interval symmetric with respect
to $0$ containing the range of $\hat{\Delta}$, and we also assume that 
$\frac{6}{5}\hat{Q}> \hat{Q}' \geq \hat{Q}$.  

In addition, the following assumptions are a part of the definition:
\begin{em}
\begin{itemize}
\item
$\Upsilon((-Q,Q')) = (-\hat{Q},\hat{Q}')$, 
\item
$\Upsilon(\overline{z}) = \overline{\Upsilon(z)}$ and
$\Upsilon(-\overline{z}) = -\overline{\Upsilon(z)}$ for every $z\in {\bf
C}$,
\item
outside of the disc 
$B(0,\frac{4}{3}Q)$ the map $\Upsilon$ has the form
\[ z\rightarrow \frac{\hat{Q}'-\hat{Q}}{Q'-Q}z\; ,\]  
\item
the mapping $\Upsilon$ restricted to the set
\[{\bf R}\setminus (-Q,Q')\]
 has a global $\lambda$-quasiconformal extension,
\item
the branches $\Delta$ and $\hat{\Delta}$ are
$\epsilon$-extendible,
\item
if $\Delta$ is folding, then its range is smaller than
$(-Q,Q')$, but the length of the range  is at least
$\epsilon_{1}Q$; likewise, if $\hat{\Delta}$ is folding, then its range
does not fill  the interval $(-\hat{Q}, \hat{Q}')$, but the length of
the range
is at least $\epsilon_{1}\hat{Q}$.
\end{itemize}
\end{em}
\end{defi}

The numbers $\epsilon$, $\epsilon_{1}$ and $\lambda$ will be
called {\em parameters} of the pull-back ensemble.  

\subparagraph{Vertical squeezing.}
\begin{defi}\label{defi:12fp,1}
Suppose that two parameters $s_{1}\geq 2s_{2}>0$ are given. Consider a
differentiable monotonic function $v: {\bf R}\rightarrow {\bf R}$ with
the following properties:
\begin{itemize}
\item
$v(-x) = -x$ for every $x$, 
\item
if $0<x<s_{2}$, then $v(x)=x$,
\item
if $s_{1}<x$, then $v(x) = x-s_{1}+2s_{2}$.
\end{itemize}

Each $v$ defines a homeomorphism of the plane $\cal V$ defined by
\[ {\cal V}(x+iy)=x+iv(y) \; .\]
We want to think of $\cal V$ as being quasiconformal end depending on
$s_{1}$ and $s_{2}$ only. To this end, for a given $s_{1}$ and $s_{2}$
pick some $v$ which minimizes the maximal conformal distortion of
$\cal V$. We will call this $\cal V$ the {\em vertical squeezing map}
for parameters $s_{1}$ and $s_{2}$. 
\end{defi}   
\subparagraph{The Sewing Lemma.}
\begin{fact}\label{fa:14ga,1}
Consider a pull-back ensemble with $Q=Q'$.  Let  $\cal Y$ be a
branchwise equivalence with $D$ as a domain, and suppose that
restricted to $D$ it 
 replaces $\Upsilon$ on $(-Q,Q)$ with distortion $M$. Suppose that 
for some   $R>0$ the map $\Upsilon$
transforms  the diamond neighborhood with
height $R$ of $(-Q,Q)$ exactly on the diamond neighborhood with height 
$R$ of $(-\hat{Q}, \hat{Q})$. Choose $r>0$ and $C\leq 1$ arbitrary
and assume that the
diamond neighborhood with height $r$ of $D$ is mapped by $\cal Y$
onto the diamond neighborhood with height $C r$ of $\hat{D}$
Assuming that $R\leq R_{0}$ where $R_{0}$ only depends on the
parameter $\epsilon$ of the pull-back ensemble, for
every such choice of $C,M,r$, $R$, and a set of parameters of the
pull-back ensemble, numbers $K$ and $L$, parameters $s_{1}, s_{2},
\hat{s}_{1}, \hat{s}_{2}$, a mapping $\tilde{\Upsilon}$ as well as a 
branchwise equivalence $\tilde{\cal Y}$ exist so that:
\begin{itemize}
\item
if $\cal V$ is the vertical squeezing with parameters $s_{1}, s_{2}$,
and $\hat{\cal V}$ is the vertical squeezing map with parameters
$\hat{s}_{1}$ and $\hat{s}_{2}$, then $\tilde{\cal Y}$ has the form
\[ \tilde{\cal Y} = \hat{\cal V}\circ \tilde{\Upsilon}\circ {\cal
V}^{-1}\; \] on the image of the domain  of $\tilde{\Upsilon}$ by
$\cal V$, 
\item
the domain and range of $\tilde{\Upsilon}$ are contained in diamond
neighborhoods with height $1/10$ of $D$ and $\hat{D}$ respectively,
\item
on the diamond neighborhood with height $K$ of $D$
\[ \tilde{\Upsilon} = \Delta^{-1} \circ \Upsilon \circ\Delta\]
which means the lift to branched covers, order-preserving on the real 
line in the case of $\Delta$ folding,   
\item
outside of this diamond neighborhood the conformal distortion of
$\tilde{\Upsilon}$ is bounded by the sum of $L$ and the conformal
distortion of $\Upsilon$ outside of the image of the diamond neighborhood of
$D$ with height $K$ by $\Delta$,   
\item
$\tilde{\cal Y}$ coincides with $\cal Y$ outside of the diamond
neighborhood with height
$r$ of $D$, 
\item
on the set-theoretical difference between the diamond neighborhoods
with heights $r$ and $r/2$, the map $\tilde{\cal Y}$ is quasiconformal
and its conformal distortion is bounded as the sum of $L$ and the
conformal distortion of $\cal Y$ on the diamond neighborhood of $D$
with height $r$,
\item
outside of the diamond neighborhood with height $r/2$ of $D$, the
mapping $\tilde{\cal Y}$ is independent of $\Upsilon$, 
\item
on the set-theoretical difference between the diamond neighborhood of
$D$ with height $r/2$ and $U$ 
the map $\tilde{\cal Y}$ has conformal distortion bounded almost
everywhere by $L$.
\end{itemize}
\end{fact}
\begin{proof}
This fact is a consequence of Lemmas 4.3 and 4.4 of~\cite{kus}. The
context of~\cite{kus} is slightly different from our situation, since
that paper works in the category of S-unimodal mappings and uses their
``tangent extensions'' instead of analytic continuation. However,
proofs given in~\cite{kus} work in our situation with only semantic
modifications, so we do not repeat them here. 
\end{proof} 

As a consequence of the Sewing Lemma, we get
\begin{fact}\label{fa:15gp,1}
Suppose that a pull-back ensemble is given so that $\Upsilon$
transforms the diamond neighborhood with height $R$ of $(-Q,Q)$
exactly to the homothetic diamond neighborhood of
$(-\hat{Q},\hat{Q})$. Assume also that $R\leq R_{0}$ as required by
the hypothesis of the Sewing Lemma, and that $\cal Y$ restricted to
$D$ replaces $\Upsilon$ on $(-Q,Q)$ with distortion $M$. Suppose
finally that both $\Upsilon$ and $\cal Y$ are $Q$-quasiconformal. Then
there is a map $\tilde{\cal Y}$ which differs from $\cal Y$ only on
the diamond neighborhood of $D$ with height $1/2$, is equal to
\[ \hat{\Delta}^{-1}\circ \Upsilon\circ\Delta\]
on $D$, and is $L$-quasiconformal. The number $L$ only depends on $M$,
$Q$ and the parameters of the pull-back ensemble.  
\end{fact}
\begin{proof}
This follows directly from Fact~\ref{fa:14ga,1} when one chooses
$r=1/2$ and observes that $C$ is bounded as a function of $Q$. 
\end{proof}

\paragraph{Beginning the proof of Proposition~\ref{prop:13gp,1}.}
In order to be able to use Fact~\ref{fa:14ga,1} we need a lemma that
will allow us to build pull-back ensembles. 
\begin{lem}\label{lem:14ga,1}
Suppose that a branchwise equivalence $\upsilon$ between topologically
conjugate box mappings $\phi$ and $\hat{\phi}$ is given on $(-Q,Q)$. 
Assume that this branchwise equivalence is
$Q$-quasisymmetric and satisfies the standard replacement
condition with distortion $K$. Suppose also that a pair of branches
$\Delta$ and $\hat{\Delta}$ are given so that
$\Delta$ and $\hat{\Delta}$ range through boxes that correspond by
$\upsilon$. Moreover, assume that if $\Delta$ one is folding, the
other one is, too, and the critical values correspond by the
topological conjugacy. Suppose also that $\phi$ and $\hat{\phi}$ are
both $\alpha$-fine and all of their branches are
$\epsilon$-extendible. 
If $\Delta$ and $\hat{\Delta}$ are also $\epsilon$
extendible, then there  there is a mapping $\upsilon'$ which coincides
with $\upsilon$ on $(-\frac{6}{5}Q, \frac{6}{5}Q)$ except on the
domain which contains the critical value
of $\Delta$, and an $L$-quasiconformal extension of
$\upsilon'$, called $\Upsilon$, which together with $\Delta$ and
$\hat{\Delta}$ gives a pull-back ensemble with parameters $\epsilon$,
$\epsilon_{1}$ and $\lambda$. Numbers $L$ and $\lambda$ depend on $Q$
and $K$ only. The mapping $\Upsilon$ transforms the diamond
neighborhood with height $R_{0}$ of $(-Q,Q)$ exactly onto the diamond
neighborhood of $(-\hat{Q}, \hat{Q})$ with the same height, where
$R_{0}$ is determined in terms of $\epsilon$ by Fact~\ref{fa:14ga,1}.
\end{lem} 
\begin{proof}
We have two tasks to perform: one is to build $\upsilon'$ so that it
maps the critical value of $\Delta$ to the critical value of
$\hat{\Delta}$, and second is to construct the proper quasiconformal
extension. Let us address the second problem first. We use this fact.
\begin{fact}\label{fa:15gp,2}
Let $\upsilon$ be a $Q$-quasisymmetric mapping of the line into
itself, and let $(-\hat{Q},\hat{Q})=\upsilon((-Q,Q))$. For every
$R<1$ there an $L$-quasiconformal homeomorphism $\Upsilon$ of the plane with
these properties: 
\begin{em}
\begin{itemize}
\item
$\Upsilon$ maps the diamond neighborhood of $(-Q,Q)$ with height $R$
exactly onto the diamond neighborhood of $(-\hat{Q},\hat{Q})$ with the
same height, 
\item
$\Upsilon$ restricted to $(-\frac{6}{5}Q, \frac{6}{5}Q)$ equals
$\upsilon$,
\item
$\Upsilon$ maps $B(0,\frac{4}{3}Q)$ exactly onto
$B(0,\frac{4}{3}\hat{Q})$ and is affine outside of this ball.
\end{itemize}
\end{em}
The bound $L$ depends on $R$ (continuously) and $Q$ only.
\end{fact}

Fact~\ref{fa:15gp,2} follows from the construction 
done in~\cite{kus} in the proof of Lemma 4.6. From
Fact~\ref{fa:15gp,2} we see that once $\upsilon'$ is constructed with
the desired properties, Lemma~\ref{lem:14ga,1} follows.

To get $\upsilon'$, we first construct take any point $c$ and
construct $\upsilon''$ which maps
$c$ to $H(c)$, where $H$ is the topological conjugacy, and
$\upsilon''=\upsilon$ outside of the box that contains $c$. To
construct $\upsilon''$, we consider two cases. If the
$c$ is not in an external domain, one can compose
$\upsilon''$ with a diffeomorphism of bounded distortion which moves 
$\upsilon(c)$ to $H(c)$ and is the identity outside
of the box. If $c$ is an external domain, map it by the
external branch into the phase space of a version of $\phi'$ refined
at the boundary. By the previous argument, this image requires only a
push by a diffeomorphism of bounded distortion. This way we get some
branchwise equivalence $\upsilon_{1}$. Using Fact~\ref{fa:15gp,2}, we
build the pull-back ensemble which consists of the external branch, its
counterpart in the phase space of $\hat{\phi'}$, and the appropriate
extension $\Upsilon_{1}$ of $\upsilon_{1}$. To this pull-back ensemble
we can apply Fact~\ref{fa:15gp,1} and get $\upsilon''$ on the real
line. By the way, in this case $\upsilon''$ can immediately be
taken as $\upsilon'$. Generally, in order to obtain $\upsilon'$ from 
$\upsilon''$ we apply
a similar procedure. Take the branch which contains the critical
value of $\Delta$ and use the image of the critical value by this
branch  as the
$c$ to build $\upsilon''$. Then construct the pull-back ensemble by
Fact~\ref{fa:15gp,2} and
pull $\upsilon''$ back to get $\upsilon'$. Again, $\upsilon'$ is
quasisymmetric from Fact~\ref{fa:15gp,1}
\end{proof}  

We will next prove that Lemma~\ref{lem:14ga,1} and the Sewing Lemma
allow us to construct $\upsilon_{\infty}$ and $\upsilon_{\infty,b}$
which are quasisymmetric as needed. The construction of
$\upsilon_{\infty,b}$ is quite the same as for $\upsilon_{\infty}$,
only we start with $\phi_{b}$. So we only focus on the construction of
$\upsilon_{\infty,b}$.  

\paragraph{Bounded cases of inducing.}
Suppose that we are in the situation of Proposition~\ref{prop:13gp,1}. 
In the case of a simple critical pull-back and simultaneous monotone
pull-back, the quasisymmetric estimate for $\upsilon_{\infty}$ follows
directly. By Lemma~\ref{lem:14ga,1} for each branch being refined we
construct a pull-back ensemble, and get $\tilde{\cal Y}$ which is
quasiconformal with an appropriate bound by Fact~\ref{fa:15gp,1}. Note
that in the case of the simultaneous pull-back the operations on
various branches do not interfere, since each modifies $\cal Y$ only in
the diamond neighborhood with height $1/2$. 

\subparagraph{The case of filling-in.}
We begin by constructing the extension $\Upsilon$ of $\upsilon'_{0}$
using Lemma~\ref{lem:14ga,1}. For this, choose $R$ equal to $R_{0}$
which is given by the Sewing Lemma depending on $\epsilon$, and
$(-Q,Q)$ equal to the box $B_{1}$. Denote $\Upsilon_{0}=\Upsilon$. 
Also, take as $\cal Y$ any quasiconformal extension of $\upsilon_{0}$
with a norm bounded as a function of $Q$. We will proceed inductively
and obtain $\Upsilon_{i+1}$ by applying the Sewing Lemma always with
the same $\cal Y$, $\Delta$ ranging over the set of all branches which
map onto $B_{1}$, and $\Upsilon_{i}$.  

Let $\delta$ be the generic notation for the domains of branches of
$\phi$ mapping onto $B_{1}$, and let $\delta^{-m}$ denote similar
domains of $\phi^{t}_{m}$.
Let us choose $r$ so small that the diamond neighborhoods of domains
$\delta$ of $\phi$ are contained in the diamond neighborhood with
height $R$ of $(-Q,Q)$. In our inductive construction, this will
assure that $\Upsilon_{i}$ for any $i$ is the same as $\Upsilon$ on
the boundary of the diamond neighborhood of $(-Q,Q)$ with height $R$,
so that this diamond is always mapped on the homothetic diamond.     
In particular,  parameters $\epsilon$, $\lambda$ and $C$ remain fixed
since only the $\Upsilon$ component changes in pull-back ensembles.  
Note that $C$ is determined by the conformal distortion of ${\cal Y}$,
i.e. by the parameter $Q$.  
Also, since $\Upsilon_{i-1}$ coincides with $\Upsilon$ on the real
line outside of $(-Q,Q)$ the condition that $\Upsilon$ restricted to a
domain replaces $\Upsilon_{i-1}$ on $(-Q,Q)$ is satisfied with the
same distortion $M$. Thus, for  all $i$ the pull-back ensemble used to
construct $\Upsilon_{i}$ satisfies the
hypotheses of the Sewing Lemma with the same
parameters. So we will regard estimates claimed by these Lemmas as
constants. 

Next, look at the neighborhood of $D$ where $\Upsilon_{i}$ has the
form 
\[ \Upsilon_{i} = \hat{\cal
V}\circ\hat{\Delta}^{-1}\circ\Upsilon_{i-1}\circ \Delta\circ{\cal
V}^{-1}\; .\]  By Fact~\ref{fa:13gp,1}, this contains a diamond
neighborhood with height $K$ which is mapped by $\Delta$ with bounded
distortion.  In particular, by choosing $r$ possibly even
smaller, but still controlled by $K$, we can make sure that the
diamond neighborhoods with height $r$ of all domains $\delta^{-0}$  
are inside the image of this region (called the {\em inner  pull-back
region}) by $\Delta$. Now, address the issue of where
$\Upsilon_{i-1}$ and $\Upsilon_{i}$ differ. If $i=1$, they differ
only on the union of diamond neighborhoods with height $r$ of domains 
$\delta^{-0}$. For $i>1$, $\Upsilon_{i}$ and
$\Upsilon_{i-1}$ differ on the preimage of the set where
$\Upsilon_{i-2}$ and $\Upsilon_{i-1}$ differ by maps 
\[ \Delta\circ{\cal V}^{-1}\]
with $\Delta$ ranging over the set of all branches mapping onto $B_{1}$.

Next, we see inductively that the set on which $\Upsilon_{i}$ and
$\Upsilon_{i-1}$ differ is contained in the union of diamond
neighborhoods with fixed height of domains $\delta^{-i+1}$ of
$\phi_{i-1}$. 
This is clearly so for $i=1$. In the general case, we need to consider
the preimage of the set where $\Upsilon_{i-1}$ differs from
$\Upsilon_{i-2}$ by $\Delta\circ{\cal V}^{-1}$. 
Suppose that $\Upsilon_{i-2}$ and
$\Upsilon_{i-1}$ differ only on the union of diamond neighborhoods
with height $r_{i-2}$ of domains $\delta^{-i+2}$. Pick some 
domain $\delta^{-i+2}$ and  observe
that $\Delta$ extends as a diffeomorphism onto $B_{1}$. It is easily
seen that sizes of domains $\delta^{-m}$ decrease with $m$ at a
uniform exponential rate. 
Therefore that $\Delta$ restricted to the preimage of $\delta$ is 
$x$-extendible with $|\log x|$ going up exponentially fast with $i$. 
So, $\Delta$ restricted to $\Delta^{-1}(\delta)$ is almost affine with
distortion going down exponentially fast with $i$ (which follows from
K\"{o}be's distortion lemma.)
Next, provided that $r_{i-2}<1$, the height
of the diamond neighborhood of $\delta$ with height $r_{i-2}$ with
respect to $\delta_{0}$ goes down exponentially with $i$. Since the
distortion goes down exponentially fast, the preimage of this diamond
neighborhood by $\Delta$ with fit inside the diamond neighborhood of 
$\delta^{-i+1} = \Delta^{-1}(\delta^{-i+2}$ with height
$r_{i-2}(1+b(i))$ where $b(i)$ decreases exponentially fast with $i$.  
Thus, by choosing $r$ small
enough, we can ensure that $r_{i}<1$ for all $i$. The same reasoning
can be conducted for the phase space of $\hat{f}$ to prove that the
images of these diamonds by $\Upsilon_{i}$ are contained in the
diamond neighborhoods with height $1$ of the corresponding 
domains $\delta^{-i+1}$ of  $\hat{\phi}$.  

For any branch $\Delta$ of $\varphi$, consider the region $W$
defined as the intersection of the inner pull-back region with the set
on which $\cal V$ is the identity (i.e. the horizontal strip of width
$2s_{2}$.) This contains a diamond neighborhood with fixed height of
the domain of $\Delta$. So, for $i>i_{0}$ ($i_{0}$ depends on how fast
the sizes of $\delta^{-m}$ decrease with $m$ and can be bounded
through $\epsilon$ and $\alpha$), the set on which
$\Upsilon_{i-2}$ and $\Upsilon_{i-1}$ differ is contained in the
image of $W$ by $\Delta$, and its image by $\Upsilon_{i-1}$ is
contained in the image of $\hat{W}$ by $\\hat{\Delta}$. 
By Fact~\ref{fa:15gp,1} applied $i_{0}$ times, $\Upsilon_{i_{0}}$ is
still quasiconformal
because each step adds only a bounded amount of distortion. For $i\geq
i_{0}$, $\Upsilon_{i}$ is obtained on this region by replacing 
$\Upsilon_{i-1}$ with 
\[ \hat{\Delta}^{-1}\circ\Upsilon_{i-1}\circ\Delta\; .\] 
This means that for $i\geq i_{0}$ the conformal distortion of
$\Upsilon_{i}$ is the same as the conformal distortion of
$\Upsilon_{i_{0}}$.              
Thus, the
limit of $\Upsilon_{\infty}$ exists and is quasiconformal with the
same bounded norm. 
\subparagraph{Critical pull-back with filling-in.}
This is very similar to the filling-in case and once we have built the
initial pull-back ensemble, the proof goes like in Lemma 4.5
of~\cite{kus}. We leave the details out. 

\paragraph{Almost parabolic critical pull-back.}
This case requires a different set of tools. Let us first set up the
core problem is abstract terms. 
\begin{defi}\label{defi:16gp,1}
We define an {\em almost parabolic map} $\varphi$ be the following
properties. 
\begin{em}
\begin{itemize}
\item
The map $\varphi$ is defined on and interval $[0,a)$ with $a<1$.
\item
$\varphi(x) = g(x^{2})$ where $g$ is an orientation-preserving
diffeomorphism onto the image of
$\varphi$, has negative Schwarzian, and is $\epsilon$-extendible, $\epsilon>0$.
\item
$\varphi(a)=1$ and $\varphi(x)>x$ for every $x$.
\item
On the set $(c,a)$, $S\varphi \geq -K$.
\end{itemize}
\end{em}

The numbers $a,\epsilon, K$ will referred to as {\em parameters} of
$\varphi$.
\end{defi}

\begin{prop}\label{prop:16gp,1}
Suppose that two almost parabolic maps $\varphi$ and $\varphi'$ are given. 
Suppose that the parameters $a$ are both bounded from below by
$alpha>0$ and from above by $\beta<1$, likewise the parameter
$\epsilon$ is 
equal for both
mappings, and the parameters $K$ are both bounded from above by
$K_{0}$. Consider the finite sequence $a_{i}$ defined by $a_{0}=a$ and
$\varphi(a_{i})=a_{i-1}$ for $i>0$ or the sequence ends
when such $a_{i}$ cannot be found.    Define $a'_{i}$ in the
same manner using $\varphi'$ and assume that the lengths of both
sequences are equal to the same $l$. Define a homeomorphism $u$ from 
$(a_{l},a)$ onto
$(a'_{l},a')$ by the requirement that $u(a_{i})=a'_{i}$ for every $i$
and that $u$ is affine on each segment $(a_{i},a_{i-1})$. 
Then $u$ is a $Q$-quasi-isometry with $Q$ only depending on $\alpha$, $\beta$,
$K_{0}$ and $\epsilon$. 
\end{prop}
\subparagraph{A comment about Proposition~\ref{prop:16gp,1}.}
What really matters is that $u$ is quasisymmetric, which follows from
its being a quasi-isometry. Proposition~\ref{prop:16gp,1} is easily
accepted by specialists in the field, perhaps because it is an easy
fact when $\varphi$ is known to be complex polynomial-like. However, I
was unable to find a fair reference in literature regarding the
negative Schwarzian setting, though I am aware of an unpublished work
of J.-C. Yoccoz in which a similar problem was encountered and solved
in the study of critical circle mappings of unbounded type. The
approach we use here owes to the work~\cite{harm}.

\paragraph{Easy bounds.}
We now assume that $\varphi$ is a mapping that satisfies the
hypotheses of Proposition~\ref{prop:16gp,1}. Throughout the proof we
will refer to bounds
that only depend on $\alpha$, $\beta$, $\epsilon$ and $K_{0}$ as 
{\em constants.}    
And so we observe that the derivative of $g$ is bounded from both
sides by positive constants. This means that the derivative of
$\varphi$ is similarly bounded from above. Next, $c$ is bounded from
below by a positive constant, or it would be impossible to maintain
$\varphi(x)>x$. 
The second derivative of $\varphi$ is also bounded
in absolute value by a constant, from the Real K\"{o}be's Lemma. 
Next, there is exactly one point $z$ where $\varphi(x)-x$ attains a
minimum. That is since because of the negative Schwarzian there are at most two
points where the derivative of $\varphi$ is one. Also observe that is
$l$ is large enough, then the distances from $z$ to $a_{i}$ and
$a_{0}$ are bounded from below by constants. That is because the
shortest of intervals $(a_{i}, a_{i-1})$ is the one containing $z$, or
the adjacent one. If
$l$ is large this becomes much smaller that the constants bounding
from below the lengths of $(a_{j}, a_{j-1})$ or $(a_{l-j}, a_{l-j-1})$
for $j<10$. A somewhat deeper fact is this.

\begin{lem}\label{lem:16gp,1}
There are positive constants $l_{0}$, $K_{1}$ and $K_{2}$ so that if
$\varphi(x)-x<K_{1}$ and $l\geq l_{0}$, then $\varphi''(x)/\varphi'(x)
\geq K_{2}$.
\end{lem}
\begin{proof}
The proof almost copies the argument used in the demonstration of
Proposition 2 in~\cite{harm}. We use ${\cal N}g$ to mean $g''/g'$. 
Take into account the differential equation
\begin{equation}\label{ident}
 D{\cal N}g = Sg + 1/2 ({\cal N}g)^{2}\; .
\end{equation}
which is satisfied by every $C^{3}$ diffeomorphism $g$ (a direct check,
also see~\cite{mest} page 56.) Then consider an abstract class of 
diffeomorphisms
${\cal G}(w,L)$
defined as the set of functions defined in a neighborhood of $0$ 
and having the
following properties:
\begin{enumerate}
\item
their Schwarzian derivatives are negative and bounded away from $0$ by some
$-\beta$.
\item
for any $g\in {\cal G}$, $g(0)=w$,     
\item
there is no $x\in (-L,L)$ where $g$ is defined and  $g(x) \leq x$ 
\end{enumerate}
Observe that the
function 
\[ g_{0}(y)= \varphi(y+x)\]
belongs to ${\cal G}(w,L)$ with 
\[ w:=\varphi(x)-x\]
and $L$ a constant equal to the minimum of distances from $x$ to $c$
and from $x$ to $a$. This is a constant provided $l_{0}$ is large
enough and $K_{1}$ is sufficiently small. So, it suffices
to show that for every $L>0$ there is a constant $K_{1}>0$ so that
$w<K_{1}$  implies ${\cal N} g(0) > K_{2} > 0$.   

We observe that every function from $g\in {\cal G}(w,L)$ is uniquely 
determined by
three parameters: a continuous function $\psi=Sg$ and two numbers
$\nu$ and $\mu$ equal to ${\cal N} g(0)$ and $g'(0)$ respectively. Indeed,
given $\psi$ and $\nu$, ${\cal N} g$ is uniquely determined by the
differential equation~(\ref{ident}), this together with $\mu$
determines $g'$, and finally $g$ is also defined by $w$. Observe that
with $\mu$ fixed, $g$ is an increasing function of $\psi$ and $\nu$.
Indeed, a look at the equation~(\ref{ident}) reveals that if
$\psi_{1}\geq \psi_{2}$ with the same $\nu$, then the solution ${\cal N}
g_{1}(x)\geq {\cal N} g_{2}(x)$ for $x\geq 0$ while ${\cal N} g_{1}\leq {\cal N} g_{2}$
for $x\leq 0$. This is immediate if $\psi_{1}>\psi_{2}$ since we see
that at every point where the solutions cross ${\cal N} g_{1}$ is bigger on
a right neighborhood and less on a left neighborhood. Then we treat
$\psi_{1}\geq \psi_{2}$ by studying $\psi' = \psi_{1}+c$ where $c$ is
a positive parameter and using continuous dependence on parameters. As $g'$ is
clearly an increasing function of ${\cal N} g$ and $\nu$, the monotonicity
with respect to $\psi$ and $\nu$ follows. So, if we can show that for
some $\tilde{\psi}$ and $\tilde{\nu}$ and every $\mu$ the condition 
$g(x) \geq x$ is violated on $(-L, L)$, it follows that there exists 
$\delta>0$ so that for every $\psi\leq
\tilde{\psi}$, we must have $\nu > \tilde{\nu}+\epsilon$ if the function is in
${\cal G}(w,L)$. 

Pick $\tilde{\psi} =-\beta$ and $\nu=0$. The problem becomes quite
explicit. From another well-known differential formula
$u'' = Sg\cdot u$ satisfied by $u=1/sqrt{g'}$ we find 
\[ g'(x) = \frac{\mu}{\sqrt{\cosh \beta x}}; .\] 
Let $w_{n}\rightarrow 0$ and pick $\mu_{n}$ so that the corresponding
$g$ satisfies $g(x)\geq x$, or escapes to $+\infty$, on $(-L,L)$.
Observe that $\mu_{n}$ must be a bounded  sequence, since we have 
$g(x)\leq \frac{\mu}{C(L,\beta)}x + w$ for $x<0$ where $C(L,\beta)$ is
the upper bound of $\cosh \beta x$ on $[-L,0]$. Thus if such a
sequence existed, we could take a limit parameter $\mu_{\infty}$ which
would preserve $g(x) \geq x$ even for $w=0$, and this cannot be. 
\end{proof}   

\paragraph{Approximation by a flow.}
For $l$ large, the condition $\varphi(x)-x<K_{1}$ of
Lemma~\ref{lem:16gp,1} is satisfied on a neighborhood $W$ of $z$. On this
neighborhood, we can bound
\[ B(x-z)^{2} + \delta \leq \varphi(x)-x \leq A(x-z)^{2} + \delta\]  
where $\delta = \varphi(z) - z$. The numbers $A$ and $B$ are
constants, and while the upper bound is easy and satisfied in the
entire domain of $\varphi$, the lower one follows from
Lemma~\ref{lem:16gp,1} is guaranteed to hold only on $W$. Note that
the number of points $a_{i}$ outside of $W$ is bounded by a constant.
So we can change the definition of $W$ a bit so that it is preserved
by the map $u$ between two almost parabolic maps. That is, we can make
points $a_{m}$ and $a_{l-m}$ the endpoints of $W$, and if $m$ is big
enough, the condition $\varphi'(x)-x$ will also hold on $u(W)$.
Incidentally, for $l$ large, this means that $z'\in W'$. The map $u$
is clearly a quasi-isometry outside of $W$, so we have reduced the
problem to considering it on $W$. Let $T=l-2m$, that is the number of
iterations an orbit needs to travel through $W$, and let $T_{0}$ be
chosen so that $[a_{T_{0}}, a_{T_{0}+1})$ contains $z$. If analogous
times are considered for $\varphi'$, note that $T'=T$, but $T_{0}$ and
$T'_{0}$ have no reason to be the same.  

\begin{lem}\label{lem:20gp,1}
Let $x_{1}, x_{2}\in W$ with $\varphi^{t}(x_{1})=x_{2}$. Then, if 
$z-x_{1} < K$, $K$ constant
\[ \sqrt{\frac{1}{4A\delta}}(\tan^{-1}
(\sqrt{\frac{A}{\delta}}(x_{2}-z)) -      
\tan^{-1}(\sqrt{\frac{A}{\delta}}(x_{1}-z))) - 1 \leq t\]
\[  \leq
\sqrt{\frac{4}{B\delta}}(\tan^{-1}(\sqrt{\frac{B}{\delta}}(x_{2}-z)) -      
\tan^{-1}(\sqrt{\frac{B}{\delta}}(x_{1}-z))) + 1\; .\]
\end{lem}
\begin{proof}
To get the lower estimate, compare $\varphi$ with the time one map of
the flow 
\[ \frac{dx}{dt} = 2A(x-z)^{2} + \delta\; .\]
We claim that iterations of the time one map of this flow overtake
iterations of the map 
$x\rightarrow x+A(x-z)^{2}$ provided that $\delta$ is small enough,
i.e. if $l$ is sufficiently large. Let $x\in (x_{1},x_{2})$. We want
\[ 2A(x+A(x-z)^{2}-z)^{2}  \geq A(x-z)^{2}\; .\]
Observe that there is nothing to prove if $x>z$. 
Choose $K$ so that $KA<1/2$ to conclude the proof. 

Then the lower estimate follows by integrating the flow. 
The upper estimate
can be demonstrated in analogous way.  
\end{proof}

We will assume that the assumption $z-x_{1}<K$ is always satisfied,
perhaps by making $W$ even smaller.  
As a corollary to Lemma~\ref{lem:20gp,1}, we note that for $l$ large, 
$\delta\sim
T^{-2}$ where the $\sim$ sign means that the ratio of the two quantities
is bounded from both sides by positive constants. Indeed, taking
$x_{1}-a_{m}$ and $x_{2}-a_{l-m}$ we see that the arguments of the
arctan terms are very large for large $l$, meaning that they values
are close to $\pi/2$, or $-\pi/2$ respectively. In particular, 
$\delta\sim\delta'$ for two different almost parabolic maps.
 
\begin{lem}\label{lem:21gp,1}
For every $b>0$ there are constants $l_{0}$ and $\beta>0$ with the
property that  if $|a_{p} - z| \leq b\sqrt{\delta}$, then    
\[ \frac{p-m}{T} > \beta \mbox{ and } \frac{l-m-p}{T}\beta \; .\]
\end{lem}
\begin{proof}
Note that $z$ is separated from the boundary of $W$ by a distance
which is at least a positive constant. This is because the ``step'' of
the orbit of $a$ is still quite large at the boundary of $W$. When
$\delta$ is sufficiently small (depending on $b$) this means that the
entire neighborhood $(z-b\sqrt{\delta}, z+b\sqrt{\delta})$ is in a
positive distance from the boundary of $W$. This means that for
$x_{1}=a_{m}$ and $x_{2}=a_{p}$ the time of passage is
$\sim\sqrt{1/\delta}$. The same reasoning applies to the passage from
$a_{p}$ to the other end of $W$. The claim follows from
Lemma~\ref{lem:20gp,1}. 
\end{proof}      

The approximation formula given by Lemma~\ref{lem:20gp,1} works well
for points which are within a distance of the order of $\sqrt{\delta}$
from $0$. We need a different formula for points further away. 
\begin{lem}\label{lem:21gp,2}
Choose a point $a_{p}$ so that the distance from $a_{p}$ to $z$ is
bigger than $b\sqrt{\delta}$. Let $t$ denote $p-m$ if $a_{p}<z$ or
$l-m-p$  otherwise.
There are constants $l_{0}$, $b_{0}$, $C>0$ and $c>0$ so that if $l\geq
l_{0}$ and $b\geq b_{0}$, then 
\[ \frac{c}{|z-x|} \leq t \leq \frac{C}{|z-x|} \; .\]
\end{lem}
\begin{proof}
We begin by comparing the iterations of $\varphi$ with flows like in
the proof of Lemma~\ref{lem:20gp,1}.  Given a flow 
\[ \frac{dx}{dt} = \delta + A(x-z)^{2}\] we compare it with the flow 
\[ \frac{dx}{dt} = A(x-z)^{2}\; .\]
The discrepancy between time $t<T$ maps of these flows is less than 
$T\delta \sim \sqrt{\delta}$. So, for $b_{0}$ large, approximating by
the simpler flow modifies 
$|z-x|$ only by a bounded factor. Integration gives the desired
estimate. 
\end{proof}

Observe that Lemma~\ref{lem:21gp,2} implies a ``converse'' of
Lemma~\ref{lem:21gp,1}. Namely,
\begin{lem}\label{lem:21gp,3}
Under the hypotheses and using notations of Lemma~\ref{lem:21gp,2}, if
$t/T > \beta > 0$, then there is a number $q$ depending on $\beta$ so
that \[ |a_{p}-z| \leq q\sqrt{\delta}\; .\]
\end{lem}
\begin{proof}
From Lemma~\ref{lem:21gp,2}, we get 
\[ x\sim t^{-1}\sim \beta^{-1}T^{-1} \sim \beta^{-1}\sqrt{\delta} \;
.\]
\end{proof}

\paragraph{Proof of Proposition~\ref{prop:16gp,1}.}
Let us conclude the proof. If $|a_{p}-z|\leq b\sqrt{\delta}$, then by
Lemmas~\ref{lem:21gp,3} and~\ref{lem:21gp,1} we get $|a'_{p}-z'|\leq
b'\sqrt{\delta}$ where $b'$ depends on $b$. Then on the interval
$(a_{p}, a_{p+1})$ the function $u$ indeed has a slope bounded from
both sides by positive numbers depending on $b$. This follows from the
approximation given by Lemma~\ref{lem:20gp,1} used with
$x_{1}:=a_{p}$ and $x_{2}:=a_{p+1}$ or $x_{1}:=a'_{p}$ and
$x_{2}:=a'_{p+1}$ respectively. Since the arguments of the arctan
functions are bounded in all cases, the estimate follows. Because of
the symmetry of the problem, one can also bound the slope of $u$ on 
$(a_{p}, a_{p+1})$ assuming that $|a'_{p}-z'|\leq b'\sqrt{\delta'}$. 
The case which remains is when both $|a_{p}-z|$ and $|a'_{p}-z'|$ are
large compared with $\delta$. But then Lemma~\ref{lem:21gp,2} is applicable 
showing that $|a_{p}-z|\sim |a'_{p}-z'|$ and since 
\[ |a_{p}-a_{p+1}|\sim (a_{p}-z)^{2}\sim (a'_{p}-z')^{2}\sim
|a'_{p}-a'_{p+1}| \]
the slope is bounded as needed. 

We finally remark that all was done
assuming $l$ sufficiently large. In the case of bounded $l$
Proposition~\ref{prop:16gp,1} is easy. So, the proof has been
finished. 
\paragraph{Construction of the branchwise equivalence.}
Faced with a case of almost parabolic critical pull-back, we first
``mark'' the branchwise equivalence so that upon its first exit from
the central domain the critical value of $\phi$ is mapped onto the
critical value of $\hat{\phi}$. This will add only bounded conformal
distortion by Lemma~\ref{lem:14ga,1}. Next, we build a quasisymmetric
map $u$ which transforms $B$ (the box through which the central branch
ranges) onto $\hat{B}$, and the backward orbit of the endpoint $a$  of
the central domain by the central branch onto the corresponding
orbit in the phase space of $\hat{\phi}$. The hard part of the proof
is that this $u$ is uniformly quasisymmetric, and that follows from
Proposition~\ref{prop:16gp,1}. Then we change $u$ between the boundary
of $B$ and the central domain by $\upsilon_{0}$, and between $a_{i}$
and $a_{i+1}$, $i\leq l-2$ by 
\[ \hat{\phi}^{-i-1}\circ\upsilon_{0}\circ\phi^{i+1}\; .\]
By Lemma 3.14 of~\cite{kus}, this will give a uniformly quasisymmetric
map on $J\setminus (a_{l-2}, -a_{l-2})$ provided that it is
quasisymmetric on any interval $(a_{i-1},a_{i+1})$. This is clear
since we are pulling-back by bounded diffeomorphisms, with the
exception of the last interval $(a_{l-2},a_{l-2})$ where this is a
quadratic map composed with a diffeomorphism. Still, it remains
quasisymmetric. 

\paragraph{Conclusion of the proof of Proposition~\ref{prop:13gp,1}.}
To finish the proof of Proposition~\ref{prop:13gp,1}, we still need to
check that the standard replacement condition is satisfied with
uniform norm. The fact that $\upsilon_{\infty}$ restricted to each
individual domain replaces $\upsilon_{\infty}$ on $J$ is implied by
the fact that on each newly created domain $\upsilon_{\infty}$ is a
lift by diffeomorphism or folding branches, always $\epsilon$-extendible, 
of $\upsilon'$ from another branch. The technical details of this fact
are provided in~\cite{kus}, Lemma 4.6. Then the fact that
$\upsilon_{\infty}$ restricted to each box also replaces
$\upsilon_{\infty}$ from $J$ follows automatically. We just notice
that boxes of $\phi_{\infty}$, with the exception of the central
domain, are the same as the boxes of $\phi$. The construction of the
Sewing Lemma can be conducted for $\cal Y$ extending $\upsilon_{0}$
from one box only, so it will not affect the replacement condition.
This concludes the proof of Proposition~\ref{prop:13gp,1}.

\paragraph{Theorem about initial inducing.}
\begin{defi}\label{defi:22ga,1}
A real box map is called {\em suitable} if the critical value of the
central branch is in the central domain and stays there forever under
iterations of the central branch.
\end{defi}

The inducing construction described earlier in this section hangs up
on a suitable map. In is easy to see that if a map is suitable, there
must be an interval symmetric with respect to $0$ which is mapped
inside itself by the central branch. In fact, this interval must be a
restrictive interval of the original $f$.   

\noindent {\bf Theorem 2}\begin{em}
Suppose that $f$ and $\hf$ belong to ${\cal F}_{\eta}$ for some
$\eta>0$ and are topologically conjugate. Then, for every $\alpha>0$, 
we claim the existence
of conjugate  box mappings $\phi_{s}$ and $\hat{\phi}_{s}$ on
the respective fundamental inducing domains $J$ and $\hat{J}$, and
a branchwise equivalence $\upsilon_{s}$ between them so that a number
of properties are satisfied. 
\begin{em}
\begin{itemize}
\item
$\phi_{s}$ and $\hat{\phi}_{s}$ are either full, or of type I (both of
the same type.)
\item
$\phi_{s}$ and $\hat{\phi}_{s}$ both possess standard
$\epsilon$-extendibility in the sense given by
Definition~\ref{defi:10gp,1}.  
\item
$\upsilon_{s}$ satisfies the standard replacement condition with
distortion $K_{2}$ in the
sense of Definition~\ref{defi:13gp,1}. 
\item
$\upsilon_{s}$ is $K_{1}$-quasisymmetric. 
\end{itemize}
\end{em}
In addition, $\phi_{s}$ and $\hat{\phi}_{s}$ {\em either are both
suitable} or this set of conditions is satisfied:
\begin{em}
\begin{itemize}
\item
for any domain $D$ of $\phi_{s}$ or $\hat{\phi}_{s}$ which does not
belong to  monotone branch which ranges through the fundamental
inducing domain, 
\[ \frac{|D|}{\dist(D,\partial J)} \leq \alpha\; ,\]
\item
if the mappings are of type I, then 
\[ \frac{|B'|}{\dist(B,\partial J)}\leq \alpha\]
and the analogous condition holds for $\hat{\phi}_{s}$,
\item
they both have extensions as complex box mappings with hole
structures satisfying the geometric bound $K_{3}$.
\item
on the boundary of each hole, the mapping onto the box is $K_{4}$-
quasisymmetric. 
\end{itemize}
\end{em}

The bounds $\epsilon$ and $K_{i}$ depend only on $\eta$ and $\alpha$.
\end{em}

\subparagraph{An outline of the proof of Theorem 2.}
Theorem 2 is going to be the main result of this section. The strategy
will be to obtain $\phi_{s}$ and its counterpart in a bounded (in
terms of $\eta$ and $\alpha$) number of basic inducing operations.
Then $\epsilon$-extendibility will be satisfied by construction. The
main difficulty  will be to obtain the bounded hole structure. Next,
the branchwise equivalence will be constructed based on
Proposition~\ref{prop:13gp,1}, so not much extra technical work will
be required other than checking the assumptions on this Proposition.  

\subsection{Inducing process}
\paragraph{The starting point.} 
\begin{fact}\label{fa:22ga,1}
Under the hypotheses of Theorem 2, conjugate induced mappings
$\phi_{0}$ and $\hat{\phi}_{0}$ exist with a branchwise equivalence
$\upsilon_{0}$. Also, versions of $\phi_{0}$ and $\hat{\phi}_{0}$
exist which are refined at the boundary with appropriate branchwise
equivalences. In particular, $\upsilon_{0,b}$ is infinitely refined at
both endpoints. The following conditions are fulfilled for $\phi_{0}$,
$\hat{\phi}_{0}$ and all their versions refined at the boundary. For
simplicity, we state them for $\phi_{0}$ only. 
\begin{em}
\begin{itemize}
\item
$\phi_{0}$ has standard $\epsilon$-extendibility,
\item
every domain except for two at the endpoints of $J$ is adjacent to two
other domains and the ratio of lengths of any two adjacent domains is
bounded by $K_{1}$,
\item
$\upsilon_{0}$ restricted to any domain replaces $\upsilon_{0,b}$ with
distortion $K_{2}$ and coincides with $\upsilon_{0,b}$ outside of $J$,
\item
$\upsilon_{0}$ is $K_{3}$-quasisymmetric,
\item
there are two monotone branches which range through $J$ with domains
adjacent to the boundary of $J$, and the ratio of length of any of
them to the length of $J$ is at least $\epsilon_{1}$.
\end{itemize}
\end{em}

The estimates $\epsilon$, $\epsilon_{1}$ and $K_{i}$ depend only on
$\eta$. 
\end{fact}
\begin{proof}
This fact is a restatement of Proposition 1 and Lemma 4.1 of~\cite{kus}.      
\end{proof}

\paragraph{The general inducing step.}
Suppose that a real box mapping $\phi$ is given which is either full, of type
I, or of type II and not suitable. Let us also assign a {\em rank} to
this mapping
which is $0$ is the $\phi$ is full. By the assumption that the
critical value is recurrent, the critical value $\phi(0)$ is in the
domain of $\phi$. There are three main distinctions we make. 
\newline{\bf Close and non-close returns.}
If the critical value is in the central domain, the case is classified
as a {\em close return}. Otherwise, it is described as {\em
non-close}. For a close return, we look at the number of iterations of
the central branch needed to push the critical value out of the
central domain. We call it the {\em depth} of the close return. 
\newline{\bf Box and basic returns.}
The situation is described as {\em basic} provided that upon its first
exit from the central domain the critical value falls into the domain
of a monotone branch which ranges through $J$. Otherwise, we call it a
{\em box} case. 
\newline{\bf High and low returns.}
We say that a return is {\em high} if $0$ is the range of the central
branch, otherwise the case is classified as {\em low}. 
Clearly, all eight  combinations can occur, which together with three
possibilities for $\phi$ (full, type I, type II) gives us twenty-four
cases. In the description of procedures, we use the notation of
Definition~\ref{defi:9mp,1} for boxes of $\phi$. 

In all cases, we begin with minimal boundary-refinement. 
\subparagraph{A low, close and basic return.}
In this case the almost parabolic pull-back is executed. The
resulting box mapping is usually not of any classified type. This is
followed by the filling-in of all
branches ranging through $B'$ or $B$. This gives a type I mapping. 
Finally, the critical pull-back is applied to give us a full map.
\subparagraph{A low, close and box return.}    
Again, we begin by applying almost parabolic pull-back. Directly
afterwards, we apply simple critical pull-back. The result will always
be a type II mapping. 
\subparagraph{A high return, $\phi$ full or type I.}
Critical pull-back with filling-in is applied. The outcome is of type
I in the box case, and full in the basic case. 
\subparagraph{A high return, $\phi$ of type II.}
Apply filling-in to obtain a type I map, then follow the previous
case. 
\subparagraph{A low return, non-close and basic, $\phi$ not of type II.} 
Apply critical pull-back with filling-in. The outcome is a full map. 
\subparagraph{A low return, non-close and basic, $\phi$ of type II.}
Use filling-in to pass to type I, then follow the preceding step. 
\subparagraph{A low and box return.}
Apply simple critical pull-back. The outcome is a type II map. 

This is a well-defined algorithm to follow. The rank of the resulting
mapping is incremented by $1$ in the box case, and reset to $0$ in the
basic case.   

\subparagraph{Notation.}
A box mapping is of rank $n$, we will often denote $B_{n}:=B$ and
$B_{n'}:=B'$.  

\paragraph{Features of general inducing.}
\begin{lem}\label{lem:10ma,1}
If $\phi$ is a type I or II induced box mapping of rank $n$ derived from some
$f\in {\cal F}$ by a sequence of generalized inducing steps,  then 
\[ \frac{|B_{n}|}{|B_{n'}|} \leq 1- \epsilon\]
where $\epsilon$ is a constant depending on $\eta$ only. 
\end{lem}
\begin{proof} 
If $\phi$  is full, the ratio is indeed bounded away from $1$ since a
fixed proportion of $B_{0}$ is occupied by the domains of two branches
with return time $2$. In a sequence of box mappings this ratio remains
bounded away from $1$ by Fact 4.1 of~\cite{indue}.   
\end{proof}

\begin{lem}\label{lem:22gp,1}
Let $\phi_{n}$ denote the mapping obtained from $phi_{0}$ given by
Fact~\ref{fa:22ga,1} in a sequence of $n$ general inducing steps. 
There are sequences of positive estimates $\epsilon(n)$ and $\epsilon_{1}(n)$  
depending otherwise only on $\eta$ so that the following holds:
\begin{em}
\begin{itemize}
\item
$\phi_{n}$ has standard $\epsilon(n)$-extendibility,
\item
the margin of extendibility of branches ranging through $J$ remains
unchanged in the construction,
\item
all branches are $\epsilon(n)$-extendible,
\item
for every box $B$ smaller than $J$, both external domains exist,
belong to monotone branches mapping onto $J$, and the ratios of their
lengths to the length of $B$ are bounded from below by
$\epsilon_{1}(n)$.
\end{itemize}
\end{em}
\end{lem}
\begin{proof}
This lemma follows directly from analyzing the multiple cases of the
general inducing step. The second claim follows directly since we use
minimal boundary refinement. For the first claim, in most cases one
uses Lemma~\ref{lem:12gp,1}. The more difficult situation is encountered if a
simple     
critical pull-back occurs. In this situation, the extendibility of
short monotone branches may not be preserved. However, one sees by
induction that a short monotone domain is always adjacent to two long
monotone domains. The ratio of their lengths can be bounded in terms
of $\eta$ and $n$. So, if the extendibility of this short monotone
branch is drastically reduced, then the critical value must be in this
adjacent long monotone domain. But then a basic return has occurred in
which the general inducing step always uses critical pull-back with
filling-in, which preserves extendibility.
\end{proof}

\paragraph{Main proposition}
\begin{prop}\label{prop:28gp,1}
For every $\alpha>0$, there are a fixed integer $N$ and $\beta_{0}>0$,
both independent of
the dynamics, for which
the following holds. Given a mapping $\phi_{0}$ obtained from
Fact~\ref{fa:22ga,1}  there is a sequence
of no more than $N$ general  inducing steps followed by no more than
$N$ steps of simultaneous monotone pull-back which 
give an induced map $\phi$ that
satisfies at least one of these conditions:
\begin{itemize}
\item
$\phi$ is of type I and suitable,
\item
$\phi$ is full and $\alpha$-fine,
\item
$\phi$ is of type I, has a hole structure of a complex box mapping
with geometrical norm less than
$\beta_{0}$, and for
every domain of a short monotone branch or the central domain, the
ratio of its length by the distance from the boundary of $I$ is less
than $\alpha$.    
\end{itemize}
\end{prop}

Proposition ~\ref{prop:28gp,1} is a major step in proving Theorem 2.
Observe that the first and third possibilities are already acceptable as
$\phi_{s}$ in Theorem 2. 
\paragraph{Taking care of the basic case.}
We construct an induced
map $\phi(\alpha)$ according by applying the general inducing step 
until the central domain and all short monotone domains $D$ satisfy 
$|D|/\dist(D,\partial J)\leq \alpha$. This will take a bounded number
of general inducing steps because the sizes of central domains shrink
at a uniform exponential rate (see~\ref{lem:10ma,1}.)  Of course, we
may be prevented from getting to this stage by encountering a suitable
map before; in that case Proposition~\ref{prop:28gp,1} is already
proven. From the stage of $\varphi(\delta)$ on, whenever we hit a
basic return we are in the second  case of
Proposition~\ref{prop:28gp,1}. Indeed, the ratio of length of the
central domain or any short monotone domain to the distance from the
boundary of $J$ is small. What is left is only to shorten long
monotone, non-external domains. This is done in a bounded number of
simultaneous monotone pull-back steps.  
So we can assume that box returns
exclusively occur beginning from $\varphi(\alpha)$. 

\subsection{Finding a hole structure}
We consider the sequence $(\varphi_{k})$ of consecutive box mappings
obtained from $\varphi_{0} := \varphi(\alpha)$ in a sequence of
general inducing steps. The objective of this section is to show that for
some $k$ which is bounded independently of everything else in the
construction, a uniformly bounded hole structure exists which extends
$\varphi_{k}$ as a complex box mapping in the sense of
Definition~\ref{defi:10ga,1}.

\paragraph{The case of multiple type II maps.}
We will prove the following lemma:
\begin{lem}\label{lem:28gp,5}
Consider some $\varphi_{m}$ of rank $n$. There is a function $k(n)$
such that if the mappings $\varphi_{m+1},\ldots,\varphi_{m+k(n)}$ are
all of type II, then $\varphi_{m+k(n)}$ has a hole structure
which makes it a complex box mapping. 
The geometric norm of this hole structure is bounded 
depending solely on $n$ and $\eta$.   
\end{lem}
\begin{proof}
If a sequence of type II mappings occurs, that means that the image of
the central branch consistently fails to cover the critical point. The
rank of all branches is fixed and equal to $n$. The central domain
shrinks at least exponentially fast with steps of the construction at
a uniform rate by Lemma~\ref{lem:10ma,1}. 
Thus, the ratio of the length of the central domain of $\varphi_{m+k}$
to the length of $B_{n}$ is bounded by a function of $k$ which also
depends on $n$, and for a fixed $n$ goes to $0$ as $k$ goes to
infinity.

This means that we will be done if we show that a small enough value
of this ratio ensures the existence of a bounded hole structure.   
We choose two symmetrical circular arcs which intersect the line in
the endpoints of $B_{n}$ at angles $\pi/4$ to be the boundary of the
box.($\alpha$
will be chosen in a moment, right now assume $\alpha<\pi/2$).
We take the
preimages of the box by the monotone branches of rank $n$. They are
contained in similar circular sectors circumscribed on their domains
by Poincar\'{e} metric considerations given in~\cite{miszczu}.    
Now, the range of the central branch is not too short compared to the
length of $B_{n}$ since it at
least covers one branch adjacent to the boundary of the box (otherwise
we would be in the basic case). We claim that the
domain of that external branch constitutes a proportion of the
box bounded depending on the return time of the central branch.
Indeed, one first notices that the return time of the external branch
is less than the return time of the central branch. This follows
inductively from the construction. But the range of the external
domain is always the whole fundamental inducing domain, so the domain
of this branch cannot be too short. On the other hand, the size of the
box is uniformly bounded away from $0$. 
To obtain
the preimage of the box by the complex continuation of the central
branch, we write the central branch as $h(z-1/2)^{2}$. The preimage by
$h$ is easy to handle, since it will be contained in a similar
circular sector circumscribed on the real preimage. Since the
distortion of $h$ is again bounded in terms of $n$, the real range of
$(x-1/2)^{2}$ on the central branch will cover a proportion of the
entire $h^{-1}(B_{n})$ which is bounded away from $0$ uniformly in
terms of $n$. Thus, the preimage of the complex box by the central
branch will be contained in a star-shaped region which meets the real
line at the angle of $\pi/4$ and is contained in a rectangle built on
the central domain of modulus bounded in terms of $n$.
    
It follows that this preimage will be contained below in the complex box
if the middle domain is sufficiently small. As far as the geometric norm is
concerned, elementary geometrical considerations show that it is
bounded in terms of $n$.
\end{proof}
\paragraph{Formation of type I mappings.}
We consider then same sequence $\varphi_{k}$ and we now analyze the
cases when the image of the central branch covers the critical point.
Those are exactly the situations which lead to type I maps in general 
inducing with filling-in.  

\subparagraph{A tool for constructing complex box mappings.}  
The reader is warned that the notations used in this technical
fragment are ``local'' and should not be confused with symbols having
fixed meaning in the rest of the paper.
The construction of complex box mappings (choice of a bounded hole
structure ) in the remaining cases will be based on the technical work
of~\cite{limo}. We begin with a lemma which appears there without
proof.
\begin{lem}\label{lem:2ha,1}
Consider a quadratic polynomial $\psi$ normalized so that
$\psi(1)=\psi(-1)=-1$, $\psi'(0)=0$ and $\psi(0)=a\in (-1,1)$.
The claim is that if $a<\frac{1}{2}$, then $\psi^{-1}(D(0,1))$ is
strictly convex.
\end{lem}
\begin{proof}
This is an elementary, but somewhat complicated computation. We will
use an analytic approach by proving that the image of the tangent line
to $\partial\psi^{-1}(D(0,1))$ at any point is locally strictly
outside of $D(0,1)$ except for the point of tangency.   
We represent points in $D(0,1)$ in polar coordinates $(r,\phi)$
centered at $a$ so that $\phi(1)=0$, while for points in preimage we
will use similar polar coordinates $r',\phi'$. By school geometry we
find that the boundary of $D(0,1)$ is given by
\[ (r+a\cos\phi)^{2} + a^{2}\sin^{2}\phi = 1\; .\]
By symmetry, we restrict our considerations to $\phi\in[0,\pi]$. We
can then change the parameter to $t=-a\cos\phi$, which allows us to
express $r$ as a function of $t$ for boundary points, namely
\begin{equation}\label{equ:2ha,1}
r(t) = \sqrt{1-a^{2}+t^{2}} + t\; .
\end{equation}

Now consider the tangent line at the preimage of $(r(t),t)$ by $\psi$.
By conformality, it is perpendicular to the radius joining to $0$, so
it can be represented as the set of points $(r',\phi')$
\[ r'=\frac{\sqrt{r(t)}}{\cos^{2}(\theta/2)}\:,\; \phi' = \frac{\pi}{2}
- \frac{\phi}{2} + \frac{\theta}{2} \]
where $\theta$ ranges from $-\pi$ to $\pi$. The image of this line is
given by $(\hat{r}(\theta),\phi-\theta)$ where 
\[  \hat{r}(\theta) = \frac{2r(t)}{1 + \cos\theta} \; .\]
Let us introduce a new variable $t(\theta) := -a\cos(\phi-\theta)$ so
that $t=t(0)$. Our task is to prove that
\begin{equation}\label{equ:2ha,2}
r(t(\theta)) < \hat{r}(\theta)
\end{equation} for values of $\theta$ is some
punctured neighborhood of $0$. This will be achieved by comparing the
second derivatives with respect to $\theta$ at $\theta=0$.  
By the formula~[\ref{equ:2ha,1}]
\[ r(t(\theta)) = \sqrt{1-a^{2}+t(\theta)^{2}} + t(\theta)\; .\]
The second derivative at $\theta=0$ is
\[ \frac{1-a^{2}}{(\sqrt{1-a^{2}+t(\theta)^{2}})^{3}}(a^{2} - t^{2}) - t -
\frac{t^{2}}{\sqrt{1-a^{2}+t(\theta)^{2}}}\; = \]
\[ =-\sqrt{1-a^{2}+t(\theta)^{2}}-t +
\frac{1-a^{2}}{(\sqrt{1-a^{2}+t(\theta)^{2}})^{3}}\; .\] 
The second derivative of the right-hand side of the desirable
inequality~[\ref{equ:2ha,2}] is more easily computed as
\[ \frac{\sqrt{1-a^{2}+t(\theta)^{2}}+t}{2} \; .\]
Thus, the proof of the estimate~[\ref{equ:2ha,2}], as well as the
entire lemma, requires showing that
\[ \frac{3}{2}(\sqrt{1-a^{2}+t(\theta)^{2}}+t) -
\frac{1-a^{2}}{(\sqrt{1-a^{2}+t(\theta)^{2}})^{3}}  > 0\]
for $|a|<1/2$ and $|t|\leq |a|$.  For a fixed $a$, the value of this
expression increases with $t$. So, we only check $t=-a$ which reduces
to 
\[ \frac{3}{2}(1-a) - (1 - a^{2}) > 0\]
which indeed is positive except when $a\in [1/2,1]$.
\end{proof}
\subparagraph{The main lemma.}
Now we make preparations to prove another lemma, which is
essentially Lemma 8.2 of~\cite{limo}.
Consider three nested intervals $I_{1}\subset I_{0}\subset I_{-1}$
with the common midpoint at $1/2$. Suppose that a map $\psi$ is
defined on $I_{1}$ which has the form $h(x-1/2)^{2}$ where $h$ is
a polynomial diffeomorphism onto $I_{-1}$ with non-positive Schwarzian
derivative. We can think $\psi$ as the central branch of a generalized
box mapping.  We denote 
\[ \alpha := \frac{|I_{1}|}{|I_{0}|}\]. Next, if $0 < \theta\leq \pi/2$
we define $D(\theta)$ to be the union of two regions symmetrical with
respect to the real axis. The upper region is defined as the
intersection of the upper half plane with the disk centered in the
lower $\Re = 1/2$ axis so that its boundary crosses the real line at
the endpoints of $I_{0}$ making angles $\theta$ with the line. So,
$D(\pi/2)$ is the disk having $I_{0}$ as diameter.
\begin{lem}\label{lem:7np,1}
In notations introduced above, if the following conditions are satisfied:
\begin{itemize}
\item
$\psi$ maps the boundary of $I_{-1}$ into the boundary of $I_{0}$,
\item
the image of the central branch contains the critical point,
\item 
the critical value inside $I_{0}$, but not inside $I_{1}$,
\item
the distance from the critical value to the boundary of $I_{0}$ is no
more than the (Hausdorff) distance between $I_{-1}$ and $I_{0}$,
\end{itemize}
then $\psi^{-1}(D(\theta))$ is contained in $D(\pi/2)$ and the
vertical strip based on $I_{1}$.
Furthermore, for every $\alpha < 1$ there is a choice of
$0<\theta(\alpha)<\pi/2$ so that 
\[\psi^{-1}(D(\theta(\alpha)))\subset D(\theta(\alpha))\] with a
modulus at least $K(\alpha)$, and $\psi^{-1}(D(\theta(\alpha)))$ is
contained in the
intersection of two convex angles with vertices at the endpoints of
$I_{1}$ both with measures less than $\pi-K(\alpha)$. 
Here, $K(\alpha)$ is a continuous positive function. 
\end{lem}
\begin{proof}
By symmetry, we can assume that the critical
value, denoted here by $c$, is on the left of $1/2$. 
Then $t$ denotes the right endpoint of $I_{0}$, and $t'$ is the other
endpoint of $I_{0}$. Furthermore, $x$ means the right endpoint of
$B_{n-1}$. By assumption, $h$ extends to the
range $(t',x)$. 
To get the information about the preimages of points 
$t,t',c,x$ one considers their cross-ratio 
\[ C = \frac{(x-t)(c-t')}{(x-c)(t-t')} \geq \frac{1+\alpha}{4}\]
where we used the assumption about the position of the critical value
relative $I_{0}$ and $I_{-1}$. The cross ratio will not be decreased by
$h^{-1}$. In addition, one knows that $h^{-1}$ will map the disk of
diameter $I_{0}$
inside the disk of  diameter $h^{-1}(B_{n})$ by the Poincar\'{e}
metric argument of~\cite{miszczu}.  As a consequence of the
non-contracting property of the cross-ratio, we get
\begin{equation}\label{equ:1hp,1}
\frac{h^{-1}(c) - h^{-1}(t')}{h^{-1}(t) - h^{-1}(t')} <
\frac{1+\alpha}{4}\; .
\end{equation}
When we pull back the disk based on $h^{-1}(I_{0})$, we will get a
figure which intersects the real axis along $I_{1}$. Notice that by
the estimate~[\ref{equ:1hp,1}] and Lemma~\ref{lem:2ha,1}, the preimage
will be convex, thus necessarily contained in the vertical strip based
on $I_{1}$. Its height in the imaginary direction is
\begin{equation}\label{equ:1hp,2}
 \frac{|I_{1}|}{2}\sqrt{\frac{h^{-1}(t) - h^{-1}(c)}
{h^{-1}(c) - h^{-1}(t')}} <
 \frac{|I_{1}|}{2}\sqrt{\frac{3-\alpha}{1+\alpha}}\; ,
\end{equation}

 where we used the estimate~[\ref{equ:1hp,1}] in the last inequality.    
Clearly, 
\[ \psi^{-1}(D(\pi/2))\]
is contained in the disk of this radius
centered at $1/2$. To prove that 
\[ \psi^{-1}(D(\pi/2))\subset D(\pi/2)\; , \] in view of the 
relation~[\ref{equ:1hp,2}] we need
\begin{equation}\label{equ:2ha,3}
 \alpha\sqrt{\frac{3-\alpha}{1+\alpha}} < 1\; 
\end{equation}
By calculus one readily checks that this indeed is the case when
$\alpha<1$.
To prove the uniformity statements, we first observe that 
\[ \psi^{-1}(D(\theta)) \subset \psi^{-1}(D(\pi/2)) \; \]
for every $\theta < pi/2$. Since [\ref{equ:2ha,3}] is a sharp
inequality, for every $\alpha<1$ there is some range of values of
$\theta$ below $\pi/2$ for which $\psi^{-1}(D(\theta))\subset
D(\theta)$ with some space in between. We only need to check the
existence of the angular sectors. For the intersection of 
$\psi^{-1}(D(\theta))$ with a narrow strip around the real axis, such
sectors will exist, since the boundary intersects the real line at
angles $\theta$ and is uniformly smooth. Outside of this narrow strip,
even $\psi^{-1}(D(\pi/2))$ is contained in some angular sector by its
strict convexity.
\end{proof}

The assumption of
extendibility to the next larger box is always satisfied in our
construction. 
\subparagraph{The case when there is no close return.}
We now return to our construction and usual notations. We consider a
map $\varphi_{k}$, type II and of rank $n$,
 whose central branch covers the critical
point, but without a close return. Then:
\begin{lem}\label{lem:2ha,2} 
Either the Hausdorff distance from $B_{n}$ to $B_{n-1}$ exceeds the
Hausdorff distance from $B_{n-1}$ to $B_{n-2}$, or $\varphi_{k}$ has a
hole structure uniformly bounded in terms of $n$.
\end{lem}
\begin{proof}
Suppose the condition on the Hausdorff distances fails. We choose the
box around $B_{n}$ and the hole around $B_{n+1}$ by
Lemma~\ref{lem:7np,1}. Observe that the quantity $\alpha$ which plays
a role in that Lemma is bounded away from $1$ by
Lemma~\ref{lem:10ma,1}. The box is
then pulled back by these monotone branches and its preimages are
inside similar figures built on the domains of branches by the usual 
Poincar\'{e} metric argument of~\cite{miszczu}. For those monotone
branches, the desired bounds follow immediately.
\end{proof}

\paragraph{Proof of Proposition~\ref{prop:28gp,1}.}
This is just a summary of
the work done in this section. We claim that we have proved that
either a map with a box structure can be obtained from
$\varphi(\delta)$ in a uniformly bounded number of steps of general
inducing, or the Proposition~\ref{prop:28gp,1} holds
anyway. Since Lemmas~\ref{lem:28gp,5} and~\ref{lem:2ha,2} provide
uniform bounds for the hole
structures in terms of $k$ or the rank which is bounded in terms of
$k$, it follows that the hole structure is bounded or the starting
condition holds anyway.  If the
inducing fails within this bounded number of steps because of a
suitable map being reached, then the stopping time on the central
branch of the suitable map is bounded, hence
Proposition~\ref{prop:28gp,1} again follows. 

So, we need to prove that claim. If the claim fails, then by
Lemma~\ref{lem:28gp,5} the situations in which the image of the
central branch covers the critical point have to occur with definite
frequency. That is, we can pick a function $m(k)$ independent of other
elements of the construction which goes to infinity with $k$ such that
among $\varphi_{1},\ldots,\varphi_{k}$ the situation in which the
critical point is covered by the image of the central branch occurs at
least $m(k)$ times. But each time that happens, we are able to
conclude by Lemma~\ref{lem:2ha,2} that the
Hausdorff distance between more deeply nested boxes is more than
between shallower boxes. Initially, for $\varphi_{0}$ whose  rank
was $n$, the $B_{n}$ distance between and $B_{n-1}$ was a fixed
proportion of the
diameter of $B_{n}$.  So
only a bounded number of boxes can be nested inside $B_{n-1}$ with
fixed space between any two of them.
So we have a bound on
the value of $m(k)$, thus on $k$. This proof of the
claim is a generalization of the reasoning used in~\cite{limo}.
The claim concludes the proof of Proposition~\ref{prop:28gp,1}.    

\subsection{Proof of Theorem 2}
\paragraph{Immediate cases of Theorem 2.}
Given two conjugate mappings $\phi_{0}$ and $\hat{\phi}_{0}$ obtained
by Fact~\ref{fa:22ga,1}, we apply the generalized inducing process to
both. By Proposition~\ref{prop:28gp,1} after a bounded number of
generalized inducing steps we get to one of the three possibilities
listed there: a suitable map, a full map, or a type I complex box map
with a bounded hole structure. In the first and third cases all we
need in order to conclude the proof is the existence of a branchwise
equivalence $\upsilon_{s}$ with needed properties. In the second case
extra work will be required. 

We will show that generally after $n$ general inducing
steps stair from a pair of conjugate $\phi_{0}$ and $\hat{\phi}_{0}$
we get a branchwise equivalence $\upsilon$ which satisfies    
\begin{em}
\begin{itemize}
\item
$\upsilon$ coincides with $\upsilon_{0,b}$ outside of $J$,
\item
$\upsilon$ satisfies the standard replacement condition with
distortion $K_{2}(n)$ in the
sense of Definition~\ref{defi:13gp,1}. 
\item
$\upsilon_{s}$ is $K_{1}(n)$-quasisymmetric. 
\end{itemize}
The bounds $K_{1}(n)$ and $K_{2}(n)$ depend only on $\eta$ and $n$. 
\end{em}

This follows by induction with respect to $n$ from
Proposition~\ref{prop:13gp,1} using Lemma~\ref{lem:22gp,1}. 
This means Theorem 2 has been proved except in the second case of
Proposition~\ref{prop:28gp,1}. In the remaining case, we still have
the branchwise equivalence with all needed properties. 

\paragraph{The case of $\phi$ full.}
We need to construct a hole structure. 
Observe that we
can assume that the range of the central branch of  $\phi$ is not
contained in an external
branch. In that case we could compose the central branch  with this
external branch until the critical value leaves the external domain.
This would not change the branchwise equivalence, box or extendibility
structures of the map.  Also, instead of $\phi$ rather consider the
version $\phi_{b}$ refined at the endpoint not in the range of the
central domain enough times to make all domains inside this external
domain of $\phi$ shorter than some $\alpha'$.  

\subparagraph{A full map with two basic returns.}
Let us first assume that $\phi_{b}$ shows a basic return. 
Carry out a
general inducing step. The resulting full mapping $\phi_{1}$ has short
monotone branches which are $\epsilon'$-extendible and $\epsilon'$
goes to $1$ as $\alpha$ and $\alpha'$ go to $0$. Indeed, these short branches
extend with the margin equal to the central central domain of $\phi$,
and if $\alpha$ is small this is much larger than the central domain
of $\phi_{1}$. Let $\delta$ denote
the distance from the critical value of $\phi_{1}$ to the boundary of
$J$. We claim that for very $\delta_{0}>0$ there is an $\alpha>0$,
otherwise only depending on $\eta$, so that if $\phi$ was
$\alpha$-fine, then a bounded hole structure exists. Indeed, take the
diamond neighborhood with height $\beta$ of $J$ as the complex box. 
Its preimage by the central branch is quasidisc with norm depending of
$\delta$, $\beta$ and the extendibility (thus ultimately on $\eta$.)    
If $\beta$ is very small, the preimage is close to the ``cross'' which
is the preimage of $J$ by the central branch is the complex plane. In
particular, for $\beta$ small enough it fits inside a rombe
symmetrical with respect to the real axis with the central domain as a
diagonal. The diameter of this rombe is bounded in terms of $\delta$.  
If the extendibility of short monotone branch is sufficiently good,
then the preimages of this rombe by the short monotone branches will
be contained in similar rombes around short monotone domains (see
Fact~\ref{fa:13gp,1} and use K\"{o}be's distortion lemma.) This gives
us a bounded hole structure. 

Finally, we show how to modify $\phi_{1}$ so that its critical value
is in a definite distance from the boundary of $J$. If the range of
the central branch is very small, this is very easy. Just compose the
central branch with the external branch a number of times to repel the
critical value from the endpoint, but so that it is still inside the
external domain, thus giving a basic return. This may require an
unbounded number of compositions, but they will not change the
branchwise equivalence. If the range of the central branch is almost
the entire $J$, choose a version of $\phi_{1}$ which is refined to the
appropriate depth at the endpoint of $J$ not in the range of the
central branch. The appropriate depth should be chosen so that the
image of the critical value by the external branch of the refined
version is still in an external domain, but already in a definite 
distance from the boundary. The possibility of doing this follows
since we can bound the eigenvalue of the repelling periodic point in
the boundary of $J$ from both sides depending on $\eta$ (see Fact 2.3
in~\cite{kus}.) Then apply the general inducing step to this version
of $\phi_{1}$, and the construct the hole structure as indicated above
for the resulting map. 

Since in this process we only use a bounded number of inducing
operations, or operations that do not change the branchwise
equivalence, the branchwise equivalence between the maps for which we
constructed hole structures will satisfy the requirements of Theorem
2. So, in the case of a double basic return we finished the proof.      

\subparagraph{A box return.}
In this case, we apply the general inducing step once to get a type II
mapping $\phi_{1}$. Observe that the central domain of $\phi_{1}$ is
very short compared to the size of the box (the ratio goes to $0$ with
$\alpha$, independently of $\alpha'$. ) Again, if the critical value
is in a definite distance from
the boundary of the box, compared with the size of the box, then we
can repeat the argument of Lemma~\ref{lem:28gp,5} to build a hole
structure. Otherwise, the critical value is in a long monotone branch
external in the box. So one more application of the general inducing
step will give us a mapping which is twice basic, so we apply the
previous step. Again, we see that the branchwise equivalence satisfies
the requirements of Theorem 2. This means that Theorem 2 has been
proved in all cases.  

\section{Complex pull-back}
\subsection{Introduction}
We will introduce a powerful tool for constructing branchwise
equivalences while preserving their quasiconformal norm. Since the
work is done in terms of complex box mappings, we have to begin by
defining an inducing process on complex box mappings. 
\paragraph{Complex inducing.}
\subparagraph{A simple complex inducing step.}
Suppose that $\phi$ is a complex box mapping which is either full or
of type I. We will define an inducing step for $\phi$ which is the
same as the general inducing used in the previous section on the real,
with the only difference that infinite boundary-refinement is used.
First, perform the infinite boundary refinement. This has
an obvious meaning for complex box mappings. Namely, one finds bad
long monotone branches and composes their analytic continuations with
$\phi_{r}$. Call the resulting map $\phi'$. Next, replace the central
branch of $\phi'$ with the
identity to get $\phi^{r}$. Next, construct $\tilde{\phi}$ which is the
same as $\phi$ outside of the central domain, and put 
\[\tilde{\phi}:=\phi^{r}\circ\phi\]
on the central domain. Finally, fill in all short univalent branches
of $\tilde{\phi}$ to get a type I or full mapping. Observe that on the
real line this is just boundary refinement followed by critical
pull-back with filling-in, so standard extendibility is preserved. 

We also distinguish {\em complex inducing without
boundary-refinement}. This is the same as the procedure described
above only without boundary refinement.  

\subparagraph{A complete complex inducing step.}
A complete complex inducing step, which will also be called a {\em
complex inducing step} is the same as the simple inducing step just
defined provided that $\phi$ does not show a close return, that is,
the critical value of $\phi$ is not in the central domain. If a close
return occurs, a complete inducing step is a sequence of simple
inducing steps until a mapping is obtained which shows a non-close
return. Then the simple inducing step is done once again and the whole
procedure gives a complete inducing step. If simple complex inducing
steps are used without boundary refinement, we talk about a {\em
complex inducing step} without boundary refinement.

\paragraph{External marking.}
Consider two equivalent conjugate complex box mappings, $\varphi$ and
$\hat{\varphi}$, see Definition~\ref{defi:10ga,1} of complex box
mappings.  
\begin{defi}\label{defi:3ha,1}
A quasiconformal homeomorphism $\Upsilon$ is called an {\em externally
marked branchwise equivalence} if it satisfies this list of conditions:
\begin{em}
\begin{itemize}
\item
restricted to the real line, $\Upsilon$ is a branchwise equivalence in
the sense of Definition~\ref{defi:25na,1},
\item
$\Upsilon$ maps each box of $\varphi$ onto the corresponding box of
$\hat{\varphi}$,  
\item
$\Upsilon$ maps each complex domain of $\varphi$ onto a complex domain
of $\hat{\varphi}$ so that these domains range through corresponding boxes, 
\item
On the union of boundaries of all holes complex domains of $\varphi$,
the functional equation
\[ \hat{\varphi} \circ \Upsilon = \Upsilon \circ \varphi\] holds. 
\end{itemize}
\end{em}
\end{defi}

The last condition of this definition will be referred to as the {\em
external marking} condition by analogy to {\em internal marking} which will be
introduced next. 

\paragraph{Internal marking.}
\begin{defi}\label{defi:7ka,1}
Let $\upsilon$ be a branchwise equivalence. An {\em internal marking 
condition} is defined to a choice of a set $S$ so that $S$ is
contained in the union of monotone rank $0$ domains of $\upsilon$ and 
each such domain contains no more than one point of $S$. The
branchwise equivalence $\upsilon$ will be said to {\em satisfy} the
internal marking condition if $\upsilon$ coincides with the conjugacy
on $S$.
\end{defi}

\begin{defi}\label{defi:7ka,2}
We will say that a branchwise equivalence $\Upsilon$ is {\em
completely internally marked}
if for each internal marking condition $\Upsilon$
can be modified without changing it on the boundaries of holes and
diamonds so that the marking condition is satisfied. We will say that
an estimate, for example a bound on the quasiconformal norm, is
satisfied for a fully internally marked $\Upsilon$, if all those
modifications can be constructed so as to satisfy this estimate. 
\end{defi}   

\begin{defi}\label{defi:25ga,1}
Two generalized induced mappings $\phi$ and $\hat{\phi}$ are called 
{\em equivalent} provided that a
branchwise equivalence $\upsilon$ 
exists between them such that for every domain $D$ of $\phi$, we have
\[ \hat{\phi}(\upsilon(D)) = \upsilon(\phi(D)) \; .\]
\end{defi}
Equivalence  of induced mappings is weaker than their topological
conjugacy and it merely means that the branchwise equivalence respects
the box structures, that corresponding branch range through
corresponding boxes, and that the critical values are in corresponding
domains. 
   
We are ready to state out main result.\newline
{\bf Theorem 3}\begin{em}
Suppose that $\phi$ and $\hat{\phi}$ are conjugate complex box
mappings with diamonds of type I or full. Suppose that $\Upsilon$ is
$Q$-quasiconformal externally marked and completely internally marked
branchwise equivalence. Also, a branchwise equivalence $\Upsilon_{b}$
is given between the versions of $\phi$ and $\hat{\phi}$ infinitely
refined at the boundary. The map $\Upsilon_{b}$ is also
$Q$-quasiconformal, externally marked and completely internally marked
and coincides with $\Upsilon$ on the boundaries of all boxes. Suppose
also that $\Upsilon$ or $\Upsilon_{b}$ restricted to the domain of any
branch that ranges through $J$ replaces $\Upsilon_{b}$ on $J$ with
distortion $K'$ and suppose $\phi$ and $\hat{\phi}$ satisfy standard
$\epsilon$-extendibility.   
 Suppose that $\phi_{1}$ and $\hat{\phi}_{1}$
are obtained from $\phi$ and $\hat{\phi}$ respectively after several
complex inducing steps, some perhaps without boundary refinement (but
always the same procedure is used on both conjugate mappings.) 
  Then, there is an externally marked and
completely internally quasiconformal branchwise equivalence
$\Upsilon_{1}$ between $\phi_{1}$ and $\hat{\phi}_{1}$. Furthermore,
if $\phi_{1}$ is full, then $\Upsilon_{1}$ is $Q$-quasiconformal.
Otherwise, $\Upsilon$ is $Q$-quasiconformal on the complement of the
central domain and the union of short univalent domains. In boundary
refinement is used at each step, then $\phi_{1}$ and $\hat{\phi}_{1}$
still have standard $\epsilon$-extendibility and
$\Upsilon_{1}$ restricted to any domain of a branch that ranges
through $J$ replaces $\Upsilon_{1}$ on $J$ with distortion $K$. The
number $K$ only depends on $\epsilon$ and $K'$. 
If the
construction without boundary refinement is used and only the box case
occurs, the theorem remains valid if $\phi$ and $\hat{\phi}$ are
complex box mappings without diamonds. 
\end{em}

\subsection{Proof of Theorem 3.} 
\paragraph{Complex pull-back.}
\subparagraph{Historical remarks.}
The line of ``complex
pull-back arguments'' initiated by~\cite{dohu} and used by numerous
authors ever since. We only mention the works which directly preceded
and inspired our construction. In~\cite{brahu} the idea of 
complex pull-back was applied to the situation with multiple
domains and images, all resulting from a single complex dynamical system.
According to~\cite{hub}, a similar
approach was subsequently used in~\cite{yoc} in the proof of
the uniqueness theorem for non-renormalizable
quadratic polynomials. Also,~\cite{limo} used a complex pull-back
construction to study metric properties of the so called ``Fibonacci
unimodal map''. 

\subparagraph{Description.}
Given $\phi$ and $\hat{\phi}$ full or of type I with an externally
marked branchwise equivalence $\Upsilon$, we will show how to build a
branchwise equivalence between $\phi_{1}$ and $\hat{\phi}_{1}$
obtained in just one complex inducing step. The first stage is
boundary refinement. Correspondingly, on the domain of each complex
branch $\zeta$ being refined, we replace $\Upsilon$ with 
\[ \hat{\zeta}^{-1}\circ\Upsilon_{b}\circ\zeta \]
where $\Upsilon_{b}$ is the branchwise equivalence between the
versions refined at the boundary which coincides with $\Upsilon$ on
the boundaries of boxes. The map obtained in this way is called
$\Upsilon'$. 

Next comes the ``critical pull-back'' stage. 
We replace $\Upsilon$ on the central domain with the lift to branched
covers 
\[ \hat{\psi}^{-1}\circ \Upsilon'\circ \psi\]
where $\psi$ and $\hat{\psi}$ are central branches. For this to be
well-defined, we need $\Upsilon(\psi(0))=\hat{\psi}(0)$. If this is a
basic return, this condition may be assumed to be satisfied because of
internal marking. In the box case, we have to modify $\Upsilon$ inside
the complex branch which contains the critical value. We can do this
modification in any way which gives a quasiconformal mapping and
leaves $\Upsilon$ unchanged outside of this complex domain. We also
choose the lift which is orientation-preserving on the real line. 
Call the resulting map $\tilde{\Upsilon}$. 

Finally, there is the infinite filling-in.   
This is realized as a limit process. Denote
$\Upsilon^{0}:=\tilde{\Upsilon}$. Then $\Upsilon^{i+1}$ is obtained
from $\tilde{\Upsilon}$ by replacing it on the domain of each short
univalent branch $\zeta$ with $\hat{\zeta}^{-1}\circ \Upsilon^{i}\circ
\zeta$. Then one proceeds the limit almost everywhere in the sense of
measure. 

Note that the branchwise equivalence $\Upsilon_{1,b}$ between versions
of $\phi_{1}$ and $\hat{\phi}_{1}$ infinitely refined at the boundary
can be obtained in the same way using $\phi_{b}$ and $\hat{\phi}_{b}$
instead of $\phi$ and $\hat{\phi}$ respectively. Also observe that
this procedure automatically gives an externally marked branchwise
equivalence, and $\Upsilon_{1}$ and $\Upsilon_{1,b}$ are equal on the
boundaries of all boxes.  
 
\subparagraph{Complex pull-back on fully internally marked maps.}
We observe that the construction of complex pull-back we just
described is well defined not only on individual branchwise
equivalences with holes, but also on families of fully internally
marked branchwise equivalences. This follows from the recursive nature
of the construction. Suppose that $\zeta$, the dynamics inside a hole  
or a diamond is used to pull back a branchwise equivalence $\Upsilon$.
We will show that any marking condition on newly created long monotone 
domains can be satisfied by choosing $\Upsilon$ appropriately marked. 
Indeed, long monotone branches of the induced map which arises in this
pull-back step are preimages of long monotone branches of the map
underlying $\Upsilon$. Thus, if $s\in S$, one should impose the
condition $\zeta(s)$ in the image. This can only lead to some
ambiguity if $\zeta$ is 2-to-1. In this case $\zeta$ is an univalent
map $H$ followed by a quadratic polynomial. One prepares two versions of
marking which are identical on the part of the real line not in the
real image of $\zeta$, and pulls back both of them by $H$. 
When they are finally pull back by the quadratic map, they will be the
same on the vertical line through the critical point of $\zeta$ as a
consequence of their being equal on the part of real axis not in the
real image. So one can then match the two versions along this vertical
line. 	   

\paragraph{Induction.}
Now Theorem 3 is proved by induction with respect to the number of
complex inducing steps. For one step, we just proved it. In the
general step of induction, the only this which is not obvious is why 
a $Q$-quasiconformal branchwise equivalence is regained after a basic
return even though the branchwise equivalence on the previous stage
was not $Q$-quasiconformal on short monotone domains and the central
domain. To explain this point, use the same notation as in the
description of complex pull-back. The mapping $\tilde{\Upsilon}$ is
$Q$-quasiconformal except on the union of domains of {\em short
monotone branches.} The key point is to notice that it is
$Q$-quasiconformal on the central domain, because this is a preimage
of a long monotone domain (we assume a basic return!). Then for 
$\Upsilon^{i}$ the unbounded conformal distortion is supported on the
set of points which stay inside short monotone branches for at least
$i$ iterations. The intersection of these sets has measure $0$, which
is very easy to see since each short monotone branch is an expanding
in the Poincar\'{e} metric of the box. So, $\Upsilon^{i}$ form a
sequence of quasiconformal mappings which converge to a quasiconformal
limit and their conformal distortions converge almost everywhere. This
limit almost everywhere is bounded by $Q$. By classical theorems about
convergence of quasiconformal mappings, see~\cite{lehvi}, the limit of
conformal conformal distortions equal to the conformal distortion of
the limit mappings. Thus, the limit is indeed $Q$-quasiconformal.

To finish the proof of Theorem 3 we have to check extendibility and
the replacement condition in the case when boundary refinement is
being used. Extendibility follows directly from
Lemma~\ref{lem:12gp,1}. Then, $\Upsilon_{1}$ restricted to any domain
whose branch ranges through $J$ is a pull-back by an extendible
monotone or folding branch of $\Upsilon$ (or $\Upsilon_{b}$) from
another such domain. The replacement condition is then satisfied,
which follows from Lemma 4.6 of~\cite{kus}.   

\subsection{The box case}
\paragraph{Box inducing.}
Suppose that a complex box mapping $\phi$ is given which is either
full or of type I and shows a box return. We then follow the complex
inducing step without boundary refinement. This will be referred to as
{\em box inducing} and is the same as the box inducing used
in~\cite{indue}. This certainly works on complex box mappings without
diamonds. 

\paragraph{Complex moduli.}
Given a complex box mapping $\phi$ of type I or full, consider the
annulus between
the boundary of of $B'$ and the boundary of $B$. Denote its modulus 
$v(\phi)$. By our definition of box mappings, $v(\phi)$ is always
positive and finite. 

\noindent{\bf Theorem 4}\begin{em}
Let $\phi$ be complex box mapping of type II. Suppose that $\phi$ has
a hole structure with a geometric bound not
exceeding $K$. Let $\phi_{0}$ denote the type I mapping obtained from 
$\phi$ by filling-in. Assume that box inducing can be applied to
$\phi_{0}$ $n$ times, giving a sequence maps $\phi_{i}$, $i=0,\cdots,
n$. Then there is number $C>0$ only depending on the bound of the hole
structure so that 
\[ v(\phi_{i}) \geq C i\]
for every $i$. 
\end{em}

\subparagraph{Related results.} 
We state the related results which will be used in the proof. This
does not exhaust the list of related results in earlier papers, see 
``historical comments'' below. The first one will be called ``the
starting condition''. 

\begin{defi}\label{defi:10ma,1}
We say that a type I or type II box mapping of rank $n$ satisfies the {\em
starting condition} with norm $\delta$ provided that  
$|B_{n}|/|B_{n'}| < \delta$ and if $D$ is a short
monotone domain of $\phi$, then also $|D|/\dist(D,\partial B_{n'})< \delta$. 
\end{defi} 

\begin{fact}\label{fa:8mp,1}
Let $\phi_{0}$ be a real box mapping, either of type I of rank $n$ or full.  
Pick $0< \tau < 1$ and assume that  the central branch is $\tau$-extendible.
For every $\tau$, there is a positive number $\delta(\tau)$ with the following
property.     
Suppose that $\phi_{i}$, $i\geq 0$ is
sequence of type I real box mappings such that
$\phi_{j+1}$ arises from $\phi_{j}$ by a box inducing step.
Let $\phi_{0}$ satisfy the starting condition with norm $\delta(\tau)$. 
Then, 
\[ |B_{n+i}|/|B_{(n+i)'}| \leq C^{i}\; \] 
where $C$ is an absolute constant less than $1$.   
\end{fact}
\begin{proof}
This is Fact 2.2 of~\cite{indue}. A very similar statement, but for a
slightly different inducing process is Proposition 1 of~\cite{yours}. 
\end{proof}

\begin{fact}\label{fa:24ga,1}
Let $\phi_{0}$ be a complex box mapping. Let
$\phi_{0},\phi_{1},\cdots,\phi_{n}$ be the sequence in which the next
map is derived from the preceding one by a box inducing step. Let
$B_{i}$ denote the central domain of $\phi_{i}$ on the real line, and
let $B_{n'}$ denote the box on the real line through which the central
branch ranges. Suppose that for some $C<1$ and every $i$
\[ \frac{|B_{n}|}{|B_{n'}|} \leq C^{i} \]
and suppose that the hole structure of $\phi_{0}$ satisfies a
geometric bound $\beta$. For every $C<1$ and $\beta$ there is a
positive $\overline{C}$ so that 
\[ v(\phi_{i}) \leq \overline{C} i\]
for every $i$.
\end{fact}
\begin{proof}
This is Theorem D of~\cite{indue}. 
\end{proof}

\subparagraph{Historical comments.}
Theorem 4, in the way we state it, follows from Theorems B and D
of~\cite{indue}. However, we give a different proof based on Theorem
D, but not on Theorem B. Weaker results saying that 
$|B_{i}|/|B_{i'}|$ decrease exponentially fast were proved in various
cases in~\cite{yours},~\cite{limo} and later papers, including an
early preprint of this work.  In the Appendix we show a result similar
to Fact~\ref{fa:24ga,1}. Even though this result is not a necessary
step in the proof of the Main Theorem, we think that it may be of
independent interest and also its proof shows the main idea of the
proof of Fact~\ref{fa:24ga,1}

\paragraph{Artificial maps.}
\subparagraph{Introduction.} 
Artificial maps were introduced in~\cite{limo}. In our language,
\cite{limo} showed how to prove that the starting condition 
must be satisfied at some stage of the box construction for a concrete
``Fibonacci'' unimodal map. \footnote{The Fibonacci map is
persistently recurrent in the sense of~\cite{yoc}, or of infinite box
type in the sense of~\cite{yours}.} The idea
was to use an artificial map ``conjugated'' to the box map obtained
as the result of inducing. Clearly, an artificial map can be set up so
as to satisfy the starting condition. Next, one shows that the induced
map and its artificial counterpart are quasisymmetrically conjugate.
Since for the artificial map the starting condition is satisfied with
progressively better norms, after a few box steps it will be
forced upon the induced map by the quasisymmetric conjugacy. 

This strategy has a much wider range of applicability than the
Fibonacci polynomial and it is at the core of our proof of
Theorem 4. 

\subparagraph{Topological conjugacy.} 
\begin{lem}\label{lem:7ne,1}
Given a box mapping induced by some $f$ from $\cal F$ and any
homeomorphism 
of the line into itself, a
box map (artificial) can be constructed that is topologically
conjugate  to the original one and whose branches are all either
affine or quadratic and folding.
The homeomorphism then becomes a branchwise equivalence.  
\end{lem}
\begin{proof}
The domains and boxes of the artificial map that we want to construct
are given by the images in the homeomorphism. Make all monotone
branches affine and the central branch quadratic. We show that by
manipulating the critical value we can make these box maps equivalent.
We prove by induction that if two maps are similar, by just changing
the critical value of one of them we can make them equivalent. The
induction proceeds with respect to the number of steps needed to
achieve the suitable map. The initial step is clearly true by the
continuity of the kneading sequence in the $C^{1}$ topology, and the
the induction step consists in the remark that by manipulating the
critical value in similar maps we can ensure that they remain similar
after the next inducing and the critical value in the next induced map
can be placed arbitrarily.  
\end{proof}
\paragraph{The construction of an original branchwise equivalence.}
We will prove the following lemma:
\begin{lem}\label{lem:2hp,11}
Under the assumptions of Theorem 4 and for every $\delta>0$ there is
an artificial map $\Phi$ with all branches either affine or
quadratic and folding which is conjugate to $\phi$. Next, $\Phi$
also has a hole structure which makes it a complex box mapping of type
II. Also, the type I mapping obtained from $\Phi$ by filling-in
satisfies the starting condition with norm $\delta$. Also, we claim
that an externally marked branchwise equivalence exists between $\phi$
and $\varphi$ which is $Q$ quasiconformal. The number $Q$ only depends
on $\delta$ and the geometric norm of the hole structure for $\phi$
($K$ in Theorem 4). 
\end{lem} 
\begin{proof}
First construct the artificial map $\Phi$. 
To this end, we choose a diffeomorphism $h$ which is the identity outside
$B'$ and squeezes the central branch so that the after filling-in the
type I  map
satisfies the starting condition. We 
observe that the ``nonlinearity'' $h''/h'$ can be made bounded in
terms of the bound on the hole structure. 
Then by Lemma~\ref{lem:7ne,1} we can
adjust the critical value of this artificial map, without moving
domains of branches around, so that the map is equivalent to
$\phi$. This defines $\Phi$. 
 
Since $\Phi$ has affine or quadratic branches, the existence of a
bounded hole
structure equivalent to the structure already existing is clear. One
simply repeats the arguments of Lemmas~\ref{lem:2ha,2} 
and~\ref{lem:28gp,5} with ``infinite extendibility'' which makes the
problem trivial. 

The final step is to construct the branchwise equivalence by
Lemma~\ref{lem:11na,2}. First, we decide that on the real line the
branchwise equivalence is $h$. We next
define $\Upsilon$ inside the complex box around $B'$. On the
boundary of the box, we just take a bounded quasiconformal map. Again
by K\"{o}be's lemma this propagates to the holes with bounded
distortion. Since $h$ is a diffeomorphism of bounded nonlinearity
inside each hole, it can be filled with a uniformly quasiconformal
mapping. Lemma~\ref{lem:11na,2} works to build $\Upsilon$ inside
the complex box belonging to $B'$ with desired uniformity. Then the
map is extended outside of the box. This can also be achieved by
Lemma~\ref{lem:11na,2} regarding $B'$ as a tooth, holes outside of the
box as other teeth, and choosing a big mouth.   
\end{proof}

\subsection{Marking in the box case}
In the box case Theorem 3 does not imply that quasiconformal norms
stay bounded. On the other hand, internal marking does not work either
in the box case. We should a special procedure of achieving an
internal marking condition. We remind the reader that $v(\phi)$
denotes the modulus between $B'$ and $B$ when $\phi$ is of type I, or
between $B_{0}$ and $B$ when $\phi$ is full. 

\begin{prop}\label{prop:4hp,1}
Let $\tilde{\phi}$ and $\tilde{\Phi}$ be conjugate complex box
mappings either full
or of type I, not suitable  and showing a box return. Suppose that
$\phi$ and $\Phi$ are derived from $\tilde{\phi}$ and $\tilde{\Phi}$
respectively, and not suitable and show a box return. 
Let $\Upsilon$ be a $Q$-quasiconformal branchwise
equivalence acting into
the phase space of $\Phi$. Let $v$ denote the minimum of
$v(\tilde{\Phi})$ and $v(\Phi)$. 
Then there is a branchwise equivalence
$\Upsilon'$ with the following properties:
\begin{em}
\begin{itemize}
\item
$\Upsilon'$ equals $\Upsilon$ except on the complex domain of $\phi$
which contains the critical value of $\phi$,
\item
$\Upsilon(\phi(0)) = \Phi(0)$, 
\item
the quasiconformal norm of $\Upsilon'$ is bounded by 
\[ Q + K_{1}\exp(-K_{2} v\]
where $K_{1}$ and $K_{2}$ are positive constants. 
\end{itemize}
\end{em}
\end{prop}

Note that the estimate of the quasiconformal norm is independent of
the geometry of $\phi$. 

\paragraph{An auxiliary lemma.}
\begin{lem}\label{lem:4hp,1}
Let $\Phi$ and $\tilde{\Phi}$ and $v$ be as in the statement of
Proposition~\ref{prop:4hp,1}. Let $B_{0}\supset B'\supset B$ be the
box structure of $\Phi$. Choose a complex domain $b\in B'$ of $\Phi$
which is either short monotone. There exists an annulus $A$
and a constant $0<C<1$ with the following properties:
\begin{em}
\begin{itemize}
\item
$A$ surrounds $b$ and is contained in $B'$.
\item
the modulus of $A$ is at least $C v$ and the modulus of
the annulus separating $A$ from the boundary of $B'$ is at least $C
v$,
\item
the boundary of $A$ intersects the real line at four points, and these
points are topologically determined, meaning that if $A'$ is
constructed for a mapping $\Phi'$ conjugated to $\Phi$, the conjugacy
will map these four points onto the four points of intersection
between $A'$ and the real line.
\end{itemize}
\end{em}
\end{lem}
\begin{proof}
We split the proof in two cases depending on whether $\tilde{\Phi}$
showed a close return or not. In the case of a non-close return, the
outer boundary of $A$ is chosen as the boundary of the domain of the
analytic extension of the branch from $b$ onto the range $B'$. This
means that the annulus of $A$ is equal to $v(\Phi)$. On the other
hand, the domain of this extension is contained in the preimage by the
central branch of $\tilde{\Phi}$ of some domain of $\tilde{\Phi}$. 
A lower bound by $v(\tilde{\Phi})/2$ is evident. Also, the topological
character of this construction is clear. 

In the case of a close
return, however, this construction would not work because the annulus
between the extension domain and the boundary of $B'$ may become
arbitrarily small. So, let us consider the central domain of $\Phi$. 
Recall that the box inducing step in the case of a close return is a
sequence of simple box inducing steps with close returns ended by a
simple box inducing step with a non-close return. Let $\Phi_{1}$ mean
the mapping obtained for $\tilde{\Phi}$ by this sequence of simple box
inducing steps with close returns. 
The central branch of $\Phi$ is the composition of a restriction of
$\psi_{1}$ (the central branch of $\Phi_{1}$) with a short univalent
branch $b_{1}$ of $\Phi_{1}$. By construction, this short univalent branch has
an analytic continuation whose domain is contained in the range of
$\psi_{1}$ and whose range is the range of the central branch of
$\tilde{\Phi}$. Inside this extension domain, there a smaller domain $\delta$
mapped onto the central domain of $\Phi$ only. Consider the annulus
$W$ between $B$ and the boundary of $\psi^{-1}(\delta)$. The modulus
of $W$ is a half of the modulus of the annulus between $b_{1}$ and the
boundary of $\delta$, which is $v(\tilde{\Phi})\cdot 2^{-l}$ where $l$
is the number of subsequent iterations of $\psi$ which keep the
critical value inside the central domain. It follows that the modulus
of $W$ is at least $v(\tilde{\Phi})/4$. Next, look at the annulus $W'$
separating $W$ from the boundary of $B'$. This is at least a half of
the annulus separating $\delta$ from the larger extension domain (onto
the range of $\psi$.) This last annulus has modulus $v(\tilde{\Phi})$.
So $W'$ has modulus at least $v(\tilde{\Phi})/2$. Now, to pick $A$
take the extension of the branch from $b$ onto $B'$ and define $A$ to
be the preimage of $W$ by this extension. Since $A$ is separated from
the boundary of the extension domain, let alone from the boundary of
$B'$, by the preimage of $W'$, the estimate claimed by this Lemma
follows. Also, the topological character of the intersection of $A$
with the real line is clear. 
\end{proof}

\paragraph{Proof of Proposition~\ref{prop:4hp,1}.}  
Suppose that the first exit time of the critical value from the
central domain under iteration by the central branch is $l$. Call 
$C=\Phi^{l}(0)$ and $c=\phi^{l}(0)$. 
Suppose $C$ belongs to a short monotone domain $b$. 

The point $\Upsilon(c)$ must also be in $b$. 
Take the annulus $A$ found for $b$ from Lemma~\ref{lem:4hp,1}. Also,
call $A'$ the annulus separating $A$ from the boundary of $B'$. 
Perturb $\Upsilon$ to $\Upsilon_{1}$ so that
$\Upsilon_{1}(c)=C$ and $\Upsilon=\Upsilon_{1}$ outside of
$\Upsilon^{-1}(A)$. We accept as obvious that this can be done by composing
$\Upsilon$ with a mapping whose conformal distortion is bounded by 
$k_{1} \exp(-k_{2}\mod A)$. By Lemma~\ref{lem:4hp,1}, note that this is a
correction allowed by Proposition~\ref{prop:4hp,1}. 
Then, we apply  
complex pull-back by the central
branch $\psi$ to $\Upsilon_{1}$ $l$-times. Call the resulting
branchwise equivalence
$\Upsilon_{2}$. If $l=1$ it is clear that the conformal distortion of
$\Upsilon_{2}$ is the same as for $\Upsilon_{1}$. If $l>1$ it less
clear since we will have to adjust the branchwise equivalence $l-1$
times inside shrinking preimages of $B$ to make $c$ and $C$
correspond. So until just before the last simple box inducing step the
conformal distortion is not well-controlled. However, we show as in
the proof of Theorem 3 that in this last simple box inducing step, the
region where the distortion was unbounded has preimage of measure $0$
because of filling-in,
so that ultimately the distortion of $\Upsilon_{2}$ is the same as for
$\Upsilon_{1}$ almost everywhere. 

Next, construct $\Upsilon_{3}$ 
which is the same as $\Upsilon$ outside of $\Upsilon^{-1}(b)$ and
equals $\Phi^{-1}\circ \Upsilon_{2} \circ \phi$ on $b$. Both mappings
match continuously because of external marking. Inside $b$, now look
at $A_{3} = \Phi^{-l-1}(A)$ and $A'_{3} = \Phi^{-l-1}(A')$. By
construction, the mapping $\Phi^{l+1}$ in the region encompassed by
these annuli is a branched cover of degree $2$, so $\mod A_{3} =
\frac{1}{2}\mod A$ and $\mod A'_{3} = \frac{1}{2}\mod A'$. We consider
two cases. If $C$ belongs to the region encompassed by the outer
boundary of $A$, then $\Upsilon_{3}(c)$ also belongs to this region.
This is because the points of intersection of $A$, and therefore of
$A_{3}$, with the real line are topologically determined by
Lemma~\ref{lem:4hp,1}. 
So, we perturb $\Upsilon_{3}$ to leave it unchanged
outside of $\Upsilon^{-1}(b)$ and to make the critical values
correspond. Like in the previous paragraph, we claim that this
adjustment will only add $k_{1}\exp(-k_{2}\mod A'_{3})$ to the
conformal distortion. In this case, we are done with the proof of
Proposition~\ref{prop:4hp,1}. Otherwise, $C$ belong to the preimage of
some short univalent domain $b'$ of $\Phi$ by $\Phi^{l+1}$. Since $b'$
is nested inside $B'$ with a modulus at least $v(\Phi)$, then this
preimage is nested inside $b$ with a modulus at least half that large.
Also, $\Upsilon_{3}(c)$ must belong to the same preimage of $b'$.
Again, we adjust $\Upsilon_{3}$ to make the critical values correspond 
and are done with the proof.   
   
\paragraph{Proof of Theorem 4.}
Given $\phi$ from Theorem 4, construct a conjugate artificial mapping
$\Phi$ from Lemma~\ref{lem:2hp,1} choosing $\delta$ from
Fact~\ref{fa:8mp,1} for a large $\tau=1001$ (the artificial map has
arbitrary extendibility.) This will guarantee by Fact~\ref{fa:24ga,1}
that if the sequence $\Phi_{i}$ is derived from from $\Phi_{0}:=\Phi$
by box inducing, then the moduli $v(\Phi_{i})$ grow at least at a
uniform linear rate. Lemma~\ref{lem:2hp,1} also gives us an externally
marked and uniformly quasiconformal branchwise equivalence
$\Upsilon_{0}$ from the phase space of $\phi$ to the phase space of
$\Phi$. 

Now proceed by complex pull-back as defined in the proof of Theorem 3,
however do the marking corrections by Proposition~\ref{prop:4hp,1}. 
We get a sequence of uniformly quasiconformal branchwise equivalences
between $\phi_{i}$ and $\Phi_{1}$. Since quasiconformal mappings
preserve complex moduli up to constants, $v(\phi_{i}) \geq K
v(\Phi_{i})$ with $K$ only depending on the conformal distortion of
$\Upsilon_{0}$, thus ultimately (Lemma~\ref{lem:2hp,1}) only on the
geometric bound of the
hole structure of $\phi$. Theorem 4 follows.       

\section{Construction of quasisymmetric conjugacies}
Here is the main result of this section. \newline
{\bf Theorem 5}\begin{em}
Suppose that $f$ and $\hat{f}$ are topologically conjugate and belong
to some ${\cal F}_{\eta}$. Assume also that if $f$ is renormalizable,
then the first return time of its maximal restrictive interval to
itself is greater than $2$.  
Then for every $\delta>0$ there exist conjugate real box mappings,
$\phi^{e}$ and $\hat{\phi}^{e}$, either full or of type I  and
infinitely refined at the boundary, 
with a branchwise equivalence
$\upsilon$ between them so that the following list of properties
holds:
\begin{em}
\begin{itemize}
\item
$\phi$ and $\hat{\phi}$ are either suitable or $\delta$-fine,
\item
$\phi$ and $\hat{\phi}$ both have standard $\epsilon$-extendibility, 
\item
$\upsilon$ is $Q$-quasisymmetric,
\item
$\upsilon$ restricted to any long monotone domain replaces $\upsilon$
on the fundamental inducing domain with distortion $K$,
\item
$\upsilon$ restricted to any short monotone domain replaces $\upsilon$
on $B$ with norm $K$.
\end{itemize}
\end{em}

The numbers $\epsilon>0$, $Q$ and $K$ depend on $\eta$ only. 
\end{em}  

\subsection{Towards final mappings}
\paragraph{Technical details of the construction.}
The next lemma tells that given mappings $\phi_{s}$, $\hat{\phi}_{s}$
and $\upsilon_{s}$ obtained from Theorem 2, we can modify an extend
$\upsilon_{s}$ to an externally marked and fully internally marked
quasiconformal branchwise equivalence. 

\begin{lem}\label{lem:25ga,1}
Let $\phi$ and $\hat{\phi}$ be topologically conjugate complex box
mappings, full or of type I, both infinitely refined at
the boundary. 
Suppose that 
\begin{em}
\begin{itemize}
\item
both have hole structures geometrically bounded by $K'$, 
\item on
the boundary of each hole the mapping is $K''$ quasisymmetric, 
\item 
for the domain $D$ of any branch of
$\phi$, $|D|/\dist(D,\partial J)\leq\alpha$
holds with some fixed $\alpha$; the same holds for every domain  of 
$\hat{\phi}$,
\item
if the mappings are of type I, then $|B'|/\dist(B',\partial J)\leq \alpha$ and
the same holds for $\hat{\phi}$,
\item
all long monotone branches of $\phi$ and
$\hat{\phi}$ and $\hat{\phi}$ are $\epsilon$-extendible, $\epsilon>0$,
\item
$\upsilon$ exists which is a completely internally marked branchwise 
equivalence between $\phi$
and $\hat{\phi}$,
\item
$\upsilon$ is $Q$-quasisymmetric,
\item
$\upsilon$ 
restricted to any domain of $\phi$ replaces $\upsilon$ on $J$ with
distortion $K$.
\end{itemize}
\end{em}

We claim that there a bound $\alpha_{0}$ depending on
$K'$ only so that if $\alpha<\alpha_{0}$, the following holds:
\begin{em}
\begin{itemize}
\item
$\phi$ and $\hat{\phi}$ can be extended to complex box mappings with
diamonds, call them $\phi^{d}$ and $\hat{\phi}^{d}$, 
\item
 an externally and completely
internally marked branchwise equivalence $\Upsilon_{0}$ 
exists between $\phi^{d}$ and
$\hat{\phi}^{d}$,
\item
$\Upsilon_{0}$ is $L$-quasiconformal,
\item
bounds 
$L$ and $\alpha_{0}$ only depend on $\epsilon$, $K$, $K'$, $K''$ and $Q$.
\end{itemize}
\end{em}
\end{lem}  
\begin{proof}
Let us first pick $\alpha_{0}$. Recall Fact~\ref{fa:13gp,1} and choose
$\kappa$ as $K_{1}$ picked for $\epsilon$ by this Fact. Make
$\kappa\leq 1/2$ as well. 
Consider the diamond neighborhood with
height $\kappa$ of $J$. The bound $\alpha_{0}$ should be picked so as to
guarantee that all holes of $\phi$, as well as the box $B'$ in case
$\phi$ is of type I, sit inside this diamond neighborhood, moreover,
that they have annular ``collars'' of definite modulus (say $1$) which
are also contained in this diamond neighborhood.    
All holes 
quasidisks bounded in terms of $K'$. So, there is a bounded ratio
between how far
they extend in the imaginary direction and the length of the real domain
they belong to. Now it is evident that $\alpha_{0}$ small will imply
this, and $\alpha_{0}$ depends only on $K'$. 
Also, by making $\alpha_{0}$ even smaller, we can ensure that diamond
neighborhoods with height $\kappa$ of all domains inside $B'$ are
contained inside the complex box corresponding to $B'$, also with
annular margins $1$. 

Next, we choose the diamonds. 
We will take the diamond neighborhood with height $\kappa$ of $J$ as the
complex box $B_{0}$. The diamonds will simply be preimages of this 
$B_{0}$ by all monotone branches. They will be contained
in diamond neighborhoods with height $\kappa$  of corresponding domains
of branches by the  Poincar\'{e} metric argument of ~\cite{miszczu}.
Also, the diamonds will be
preimages of $B_{0}$  with bounded distortion (Fact~\ref{fa:13gp,1}
again.) We can do the same thing for $\hat{\phi}$. This gives us
$\phi^{d}$ and $\hat{\phi}^{d}$. 

The final step is to construct the branchwise equivalence by
Lemma~\ref{lem:11na,2}. First, we decide that on the real line the
branchwise equivalence is $\upsilon$. On the box around $B_{0}$ extend
it in any way that maps $B_{0}$ onto $\hat{B}_{0}$ and gives a
uniformly quasisymmetric (in terms of $Q$
and $\kappa$) mapping on the union of $J$ and the upper (lower) half
of $B_{0}$. Then  $\Upsilon_{0}$ on the boundaries of
diamonds is determined by pull-back. Observe, however, that on the
curve consisting of the upper (lower) half of the diamond on the
domain on the real line, the map is quasisymmetric, and that is
because of the replacement condition. 
It follows that diamonds can be filled with uniformly
quasiconformal  mappings as Lemma~\ref{lem:11na,2} demands. 

We first
define $\Upsilon_{0}$ inside the complex box around $B'$. On the
boundary of the box, we just take a map
transforming it onto the boundary of $\hat{B}'$ and quasiconformal
with a bounded norm (in terms of $K'$ and $Q$) on the union of the
upper (lower) half of the boundary of the complex box and the real box
$B'$. This propagates to the holes with bounded deterioration of the
quasisymmetric norm (in terms of $K''$.)
Thus, each hole can be filled with a uniformly quasiconformal
mapping. Lemma~\ref{lem:11na,2} works to build $\Upsilon_{0}$ inside
the complex box corresponding to $B'$ with desired uniformity. 
This having been achieved, Lemma~\ref{lem:11na,2} is again used inside
the entire complex box around $B_{0}$. Here, the complex box around
$B'$ is formally regarded as a tooth. Finally, $\Upsilon_{0}$ is
extended to the plane by quasiconformal reflection.
\end{proof}

\paragraph{The initial branchwise equivalence.}
Suppose now that we are in the situation of Theorem 2 with mappings
$\phi_{s}$ and $\hat{\phi}_{s}$ not suitable. We want to build an
externally marked and completely internally marked branchwise
equivalence $\Upsilon_{s}$ between them. For that, we will 
use Lemma~\ref{lem:25ga,1} with $\phi:=\phi_{s,b}$ and
$\hat{\phi}:=\hat{\phi}_{s,b}$ (the additional subscript $b$ denotes
versions infinitely refined at the boundary.) Comparing the
assumptions of Lemma~\ref{lem:25ga,1} with claims of Theorem 2, we see
that two conditions that are missing are the complete internal
marking, and $|D|/\dist(D,\partial J)$ when $D$ is a long monotone
domain. Let us first show that the second property can be had by doing
more inducing on long monotone domains. On the level of inducing, we
simply compose long monotone domains whose domains are too large with
$\phi_{s}$ ($\hat{\phi}_{s}$ respectively) until we reduce their sizes
sufficiently. This can take many inducing steps on any given domain. Note
that the hole structure can simply be pull-back. The geometric
bound of the hole structure will be worsened only in a bounded fashion
provided that $\alpha$ was small enough. This follows from
Fact~\ref{fa:13gp,1}. Also, the replacement condition will not suffer
too much because we are pulling back by maps of bounded distortion (or
one can formally use Proposition~\ref{prop:13gp,1}). The only hard
point is the quasisymmetric norm. This does not directly follow from
Proposition~\ref{prop:13gp,1}, since we may have to do a large number
of simultaneous monotone pull-backs. However, this is easily seen if
we proceed by complex pull-back (like in the proof of Theorem 3). To
this end, we pick the diamond neighborhood of $J$ with height $1/2$ as
$B_{0}$, and a homothetic neighborhood of $\hat{J}$ as $\hat{B}_{0}$. 
The diamonds are the preimages of $B_{0}$ ($\hat{B}_{0}$ resp.) by long
monotone branches. We do not have any holes. Then the argument used in
the proof of Lemma~\ref{lem:25ga,1} applies and allows us to build a
branchwise equivalence $\Upsilon'$ ``externally marked'' on the
boundaries of all diamonds (but not on the boundaries of holes.) We
can then perform the complex pull-back on long monotone branches any
number of times without increasing the quasiconformal norm. 

Next, we need to show that the complete internal marking can be
realized with a bounded worsening of the bounds. This follows directly
from the proof of Lemma~\ref{lem:14ga,1}. This Lemma shows only how to
implement the marking condition at the critical value, but the
argument works the same way for any marking condition. 
\paragraph{The final maps.}
Now we can apply Lemma~\ref{lem:25ga,1} to these modified mappings, 
and get $\phi^{d}$, $\hat{\phi}^{d}$ and an externally marked
completely internally marked branchwise equivalence $\Upsilon^{d}$
between them. Now, they satisfy the hypotheses of Theorem 3. Proceed
by complex inducing with boundary refinement starting from $\phi^{d}$
and $\hat{\phi}^{d}$. 
It might be that full mappings occur infinitely many times in this
sequence. Otherwise, we can define {\em final} mappings
$\phi^{f}$, $\hat{\phi}^{f}$ with their branchwise
equivalence $\Upsilon^{f}$ as either the initial triple $\phi^{d}$,
etc., if no full mapping occurs in the sequence derived by complex
inducing, or the last triple in this sequence with $\phi^{f}$ and
$\hat{\phi}^{f}$ full. Observe that Theorem 3 is applicable with
$\Upsilon=\Upsilon_{b}=\Upsilon^{d}$. So, $\Upsilon^{f}$ has all
properties postulated by Theorem 3 for $\Upsilon_{1}$.       

\subsection{Proof of Theorem 5}
\paragraph{Getting rid of boundary refinement.}
\begin{lem}\label{lem:2hp,1}
Let $\phi$ be a box mapping, either full or of type I, infinitely
refined at the boundary, which undergoes
$k$ steps of box inducing. Suppose
that $\phi$ has standard $\epsilon$-extendibility. Then there is a
mapping $\phi^{r}$ obtained from $\phi$ by a finite number of
simultaneous monotone pull-back steps using $\phi':=\phi$ so that
after $k$ box inducing steps starting from $\phi^{r}$ the resulting
map has standard $\epsilon$-extendibility.
\end{lem}
\begin{proof}
The point is that box inducing skips boundary refinement. However, we
show that this can be offset by doing enough ``boundary refinement''
before entering the box construction. Observe first the following
thing. Under the hypotheses of the Lemma, suppose that after $k$ box
inducing steps we get a mapping $\phi_{k}$ and then do a simultaneous
monotone pull-back on all long monotone branches of $\phi_{k}$ using 
$\phi':=\phi_{k}$. Then the same mapping can be obtained by doing a
simultaneous monotone pull-back on all long monotone branches of the
original $\phi$. The proof of this remark proceeds by induction. For
$k=1$ this is rather obvious. For the induction step from $k-1$ to $k$
consider $\phi:=\phi_{1}$ and use the hypothesis of induction. It
follows that we need to perform inducing on all long branches of
$\phi_{1}$, and for that use the fact again with $k=1$. Now the lemma
follows immediately, since each time one needs boundary refinement in a
general inducing step, the appropriately refined mapping can be
obtained by simultaneous monotone pull-back on some or all long
monotone branches of $\phi$. 
\end{proof}    

\paragraph{Main estimate.}
\begin{lem}\label{lem:8ka,2}
Let $\phi$ and $\hat{\phi}$ be a pair of topologically conjugate
complex box mapping with diamonds of type I or full which undergo $k$
steps of complex box inducing resulting in mappings $\phi_{k}$ and
$\hat{\phi}_{k}$.  Suppose that both hole structures can be assigned
the separation index $K$. 
Suppose that box are infinitely refined at the boundary and have
standard $\epsilon$ extendibility. Also, suppose that a branchwise
equivalence $\Upsilon$ exists which is $Q$-quasiconformal, externally marked
and completely internally marked, and satisfies restricted to any long
monotone domain on the real line replaces $\Upsilon$ on the
fundamental inducing domain with distortion $K'$. Then there are
numbers $L_{1}$ which only depends on $K$ and $L_{2}$ depending on
$K'$ and $\epsilon$, with complex box mappings $\Phi$ and $\hat{\Phi}$
of type $I$ and a branchwise equivalence $\Upsilon'$ between them so
that the following conditions are satisfied:
\begin{em}
\begin{itemize}
\item
$\Upsilon'$ is $Q+L_{1}$-quasiconformal,
\item
$\Upsilon'$ is externally marked and completely internally marked,
\item
$\Upsilon'$ restricted to any long monotone domain replaces
$\Upsilon'$ on the fundamental inducing domain with distortion
$L_{2}$, 
\item
$\Phi$ has the same box structure as $\phi_{k}$, while $\hat{\Phi}$
has the same box structure as $\hat{\phi}_{k}$,
\item
$\Phi$ and $\hat{\Phi}$ both have standard $\epsilon$-extendibility.     
\end{itemize}
\end{em}
\end{lem}
\begin{proof}
The mapping $\Phi$ as obtained as $\phi^{r}$ for $\phi$ from
Lemma~\ref{lem:2hp,1}. $\hat{\Phi}$ is obtained in the same way for
$\hat{\phi}$. They are topologically conjugate. By
Lemma~\ref{lem:2hp,1} this means that we should obtain some
$\varphi_{0}$ by a series of simultaneous monotone pull-backs on long
monotone branches of $\phi$, and $\hat{\varphi}_{0}$ is obtained in an
analogous way for $\hat{\phi}$. Then we perform box inducing on
$\varphi_{0}$ , to get a sequence $\varphi_{i}$ with
$\varphi_{k}=\Phi$ and the same is done for $\hat{\varphi}_{0}$ which
gives $\hat{\Phi}=\hat{\varphi}_{k}$. The branchwise equivalence is
obtained by complex pull-back.   

Among the claims of Lemma~\ref{lem:8ka,2} the extendibility is clear
and the replacement condition follows in the usual way based on Lemma
4.6 of~\cite{kus}. The hard thing is the quasiconformal estimate for
$\Upsilon'$. The procedure used in the proof of Theorem 3 does not
give a uniform estimate for mappings which are not full. However, by
Proposition~\ref{prop:4hp,1} and Theorem 4
modifications required to obtain the marking in the box
case can be done with distortions which diminish exponentially fast at
a uniform rate.  
\end{proof}

\paragraph{Conclusion.}
For the proof of Theorem 5, we begin by Theorem 2 which tells us that
either we hit a suitable map first, or we can build induced mappings 
$\phi^{s}$, $\hat{\phi}^{s}$ and $\upsilon^{s}$. If we encounter the
suitable map first, then the conditions of Theorem 5 follow directly
from Theorem 2. Note that the assumption about the return time of the
restrictive interval into itself is needed to make sure that the
suitable mapping has monotone branches, and thus can be infinitely
refined at the boundary. 

Otherwise, we proceed to obtain final maps with the
branchwise equivalence between them. To this end, we build the complex
branchwise equivalence, by Lemma~\ref{lem:25ga,1}, and proceed by
Theorem 3 to obtain the branchwise equivalence between final maps.
If final maps do not exist, it means that infinitely many times in the
course of the construction we obtain full mappings. By Theorem 3, we
get them with uniformly quasiconformal branchwise equivalences. In
this sequence of full mappings the sizes of domains other than long
monotone ones go to $0$. So, having been given a $\delta$ we proceed
far enough in the construction, and then get the $\delta$-fine mapping
by applying simultaneous monotone pull-back on long monotone domains. 
Theorem 5 follows in this case as well.    

So we are only left with the case when final maps exist. Then we pick up
the construction by Lemma~\ref{lem:8ka,2}. Observe that the assumption
about the separation index is satisfied for the following reason. The
final map is either the same as $\phi^{s}$, in which case the bound
follows directly from Theorem 2, or is full and its holes are inside
the holes of $\phi^{s}$ constructed by Theorem 2, so the separation
index is even better. If $f$ was renormalizable, we choose $k$ in
Lemma~\ref{lem:8ka,2} equal to the number of box inducing steps needed
to get the suitable map. Then the conditions of Theorem 5 follow
directly from Lemma~\ref{lem:8ka,2}. The replacement condition on
short monotone domains is a consequence of the fact the by
construction the branchwise equivalence on short monotone domains is
the pull-back of the branchwise equivalence from $B$, and short
monotone branches are extendible by Theorem 4. When $f$ is
non-renormalizable we choose a large $k$ depending on $\delta$ and
follow up with a simultaneous monotone pull-back on all long domains. 
Theorem 5 likewise follows.    

\section{Proof of Theorem 1}
\subsection{The non-renormalizable case}
Theorem 1 in the non renormalizable case follows directly from Theorem
5. Choose a sequence $\delta_{n}$ tending to $0$. The corresponding
branchwise equivalences obtained by Theorem 5 for $\delta:\\delta_{n}$
will tend to the topological conjugacy in the $C^{0}$ norm. Since they
are all uniformly quasisymmetric in terms of $\eta$, so is the limit. 
   
\subsection{Construction of the saturated map}
In the renormalizable case, the only missing piece is the construction
of saturated maps with quasisymmetric branchwise equivalences between
them. Also, we need to make sure that the branches of the saturated
map are uniformly extendible. The case when the return time of the
maximal restrictive interval into itself is $2$ is not covered by
Theorem 5. In this case, we simply state that Theorem 1 is obvious and
proceed under the assumption that the return time is bigger than $2$. 
So, Theorem 1 follows from this proposition:

\begin{prop}\label{prop:8kp,1}
Suppose that conjugate suitable real box mappings,
$\varphi$ and $\hat{\varphi}$ are given, both either full or of type I  and
infinitely refined at the boundary, 
with a branchwise equivalence
$\Upsilon$ between them so that the following list of properties
holds:
\begin{em}
\begin{itemize}
\item
$\phi$ and $\hat{\phi}$ both have standard $\epsilon$-extendibility, 
\item
$\Upsilon$ is $Q$-quasisymmetric,
\item
$\Upsilon$ restricted to any long monotone domain replaces $\upsilon$
on the fundamental inducing domain with distortion $K$,
\item
$\Upsilon$ restricted to any short monotone domain replaces $\upsilon$
on $B$ with norm $K$.
\end{itemize}
\end{em}

Then, their saturated mappings $\varphi^{s}$ and $\hat{\varphi}^{s}$
are $\epsilon$-extendible. Also a $Q'$-quasisymmetric saturated
branchwise equivalence $\Upsilon'$ exists. $Q'$ depends on $Q$, $K$,
and $\epsilon$ only.  
\end{prop} 

The proof of Proposition~\ref{prop:8kp,1} is basically quoted
from~\cite{kus} with only minor adjustments.

\paragraph{An outline of the construction.}
Let $\psi$ mean the central branch. Let $I$ denote the restrictive
interval. First, we want to pull the branches $\varphi$ into
the domain of $\psi$. We notice that each point
of the line which is outside of the restrictive interval will be
mapped outside of the domain of $\psi$ under some number of
iterates of $\psi$. We can consider sets of points for which the 
number of iterates required to escape from the domain of $\psi$ is
fixed. Each such set clearly consists of two intervals symmetric with
respect to the critical point. The endpoints of these sets form two
symmetric sequences accumulating at the endpoints of the restrictive
interval, which will be called {\em outer
staircases}. Consequently, the connected components of these sets will
be called steps. 
 
This allows us to construct an induced map
from the complement of the restrictive interval in the domain of
$\psi$ to the outside of the domain $\psi$ with branches
defined on the steps of the outer staircases.  That means, we can
pull-back $\Upsilon$ to the inside of the domain of $\psi$. 

Next, we construct the {\em inner staircases}. We notice that every
point inside the restrictive interval but outside of the fundamental
inducing domain inside it is mapped into the fundamental inducing
domain eventually. Again, we can consider the sets on which the time
required to get to the fundamental inducing domain is fixed, and so we
get the steps of a pair of symmetric inner staircases.

So far, we have obtained an induced map which besides branches
inherited from $\varphi$ has uniformly extendible monotone branches
mapping onto $I$. Denote it with
$\varphi^{1}$. We now proceed by filling-in to get rid of short
monotone branches. 
We conclude with refinement of remaining long monotone branches.
Thus, we will be left with branches mappings onto $J$ only, so this is
a saturated map. Its extendibility follows from the standard argument
of inducing. The same inducing construction is used for
$\hat{\varphi}$. Because the extendibility is obvious, we only need to
worry the branchwise equivalence $\Upsilon'$.  

\paragraph{Outer staircases.}
Suppose that the domain of $\psi$
is very short compared with the length of the the domain of
$\varphi$. This means that the domain of $\psi$ is extremely
large compared with the restrictive interval. This unbounded situation
leads to certain difficulties and is dealt with in our next lemma.

\begin{lem}\label{lem:n5,1}
One can construct a branchwise equivalence $\Upsilon^{1}$ which is a 
pull-back of $\Upsilon$ with quasiconformal norm bounded as a uniform
function of the norm of $\Upsilon$. 
Furthermore, an
integer $i$ can be chosen so that the following conditions are satisfied:
\begin{itemize} 
\item
The functional equation
\[ \Upsilon\circ\psi^{j}=\hat{\psi}^{j}\Upsilon^{1} \; \]
holds for any $0\leq j\leq i$ whenever the left-hand side is defined.
\item
The length of the interval which consists of points whose $i$
consecutive images by $\psi$ remain in the domain of $\psi$ forms a
uniformly bounded ratio with the length of the restrictive interval. 
\end{itemize}   
\end{lem}
\begin{proof}
We rescale affinely so that the restrictive intervals become $[-1,1]$
in both maps. Denote the domains of $\psi$ and $\hat{\psi}$ with $P$
and $\hat{P}$ respectively.
Then, $\psi$ can be represented as $h(x^{2})$ where
$h''/h'$ is very small provided that $|P|$ is large. We can assume that $|P|$
is large, since otherwise we can take $\Upsilon^{1}:=\Upsilon$ to
satisfy the claim of our lemma.
We consider the round  disk $B$ (``box'') whose diameter  is the box of 
$\varphi$ ranged through by $\psi$. Because $\psi$ is extendible, and
its domain was assumed to be small compared to $P$, the preimage of
$B$ by $\phi$,
called $B_{1}$ sits inside $B$ with a large annulus between them. 
Analogous objects are constructed for $\hat{\varphi}$. It
is easy to build a quasiconformal extension $\upsilon$ of $\Upsilon$
which satisfies 
\[ \hat{\varphi}\circ\upsilon = \upsilon \circ \varphi\]
on $B_{1}$. With that, we are able to perform complex pull-back by
$\psi$ and $\hat{\psi}$.  

Also, assuming that $|P|$ is large enough, we can find a uniform $r$
so that the preimages
of  $B(0,r)$ by $\psi$, $\hat(\psi)$ and $z\rightarrow z^{2}$
are all inside $B(0,r/2)$. Also, we can have $B(0,r)$ contained in
$B_{1}$ as well as $\hat{B}_{1}$. Next, we choose the largest $i$ so that
$[-r,r]\subset\psi^{-i}(P)$ . 
Then, we change $\psi$ and $\hat{\psi}$. We will only
describe what is done to $\psi$. Outside of $B(0,r)$, $\psi$
is left unchanged. Inside 
the preimage of $B(0,r)$ by $z\rightarrow z^{2}$ it is $z\rightarrow
z^{2}$. In between, it
can be interpolated by a smooth degree $2$ cover with bounded
distortion.  The modified extension will be
denoted with $\psi'$.

Next, we pull-back $\Upsilon$ by $\psi'$ and $\hat{\psi}'$ exactly $i$
times. That is, if $\Upsilon_{0}$ is taken equal to $\Upsilon$, then
$\Upsilon_{j+1}$ is $\Upsilon$ refined by pulling-back $\Upsilon_{j}$
onto the domain of $\psi$. 
Now we need to check whether $\Upsilon^{1}$ has all the properties
claimed in the Lemma. To see the functional equation
condition, we note that all branches of any $\Upsilon_{j}$ are in the
region where $\psi$ coincides with
$\psi'$. 
The last condition easily follows  from the fact that $r$ can be chosen in a
uniform fashion. 
Also, the quasiconformal norm of $\Upsilon_{i}$ grows only by a constant compared
with $\Upsilon$, since points pass through the region of
non-conformality only once.  
\end{proof}

\paragraph{The staircase construction.}
We take $\Upsilon^{1}$ obtained in Lemma~\ref{lem:n5,1} and confine
our attention to its restriction to the real line, denoted with $\upsilon_{1}$.
We rely on the fact that $\upsilon_{1}$ is a quasisymmetric map and
its qs norm is uniformly bounded in terms of the quasiconformal norm
of $\Upsilon^{1}$.

\subparagraph{Completion of outer staircases.}
We will construct a induced maps $\varphi_{2}$ and $\hat{\varphi}_{2}$
with a branchwise equivalence $\upsilon_{2}$ with following properties:
\begin{itemize}
\item
The map $\upsilon_{2}$ coincides with $\upsilon_{1}$ outside of the domain
of $\psi$. Also, it satisfies 
\[ \upsilon_{2}\circ \psi^{j} = \hat{\psi}^{j}\circ\upsilon_{2} \]
on the complement of the restrictive interval provided that
$\psi^{j}$ is defined.
\item
Inside the restrictive interval, it is the ``inner staircase
equivalence'', that is, all endpoints of the inner staircase steps are
mapped onto the corresponding points.
\item
Its qs norm is uniformly bounded as a function of the qc norm of 
$\Upsilon$.
\end{itemize}

Outer staircases constructed in Lemma\ref{lem:n5,1}
connect the boundary points of the domain of $\psi$  to
the $i$-th steps which are in the close neighborhood of the
restrictive interval. Also, the $i$-th steps are the corresponding 
fundamental domains for the inverses of $\psi$ in the proximity of the
boundary of the restrictive interval. 

By bounded geometry of renormalization, see~\cite{miszczu}, 
the derivative of $\psi$ at the boundary of
the restrictive interval is uniformly bounded away from one.
Then, it is straightforward to see that the equivariant correspondence
between infinite outer staircases which uniquely extends
$\upsilon_{1}$ from the $i$-th steps is uniformly quasisymmetric (see
a more detailed argument in the last section of~cite{kus}. 

Inside the restrictive interval, the map is already determined on the
endpoints of steps, and can be extended in an equivariant way onto each
step of the inner staircase.

\paragraph{Rebuilding a complex map.}
Due to the irregular behavior of the branchwise equivalence in the
domains of branches of $\varphi_{2}$ that map onto the central domain,
they cannot be
filled-in by critical pull-backs used in~\cite{kus}. Instead, we will
construct an externally marked branchwise equivalence and apply
complex filling-in. The external marked will be achieved by
Lemma~\ref{lem:11na,2}. As the lip, we choose the circular arc which
intersects the real line at the endpoints of the central domain and
makes angles of $\pi/4$ with the line. The teeth will be the preimages
of the lip by short branches inside the central domain. We check that
the norm of this mouth is bounded. By the geodesic property
(see~\cite{miszczu}), the teeth are bounded by corresponding circular
arc of the same angle. The only property that needs to check is the
existence of a bounded modulus between the lip and any tooth.  
   This will follow if we prove that the intersection of a tooth with
the real has a definite neighborhood (in terms of the cross-ratio)
which is still inside the central domain. Since short branches inside
the central domain are
preimages of short monotone branches from the outside by a negative
Schwarzian map, it is enough to see the analogous property for the
domains of short branches in the box. If the number of box steps
leading to the suitable map was bounded, this follows from the
estimates of~\cite{yours}, as the long monotone branches adjacent to
the boundary of the box continue to have a uniformly large size.
Otherwise, one uses Theorem 3.  

Now we construct the branchwise equivalence on the lip. To this end,
we take the straight down projection from the lip to the central
domain, and lift the branchwise equivalence from the line. The
resulting map on the lip is quasisymmetric with
the norm bounded in terms of the quasisymmetric norm of $\varphi_{2}$.   
 Now we pull-back this map
to the teeth by dynamics. By the complex K\"{o}be lemma, the maps that
we use to pull-back are diffeomorphisms of bounded nonlinearity, so
they will preserve quasisymmetric properties. Since we assumed the
replacement condition for short monotone domains, we can fill each
tooth with a uniformly quasiconformal map.   
This leaves us in a
position to apply Lemma~\ref{lem:11na,2} to fill the mouth with a
uniformly quasiconformal branchwise equivalence. Call the mouth $W$.

Finally, we have to extend the branchwise equivalence to the whole
plane. To this end, we choose a half-disk with the fundamental
inducing domain normalized to $[-1,1]$ as the diameter, and make the
branchwise equivalence identity there. Next, we regard $W$ and all its
preimages by short monotone branches as teeth. This time, it is quite
clear that the norm is bounded. The branchwise equivalence of the
teeth is pulled back from $W$ by dynamics. The same argument as we
made in the preceding paragraph show that Lemma~\ref{lem:11na,2} can
be used to construct a complete branchwise equivalence on the plane.  
     
\paragraph{Construction of the saturated map.}
We now apply filling-in to  all branches which map onto the central
domain $W$.
On the level of inducing, the only
branches still left are those with the range equal to the fundamental
inducing domain of the renormalized map, and long monotone branches
onto the whole previous fundamental inducing domain.

\subparagraph{Final refinement.}
We end with a simultaneous monotone pull-back on all long monotone
branches.  The limit exists in  $L^{\infty}$ and is the saturated map.
The branchwise equivalence we get is quasiconformal as well. This can
be seen by choosing diamonds for all long monotone branches, marking
it externally, and using complex pull-back. The paper~\cite{kus}
offers an alternative way which does require external marking and
instead relies on a version of the Sewing Lemma (Fact~\ref{fa:14ga,1}
in our paper.) This closes the proof of Proposition~\ref{prop:8kp,1},
and therefore of all our remaining theorems. 

\vskip 1in
\goodbreak
\vskip -1in

\section*{Appendix}
\subsection{Estimates for hole structures}
\paragraph{Separation symbols for complex box mappings.}
\subparagraph{Definition of the symbols.}
Now, let $\varphi$ be a type I complex box mapping of rank $n$. 
An ordered quadruple
of real non-negative numbers:
 \[ s(B) := (s_{1}(B),\cdots,s_{4}(B))\; \]
will be said to give a separation symbol for $B$ if certain annuli
exist as described below. The annuli are either open or degenerate to
curves.  
Figure~\ref{fig:5xa,2}  shows a choice of separating annuli
for domain $B$, which is the same as domain $B$ from
Figure~\ref{fig:5xa,1}. 

We first assume that there are annuli $A_{1}(B)$ and $A_{2}(B)$.
Both annuli are contained in
$B_{n'}$. The annulus $A_{2}(B)$ surrounds $B_{n}$ separating it from
the domain of the analytic extension of $B$ with range $B_{n'}$.
Then $A_{1}(B)$ separates $A_{2}(B)$ from the boundary of
$B_{n'}$. We must have
\[ s_{2}(B) \leq \mod A_{2}(B)\; \mbox{ and}\]
\[ s_{1}(B) \leq \mod A_{2}(B) + \mod A_{1}(B)\; .\]

Next, three annuli are selected around $B$ which will give the
meaning of the two remaining components of the symbol. First, the
annulus $A'(B)$ is chosen exactly equal to the difference between the
domain of the canonical extension of
the branch defined on $B$ and the domain of $B$. Then, the
existence of $A_{3}(B)$ is postulated which surrounds $A'(B)$ separating it
from $B_{n}$ and from the boundary of $B_{n'}$.  Finally, $A_{4}(B)$ separates
$A_{3}(B)$ from the boundary of
$B_{n'}$. Then
\[ s_{3}(B) \leq \mod A'(B) + \mod A_{3}(B)\; \mbox{ and}\]
\[ s_{4}(B) \leq \mod A'(B) + \mod A_{3}(B) + \mod A_{4}(B)\; .\]

The dependence on $B$ will often be suppressed in our subsequent notations.
\paragraph{Normalized symbols.}
We will now arbitrarily impose certain algebraic relations
among various components of a separation symbol.
Choose a number $\beta$, and $\alpha:=\beta/2$, together
with $\lambda_{1}$ and $\lambda_{2}$. Assume
$\alpha\geq\lambda_{1},\lambda_{2}\geq -\alpha$ and $\lambda_{1} + 
\lambda_{2}\geq 0$.   
If these quantities are connected with a separation symbol $s(B)$ as
follows
\[ s_{1}(B) = \alpha + \lambda_{1}\; ,\]
\[ s_{2}(B) = \alpha - \lambda_{2}\; ,\]
\[ s_{3}(B) = \beta - \lambda_{1}\; ,\]
\[ s_{4}(B) = \beta + \lambda_{2}\; .\]

we will say that $s(B)$ is normalized with norm $\beta$ and
corrections $\lambda_{1}$ and $\lambda_{2}$.
\subparagraph{Separation index of a box mapping.}
For a type I complex box mapping $\phi$ a positive number $\beta$ is
called its {\em separation index} provided that valid normalized
separation symbols with norm $\beta$ exist for all univalent branches.   

\paragraph{Monotonicity of separation indexes.}
The nice property of separation indexes is that they do not decrease in
the box inducing process. In fact, one could show that they increase
at a uniform rate and this will be the final conclusion to be drawn
from Theorem 4. For now, we prove
\begin{lem}\label{lem:24gp,1}
Let $\phi_{i}$, $i$ between $0$ and $m$ be a sequence of complex box
mappings of type I with the property that the next one arises from the
previous one in a simple box inducing step. If $\beta_{0}$ is a
separation index of $\phi_{0}$, then $\beta_{0}$ is also a separation
index of $\phi_{i}$ for any $i<m$. In addition, if $\phi_{i}$ arises
after a non-close return, then $v(\phi_{i})\geq \beta_{0}/4$.  
\end{lem}
\begin{proof}
The proof of Lemma has to be split into a number of
cases. As analytic tools, we will use the behavior of moduli of annuli under
complex analytic mappings. Univalent maps transport the annuli
without a change of modulus, analytic branched covers of degree $2$
will at worst halve them, and for a sequence of nesting annuli their
moduli are super-additive (see~\cite{lehvi}, Ch. I, for proofs, or
~\cite{brahu} for an application to complex dynamics.) To facilitate
the discussion, we will also need a classification of branches
depending on how they arise in a simple box inducing step. 

\paragraph{Some terminology.}
Consider a abstract setting in which one has a bunch of univalent
branches with common range $B'$ and fills them in to get branches
mapping  onto some $B\subset B'$. 
The original branches mapping onto $B'$ will be called {\em parent
branches}  of the filling-in process. Clearly, every branch after the
filling-in has a dynamical extension with range $B'$. For two
branches, the domains of these respective extensions may 
be disjoint or contain one another. 
In the first case we say that the original branches were {\em independent}. 
Otherwise, the one mapped with a smaller extension domain is called 
{\em subordinate} to
the other one. Note that if $b'$ is subordinate to $b$, then the
extension of $b$ maps $b'$ onto another short univalent domain. 

We then distinguish the set of "maximal" branches subordinate to none. They 
are mapped by their parent branches directly onto the
central domain. Therefore, the domains of extensions of maximal
branches mapping onto $B'$ are
disjoint. They also cover domains of all branches. The extensions of maximal
branches are exactly {\em parent branches} of the filling-in process.
These extensions with range $B'$ will called {\em canonical extensions.}

Now, in a simple box inducing step, the parent branches are the short
monotone branches of $\tilde{\phi}$. Among these we distinguish at
most two {\em immediate branches} which restrictions of the central
branch of $\phi$ to the preimage of the central domain. All
non-immediate parent branches are compositions of the central branch
of $\phi$ with short monotone branches of $\phi$.  
For example, in the non-close return the first filling gives a set of
parent  branches, two of which may be immediate, which later get
filled in.  In the close return filling-in is done twice, so we will
be more  careful in speaking about parent branches.
Figure~\ref{fig:5xa,1} shows examples of independent and subordinate
domains. 

We will sometimes talk of branches meaning their domains, for example
saying that a branch is contained in its parent branch.  
We assume that $\phi$ has rank $n$, so $B=B_{n}$ and $B'=B_{n'}$. 
Let $\psi$ be the central branch of $\phi$. Let 
$B_{n+1}$ denote the central domain of the newly created map
$\phi^{1}$.  Observe
that $B_{(n+1)'}=B_{n}$. Suppose that $\beta$ is a separation index of
$\phi$. 
We will now proceed to build separation
symbols with norm $\beta$ for all short univalent domains of
$\phi^{1}$. Let $g$ be a short univalent branch of $\phi^{1}$ and 
$p$ denote the parent branch $g$. The parent branch necessarily has
the form $P'\circ\psi$. Let $P$ be the branch of $\phi$ whose domain
contains the critical value.  Objects (separation annuli, components
of separation symbols) referring to $\phi^{1}$ will be marked with
bars. 

\subparagraph{Reduction to maximal branches.}
Note that it is sufficient show that symbols with norm $\beta$ exist
for maximal branches. Indeed, suppose that a separation symbol exists 
for a maximal branch $b$ and let $b'$ be subordinate to $b$. 
We can take $A_{1}(b') = A_{1}(b)$ and
$A_{2}(b')=A_{2}(b)$. Likewise, we can certainly adopt 
$A_{4}(b')=A_{4}(b)$, and $A_{3}(b')$ can be chosen to contain
$A_{3}(b)$.  The annulus $A'(b)$ is the
preimage of the annulus $B_{n'}\setminus B_{n}$ by the parent branch
of $b$. The annulus $A'(b')$ is the preimage of the same annulus by
the canonical extension of $b$, so it has the same modulus.
Since the domain of the canonical extension
of $b'$ is contained in the parent domain  (equal to the domain of the
canonical extension of $b$), the assertion follows. 

\paragraph{Non-close returns.}
Let us assume that $\phi$ makes a non-close return, that is $P'\neq
\psi$. 
\subparagraph{The case of $p$ immediate.}
Let $b$ denote the maximal branch in $p$. 
The new central hole $B_{n+1}$ is separated from the boundary of
$B_{n}$ by an annulus of modulus at least $(\beta+\lambda_{2}(B))/2$.
The annulus
$\overline{A}_{2}(b)$ around $B_{n+1}$ will be the preimage by the
central branch of the region contained in and between $A_{3}(P)$ and
$A'(P)$. Then, $\overline{A}_{1}(b)$ is the preimage of $A_{4}$. It
follows that we can take
\[ \overline{s}_{1} = \frac{\beta+\lambda_{2}(P)}{2}\; \mbox{and}\]     
\[ \overline{s}_{2} = \frac{\beta-\lambda_{1}(P)}{2} \; .\]
Of course, since components of the symbol are only lower estimates, we
are always allowed to decrease them if needed.
The annulus $\overline{A}'$ is naturally given as the preimage of the
annulus between $B_{n+1}$ and the boundary of $B_{n}$ by the central
branch, likewise $\overline{A}_{3}$ is the preimage of $A_{2}(P)$, and
$\overline{A}_{4}$ is the preimage of $A_{1}(P)$. Since the first two
preimages are taken in an univalent fashion, we get
\[ \overline{s}_{3} = \frac{\beta+\lambda_{2}(B)}{2} + \alpha -
\lambda_{2}(B)  \; \mbox{and}\]
\[ \overline{s}_{4} = \overline{s}_{3} +
\frac{\lambda_{1}(B)+\lambda_{2}(B)}{2} = \frac{\beta}{2} + \alpha +
\frac{\lambda_{1}(B)}{2}\; .\]
Thus, if we put 
\[ \overline{\lambda}_{1} = \frac{\lambda_{2}(B)}{2}\: ,\: 
   \overline{\lambda}_{2} = \frac{\lambda_{1}(B)}{2} \]
we get a valid separation symbol with norm $\beta$.
In the remaining non-immediate cases, the branch $P'$ is defined. 

\subparagraph{$P'$ and $P$ non-immediate and independent.}
To pick $\overline{A}_{2}(b)$, we take the preimage by $\psi$ of
the annulus separating $P$ from the boundary of the domain of
its canonical extension with range $B_{n'}$, i.e. $A'(P)$. 
We claim that its modulus in all 
cases is
estimated from below by
$\alpha+\delta$ where $\delta$ is chosen as the supremum of
$-\lambda_{2}(b')$ over all univalent domains $b'$ of $\phi$.  
Indeed, $P$ is carried 
onto $B_{n}$ by the extended branch, and the estimate is $\alpha$ plus the
maximum of $\lambda_{1}(b')$ with $b'$ ranging over the set of all
short univalent domains of $\phi$.  
and $\lambda_{1}(P')$ which is at least
The assertion follows  since $\lambda_{1}(b') + \lambda_{2}(b') \geq
0$ for any $b'$. 
To pick $A_{1}(b)$, consider the
annulus separating $A'(P)$ from the
boundary of $B_{n'}$, i.e. the region in and between $A_{3}(P)$ and 
$A_{4}(P)$. Pull this region 
back by the central branch to get 
$\overline{A}_{1}(b)$.
By the hypothesis of the induction, the estimates are
\[ \overline{s}_{1} = \frac{\beta+\lambda_{2}(P)}{2}\; \mbox{and}\]     
\[ \overline{s}_{2} = \frac{\alpha+\delta}{2} \; .\]

Since $b$ is maximal $\overline{A}'(b)$ is determined with modulus at least
$\overline{s}_{1}$. The annulus $\overline{A}_{3}(b)$ will be obtained as the
preimage by the central branch of $A'(P')$. This has modulus at
least $\alpha+\delta$ in all cases as argued above. The annulus 
$\overline{A}_{4}(b)$ is the preimage of the region in and between
$A_{3}(P')$ and $A_{4}(P')$. 
By induction, 
\[ \overline{s}_{3} = \frac{\beta+\lambda_{2}(P)}{2} +
\alpha+\delta\; \mbox{and}\]
\[ \overline{s}_{4} = \overline{s}_{3} +
\frac{\beta + \lambda_{2}(P') - \alpha - \delta}{2} \; .\]

We put $\overline{\lambda}_{1} = \frac{\lambda_{2}(P)}{2}$ and
$\overline{\lambda}_{2}=\frac{\alpha-\delta}{2}$. We check that 
\[ \overline{s}_{3} + \overline{\lambda_{1}} = \frac{\beta}{2} + \alpha +
\lambda_{2}(P) + \delta \geq \beta - \lambda_{2}(P) +
\lambda_{2}(P) \geq \beta \; .\]

In a similar way one verifies that
\[ \overline{s}_{4} - \overline{\lambda}_{2} \geq \beta\; .\]  
Also, the required inequalities between corrections
$\overline{\lambda}_{i}$ follow directly. 

 \subparagraph{$P'$ subordinate to $P$.}
This means that some univalent mapping onto $B_{n'}$ transforms $P$ onto
$B_{n}$ and $P'$ onto some $P''$. Consider $A_{2}(P'')$ which separates
$B_{n}$ from $P''$, and a larger annulus $A_{1}(P'')$. Their preimages
first  by the extended branch and then by
the central branch give us $\overline{A}_{2}(b)$ and 
$\overline{A}_{1}(b)$ respectively.     
The estimates are
\[ \overline{s}_{2} = \frac{\alpha-\lambda_{2}(P'')}{2} \;\mbox{and} \]
\[ \overline{s}_{1} = \frac{\alpha+\lambda_{1}(P'')}{2}\; .\]     

The annulus $\overline{A}'(b)$ is uniquely determined with modulus
$\overline{s}_{1}$, and $\overline{A}_{3}(b)$ will be the preimage of the
annulus separating $P''$ from $B_{n}$. Finally, $\overline{A}_{4}(b)$
will separate the image of $\overline{A}_{3}(b)$ from $B_{n'}$. 
The estimates are
\[ \overline{s}_{3} = \frac{\alpha+\lambda_{1}(P'')}{2} + \beta -
\lambda_{1}(P'') = \beta + \frac{\alpha-\lambda_{1}(P'')}{2} \; \mbox{and}\]
\[ \overline{s}_{4} = \overline{s}_{3} +
\frac{\lambda_{1}(P'')+\lambda_{2}(P'')}{2} = \beta +
\frac{\alpha+\lambda_{2}(P'')}{2} \; .\]

Set 
\[ \overline{\lambda}_{1} = \frac{-\alpha+\lambda_{1}(P'')}{2}\; \mbox{ and}\] 
\[ \overline{\lambda}_{2} = \frac{\alpha+\lambda_{2}(P'')}{2}\; .\]
The requirements of a normalized symbol are clearly
satisfied. 

\subparagraph{$P$ subordinate to $P'$.}
This situation is
analogous to the situation of immediate parent branch considered at the
beginning. Indeed, by mapping $P'$ to $B_{n}$ and composing with the
central branch one can get a folding branch with range $B_{n'}$ defined on
$P'$. We now see that the situation inside the domain of the canonical
extension of $P$ is analogous to the case
of immediate parent branches, except that the folding branch maps onto
a larger set
$B_{n'}$. So the estimates can only improve. 

\paragraph{A close return.}
In this case there are no immediate parent branches and we really have
only one case to consider. Fix some short univalent branch $b$ of
$\phi^{1}$, let $p$ be its parent branch, and denote $p=P'\circ\psi$. 
Consider $A_{2}(P')$ and $A_{1}(P')$. Their preimages by
the central branch give us $\overline{A}_{2}(b)$ and 
$\overline{A}_{1}(b)$ respectively.     
The estimates are
\[ \overline{s}_{2} = \frac{\alpha-\lambda_{2}(P')}{2} \;\mbox{and} \]
\[ \overline{s}_{1} = \frac{\alpha+\lambda_{1}(P')}{2}\; .\]     

The annulus $\overline{A}'(b)$ is uniquely determined with modulus
$\overline{s}_{1}$, and $\overline{A}_{3}(b)$ will be the preimage of the
annulus separating $P'$ from $B_{n}$ i.e. the annulus containing 
$A_{3}(P')$ and $A'(P')$ together with the region between them.
Finally, $\overline{A}_{4}(b)$
will be the preimage of $A_{4}(P')$ by $\psi$.
The estimates are
\[ \overline{s}_{3} = \frac{\alpha+\lambda_{1}(P')}{2} + \beta -
\lambda_{1}(P') = \beta + \frac{\alpha-\lambda_{1}(P')}{2} \; \mbox{and}\]
\[ \overline{s}_{4} = \overline{s}_{3} +
\frac{\lambda_{1}(P')+\lambda_{2}(P')}{2} = \beta +
\frac{\alpha+\lambda_{2}(P')}{2} \; .\]

Set 
\[ \overline{\lambda}_{1} = \frac{-\alpha+\lambda_{1}(P')}{2}\; \mbox{ and}\] 
\[ \overline{\lambda}_{2} = \frac{\alpha+\lambda_{2}(P')}{2}\; .\]

The requirements of a normalized symbol are clearly
satisfied. Not quite surprisingly, these are the same estimates we got
in the non-close case with $P'$ subordinate to $P$. 

\paragraph{Conclusion.}
We already proved by induction that $\beta_{0}$ remains a separation index
for all $\phi_{i}$. It remains to obtained the estimate 
$v(\phi_{i+1})\geq \beta_{0}/4$ under the assumption that $\phi_{i}$
makes a non-close return. This is quite obvious from considering the
separation symbol for the branch $P$ which contains the critical
value. Since $s_{4}(P) = \beta_{0}+\lambda_{2}(P) \ geq \beta_{0}/2$
and because of superadditivity of conformal moduli, there is an
annulus with modulus at least $\beta_{0}/2$ separating $P$ from $B'$,
and its pull-back by $\psi$ gives as an annulus with desired modulus.
This all we need to finish the proof of Lemma~\ref{lem:24gp,1}. Note,
however, that we cannot automatically claim $\overline{s}_{1} \geq 
\beta_{0}/4$  even if the preceding return was non-close. 
\end{proof}

\clearpage
\pagestyle{empty}
\begin{figure}[pt]
\psfig{figure=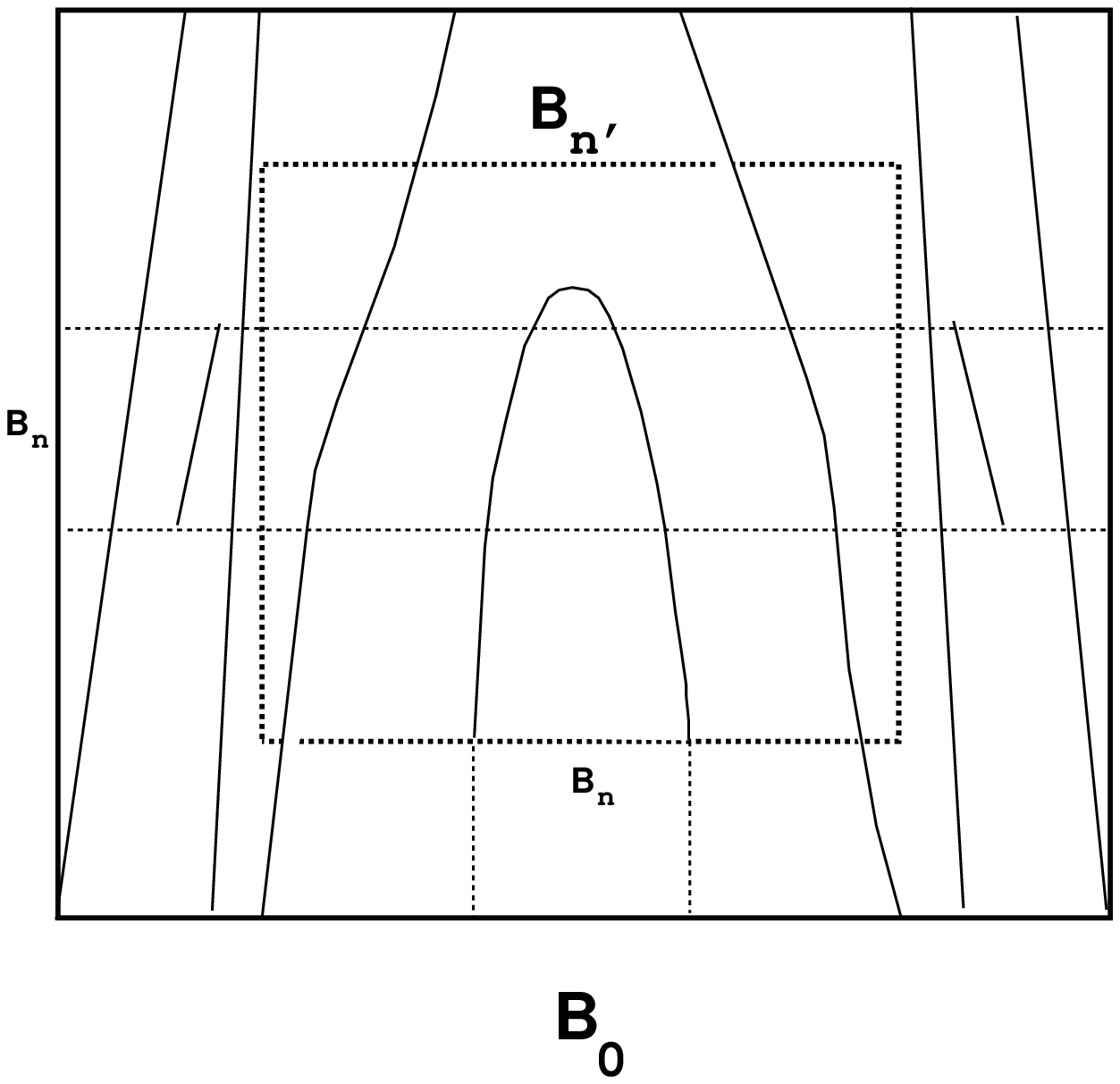}
\caption{A sample graph of a type I real box mapping. All three
permissible types of branches: central folding, long and short
monotone are shown. Be aware that typically one has infinitely many
branches.}
\label{fig:4xp,1}
\end{figure}

\begin{figure}[pt]
\psfig{figure=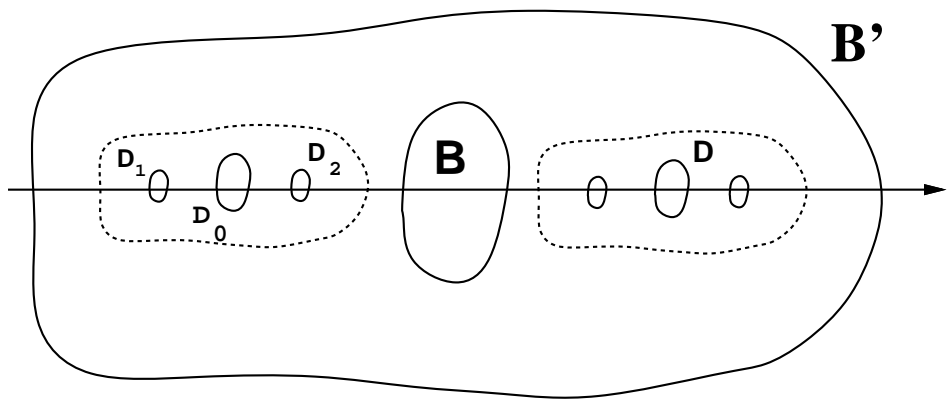}
\caption{A type I complex box mapping. Dotted lines show domains of
canonical extensions. Domains $D_{0}$ and $D$ are look like they are maximal. 
Then $D_{1}$ and $D_{2}$ are subordinate to $D_{0}$, but apparently
independent from one another as well as from $D$. $D_{0}$ and $D$ are also
independent. There may be
univalent domains outside of $B'$, not shown here.}
\label{fig:5xa,1}
\end{figure}

\begin{figure}[pt]
\psfig{figure=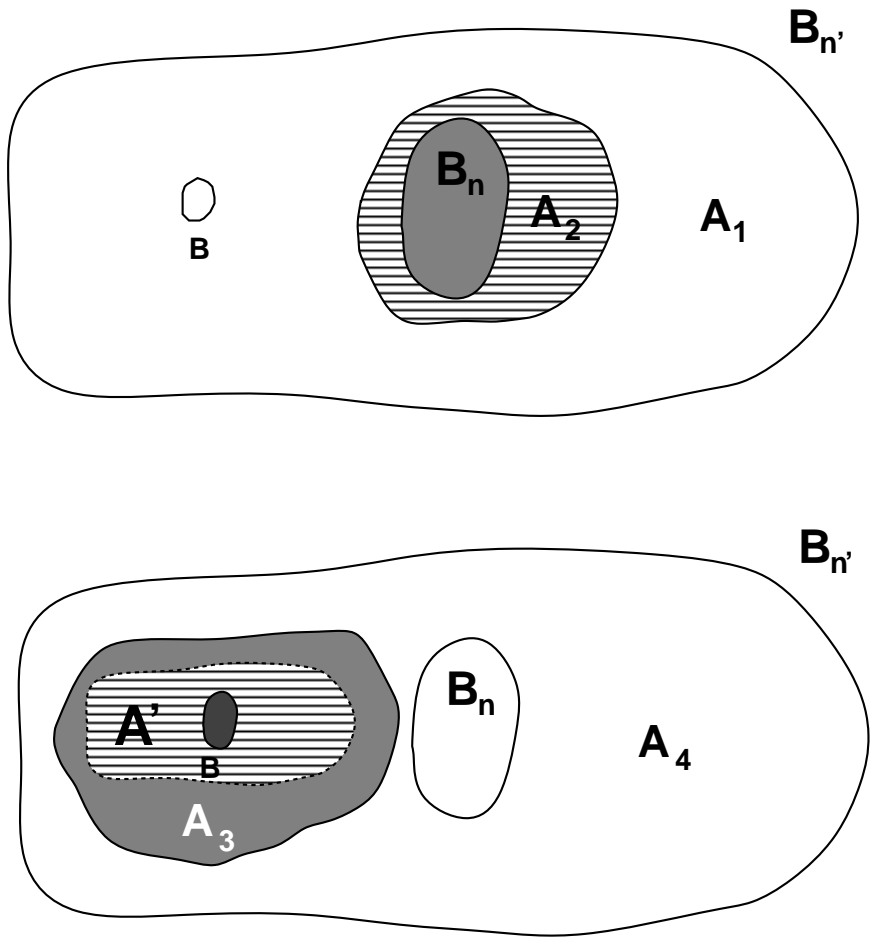}
\caption{A choice of separating annuli for $B$. Note that the
outermost annuli $A_{1}$ and $A_{4}$ are filled in white.}
\label{fig:5xa,2}
\end{figure}
\clearpage

\end{document}

%% file: imsmark.tex
\def\IMSmarkvadjust{0 pt}
\def\IMSmarkhadjust{0 pt}
\def\IMSmarkhpadding{0 pt}
\def\IMSpubltext{Published in modified form:}
\def\SBIMSMark#1#2#3{
 \font\SBF=cmss10 at 10 true pt
 \font\SBI=cmssi10 at 10 true pt
 \setbox0=\hbox{\SBF \hbox to \IMSmarkhpadding{\relax}
                Stony Brook IMS Preprint \##1}
 \setbox2=\hbox to \wd0{\hfil \SBI #2}
 \setbox4=\hbox to \wd0{\hfil \SBI #3}
 \setbox6=\hbox to \wd0{\hss
             \vbox{\hsize=\wd0 \parskip=0pt \baselineskip=10 true pt
                   \copy0 \break%
                   \copy2 \break%
                   \copy4 \break}}
 \dimen0=\ht6   \advance\dimen0 by \vsize \advance\dimen0 by 8 true pt
                \advance\dimen0 by -\pagetotal
	        \advance\dimen0 by \IMSmarkvadjust
 \dimen2=\hsize \advance\dimen2 by .25 true in
	        \advance\dimen2 by \IMSmarkhadjust

%
%
  \openin2=publishd.tex
  \ifeof2\setbox0=\hbox to 0pt{}
  \else 
     \setbox0=\hbox to 3.1 true in{
                \vbox to \ht6{\hsize=3 true in \parskip=0pt  \noindent  
                {\SBI \IMSpubltext}\hfil\break
                \input publishd.tex 
                \vfill}}
  \fi
  \closein2
  \ht0=0pt \dp0=0pt
 \ht6=0pt \dp6=0pt
 \setbox8=\vbox to \dimen0{\vfill \hbox to \dimen2{\copy0 \hss \copy6}}
 \ht8=0pt \dp8=0pt \wd8=0pt
 \copy8
 \message{*** Stony Brook IMS Preprint #1, #2. #3 ***}
}